\documentclass[10pt]{article}%

\usepackage{mathrsfs}
\DeclareMathAlphabet{\mathpzc}{OT1}{pzc}{m}{it}

\usepackage{color}
\usepackage{amsmath}
\usepackage{graphicx}
\usepackage{amsfonts}
\usepackage{amssymb}%
\usepackage{latexsym}
\usepackage{psfrag}
\usepackage{subfig}
\usepackage{tikz}
\usetikzlibrary{patterns}

\newtheorem{theorem}{Theorem}[section]
\newtheorem{corollary}[theorem]{Corollary}
\newtheorem{definition}[theorem]{Definition}
\newenvironment{proof}[1][Proof]{\noindent \emph{#1.} }
{\hfill \ \rule{0.5em}{0.5em}}
\newtheorem{lemma}[theorem]{Lemma}
\newtheorem{proposition}[theorem]{Proposition}

\newtheorem{example}[theorem]{Example}
\newtheorem{remark}[theorem]{Remark}

\newtheorem{condition}[theorem]{Condition}

\newcounter{mycount}

\numberwithin{equation}{section}

\usepackage[a4paper]{geometry}
\newcommand{\beq}{\begin{equation}}
\newcommand{\eeq}{\end{equation}}
\newcommand{\beqs}{\begin{equation*}}
\newcommand{\eeqs}{\end{equation*}}
\newcommand{\bit}{\begin{itemize}}
\newcommand{\eit}{\end{itemize}}
\newcommand{\ben}{\begin{enumerate}}
\newcommand{\een}{\end{enumerate}}
\newcommand{\bal}{\begin{align}}
\newcommand{\eal}{\end{align}}
\newcommand{\bals}{\begin{align*}}
\newcommand{\eals}{\end{align*}}
\newcommand{\bse}{\begin{subequations}}
\newcommand{\ese}{\end{subequations}}
\newcommand{\bpr}{\begin{proposition}}
\newcommand{\epr}{\end{proposition}}
\newcommand{\bre}{\begin{remark}}
\newcommand{\ere}{\end{remark}}
\newcommand{\bpf}{\begin{proof}}
\newcommand{\epf}{\end{proof}}
\newcommand{\ble}{\begin{lemma}}
\newcommand{\ele}{\end{lemma}}
\newcommand{\bco}{\begin{corollary}}
\newcommand{\eco}{\end{corollary}}
\newcommand{\bex}{\begin{example}}
\newcommand{\eex}{\end{example}}
\newcommand{\bth}{\begin{theorem}}
\newcommand{\enth}{\end{theorem}}

\newcommand{\noi}{\noindent}

\newcommand{\bT}{\mathbf{T}}

\newcommand{\cC}{{\cal C}}

\newcommand{\cL}{{\cal L}}
\newcommand{\opL}{\cL}

\newcommand{\cD}{{\cal D}}

\newcommand{\cM}{{\cal M}}

\newcommand{\bx}{\mathbf{x}}

\newcommand{\bxi}{{\boldsymbol{\xi}}}
\newcommand{\balpha}{{\boldsymbol{\alpha}}}

\newcommand{\bze}{\mathbf{0}}

\newcommand{\bnu}{\boldsymbol{\nu}}

\newcommand{\by}{\mathbf{y}}

\newcommand{\bQ}{\mathbf{Q}}
\newcommand{\bF}{\mathbf{F}}
\newcommand{\supp}{\mathrm{supp}}

\newcommand{\eps}{\varepsilon}

\newcommand{\ri}{{\rm i}}
\newcommand{\rd}{{\rm d}}

\newcommand{\Rea}{\mathbb{R}}
\newcommand{\Com}{\mathbb{C}}

\newcommand{\half}{\frac{1}{2}}

\newcommand{\esssup}{\mathop{{\rm ess} \sup}}
\newcommand{\essinf}{\mathop{{\rm ess} \inf}}

\newcommand{\loc}{\mathop{{\rm loc}}}

\newcommand{\dist}{\mathop{{\rm dist}}}

\newcommand{\diff}[2]{\frac{\rd #1}{\rd #2}}
\newcommand{\pdiff}[2]{\frac{\partial #1}{\partial #2}}

\newcommand{\Oi}{{\Omega_-}}
\newcommand{\Oe}{{\Omega_+}}
\newcommand{\OR}{{\Omega_R}}

\newcommand{\GR}{{\Gamma_R}}

\newcommand{\pD}{{\partial \domaingen}}

\newcommand{\GammaR}{\Gamma_R}

\newcommand{\dudnu}{\pdiff{u}{\nu}}
\newcommand{\dudnuw}{\partial u/\partial \nu}
\newcommand{\dudnuA}{\pdiff{u}{\nu_A}}
\newcommand{\dudnuAw}{\partial u/\partial \nu_A}

\newcommand{\dvdnu}{\pdiff{v}{\nu}}
\newcommand{\dvdnuw}{\partial v/\partial \nu}
\newcommand{\dvdnuA}{\pdiff{v}{\nu_A}}
\newcommand{\dvdnuAw}{\partial v/\partial \nu_A}

\newcommand{\nT}{\nabla_{\pD}}

\newcommand{\nus}{|u|^2}
\newcommand{\mus}{\nus}
\newcommand{\ngus}{|\nabla u|^2}
\newcommand{\nurs}{|u_r|^2}

\newcommand{\gu}{\nabla u}

\newcommand{\nvs}{|v|^2}
\newcommand{\ngvs}{|\nabla v|^2}

\newcommand{\nvrs}{|v_r|^2}

\newcommand{\gv}{\nabla v}
\newcommand{\gw}{\nabla w}
\newcommand{\gvb}{\overline{\nabla v}}

\newcommand{\vb}{\overline{v}}
\newcommand{\ub}{\overline{u}}

\newcommand{\LtGD}{{L^2(\Gamma_D)}}
\newcommand{\LtGI}{{L^2(\Gamma_I)}}

\newcommand{\LtD}{{L^2(\domain)}}

\newcommand{\HoDk}{{H^1_{0,D}(\domain_R)}}
\newcommand{\HoDkk}{{H^1_{k}(\domain_R)}}

\newcommand{\DOmegabar}{\cD(\overline{\domaingen})}

\newcommand{\tendi}{\rightarrow \infty}
\newcommand{\tendo}{\rightarrow 0}

\newcommand{\vecer}{\mathbf{e}_r}

\newcommand{\nmin}{\varmin}

\def\XXint#1#2#3{{\setbox0=\hbox{$#1{#2#3}{\int}$}
     \vcenter{\hbox{$#2#3$}}\kern-.5\wd0}}

\newcommand*{\N}[1]{\left\|#1\right\|}

\newcommand{\vertiii}[1]{{\left\vert\kern-0.25ex\left\vert\kern-0.25ex\left\vert #1 
    \right\vert\kern-0.25ex\right\vert\kern-0.25ex\right\vert}}

\usepackage{hyperref}
\definecolor{myblue}{rgb}{0,0,0.6}
\hypersetup{colorlinks=true,
linkcolor=myblue,citecolor=myblue,filecolor=myblue,urlcolor=myblue}

\allowdisplaybreaks[4]

\newcommand{\ton}{\text{ on }}
\newcommand{\tin}{\text{ in }}
\newcommand{\tfa}{\text{ for all }}
\newcommand{\tfor}{\text{ for }}
\newcommand{\tas}{\text{ as }}
\newcommand{\tand}{\text{ and }}
\newcommand{\tst}{\text{ such that }}

\newcommand{\hatx}{\widehat{\bx}}

\newcommand{\var}{{n}}
\newcommand{\varmin}{{n_{\min}}}
\newcommand{\varmax}{{n_{\max}}}

\newcommand{\princP}{\sigma(P)}

\newcommand{\domain}{\Omega}
\newcommand{\domaingen}{\Omega}

\newcommand{\Oin}{\Omega_T}
\newcommand{\Oout}{\Rea^d\setminus\overline{\Omega_T}}
\definecolor{escol}{rgb}{0,0,0.8}
\definecolor{estcol}{rgb}{0,0.8,0}

\begin{document}

\title{The Helmholtz equation in heterogeneous media: a priori bounds, well-posedness, and resonances
}

\author{
I. G.~Graham\footnotemark[1],\,\, 
O. R.~Pembery\footnotemark[1],\,\, 
E. A.~Spence\footnotemark[1]
}

\date{\today}

\renewcommand{\thefootnote}{\fnsymbol{footnote}}

\footnotetext[1]{Department of Mathematical Sciences, University of Bath, Bath, BA2 7AY, UK, \tt I.G.Graham@bath.ac.uk, O.R.Pembery@bath.ac.uk, E.A.Spence@bath.ac.uk}

\renewcommand{\thefootnote}{\arabic{footnote}}

\maketitle

\begin{abstract}
We consider the exterior Dirichlet problem for the heterogeneous Helmholtz equation, i.e.~the equation $\nabla\cdot(A \gu ) + k^2 n u =-f$  where both $A$ and $n$ are functions of position. We prove new a priori bounds on the solution under conditions on $A$, $n$, and the domain that ensure nontrapping of rays; 
the novelty is that these bounds are explicit in $k$, $A$, $n$, and geometric parameters of the domain. We then show that these a priori bounds hold when $A$ and $n$ are $L^\infty$ and satisfy certain monotonicity conditions, and thereby obtain new results \emph{both} about the well-posedness of such problems \emph{and} about the resonances of acoustic transmission problems (i.e.~$A$ and $n$ discontinuous) where the transmission interfaces are only assumed to be $C^0$ and star-shaped; the novelty of this latter result is that until recently the only known results about resonances of acoustic transmission problems were for $C^\infty$ convex interfaces with strictly positive curvature.

\paragraph{Keywords:} Helmholtz equation, heterogeneous, variable wave speed, high frequency, transmission problem, nontrapping, resolvent, uniqueness, resonance, semiclassical

\paragraph{AMS subject classifications:} 35J05, 35J25, 35B34, 35P25, 78A45

\end{abstract}

\section{Introduction}\label{sec:intro}

This paper is concerned with proving a priori bounds on the solutions of the heterogeneous Helmholtz equation 
\beq\label{eq:1}
\opL_{A,n} u:= \nabla\cdot(A \gu ) + k^2 n u =-f,
\eeq
where $k>0$ is the wavenumber, $u=u(\bx)$ with $\bx\in \Rea^d$, $d=2,3$, and where
\ben
\item $A$ is a symmetric, real-valued, positive-definite  matrix function of $\bx$, 
\item $n$ is a real-valued function of $\bx$, bounded away from zero, and 
\item both $I-A$ and $1-n$ have compact support.
\een
The PDE \eqref{eq:1} arises by taking the Fourier transform in time (with Fourier variable $k$) of the heterogeneous wave equation
\beq\label{eq:2}
 \nabla\cdot(A \nabla U ) - \frac{1}{c^2} \pdiff{^2 U}{t^2} =-F,
\eeq
for $U=U(\bx,t)$, where $n(\bx)= c(\bx)^{-2}$. Important applications of \eqref{eq:1}/\eqref{eq:2} include describing the 
transverse-magnetic (TM) and transverse-electric (TE) modes of Maxwell's equations (see, e.g., \cite[Remark 2.5]{MoSp:17}) and forming the so-called 
``acoustic approximation" of the elastodynamic wave equation (see, e.g., the derivation in \cite[\S1.2]{Ch:15} and the references therein).

We consider the \emph{exterior Dirichlet problem (EDP)}; i.e.~Equation \eqref{eq:1} is posed in the exterior of a bounded obstacle, with Dirichlet boundary conditions on the obstacle, the Sommerfeld radiation condition at infinity, and data $f$ with compact support. Observe that the assumptions that $I-A$ and $1-n$ 
have compact support mean that operator $\cL_{A,n}$ becomes the homogeneous Helmholtz operator $\cL := \Delta  +k^2$ outside a compact set.
In the case when the obstacle is the empty set, the EDP becomes the \emph{full-space problem}.

We allow one or both of $A$ and $n$ to be discontinuous, and such BVPs are usually called \emph{transmission problems}; 
e.g., when the obstacle is the empty set, and $A$ and $n$ are discontinuous on a common interface, the EDP becomes the classic problem of transmission through one penetrable obstacle.

A standard model problem in the numerical analysis of the Helmholtz equation is the \emph{truncated exterior Dirichlet problem (TEDP)}, where 
the radiation condition is approximated by truncating the (unbounded) exterior domain and applying an impedance boundary condition on the artificial boundary (see, e.g., \cite[\S3.2]{Ih:98}, \cite[\S5.1]{BaSpWu:16}, and the references therein for discussion of this approximation). 
When the obstacle is the empty set, the TEDP becomes the \emph{interior impedance problem (IIP)}. Although we focus on the EDP, we show in Appendix \ref{sec:TEDP} how our results apply to the TEDP.

\paragraph{Goal and motivation.}
Our goal is to obtain a priori bounds on the solution to the EDP that are explicit in $k>0$, $A$, and $n$. 
These a priori bounds then allow us to prove existence and uniqueness results about the solution of the EDP via Fredholm theory, and also to prove results about the location of resonances in the complex $k$-plane (via the classical link between $k$-explicit a priori bounds for $k$ real and resonance-free regions under the real-$k$ axis; see, e.g., \cite{Vo:99}). 

Obtaining $k$-explicit a priori bounds on the solutions of \eqref{eq:1} is a classic problem \cite{Bl:73, Va:75, BlKa:77, Bu:98, PeVe:99, PoVo:99a, PoVo:99, CaPoVo:99, Bu:02, Be:03a, Ca:12, CaLePa:12, NgVo:12, Sh:17} 
which has recently been the focus of renewed interest from the numerical-analysis community (with this interest mostly focused on the TEDP/IIP). 
Indeed, 
there has been sustained interest in proving a priori bounds on either the TEDP or the IIP when $A\equiv I$ and $n\equiv 1$ \cite{Me:95, CuFe:06, He:07, EsMe:12, Sp:14, BaSpWu:16} and recent interest in proving bounds when one or both of $A$ and $n$ are variable \cite{BrGaPe:15, Ch:16, BaChGo:17, OhVe:18, MoSp:17, SaTo:17, GrSa:18}.

One of our main motivations for proving bounds that are not only explicit in $k$, but also explicit in $A$ and $n$ is that such bounds can then be used to prove a priori bounds and results about the well-posedness of \eqref{eq:1} when 
$A$ and $n$ are random fields; 
this is considered in the companion paper \cite{PeSp:18}.

\paragraph{Recap of existing well-posedness results.}
In 2-d, the unique continuation principle (UCP) holds (and gives uniqueness) when $A$ is $L^\infty$ and $n\in L^p$ for some $p>1$ \cite{Al:12}. In 3-d, the UCP holds when $A$ is Lipschitz \cite{GaLi:87, Ka:88} and $n \in L^{3/2}$ \cite{JeKe:85,  Wo:92}; see \cite{GrSa:18} for these results applied specifically to Helmholtz problems.
Fredholm theory then gives existence and an a priori bound on the solution; this bound, however, is not explicit in $k$, $A$, or $n$. 

An example of an $A \in C^{0,\alpha}$ for all $\alpha<1$ for which the UCP fails in 3-d  is given in 
\cite{Fi:01}. Nevertheless, the UCP can be extended from Lipschitz $A$ to piecewise-Lipschitz $A$ 
by the Baire-category argument in \cite{BaCaTs:12} (see also \cite[Proposition 2.11]{LiRoXi:16}), with well-posedness then following by Fredholm theory as before -- we discuss this argument of  \cite{BaCaTs:12} further in \S\ref{sec:previous}.

\paragraph{Recap of existing a priori bounds on the EDP in trapping and nontrapping situations.} 
In this overview discussion, for simplicity, we consider 
the case of zero Dirichlet boundary conditions on $\partial\domain_-$, where $\domain_-$ denotes the obstacle.

When $A, n$, and $\domain_-$ are all $C^\infty$ and such that the problem is \emph{nontrapping} (i.e.~all billiard trajectories starting in an exterior neighbourhood of $\domain_+:= \Rea^d\setminus \overline{\domain_-}$
and evolving according to the Hamiltonian flow defined by the symbol of \eqref{eq:1} escape from that neighbourhood after some uniform time),
then either (i) the propagation of singularities results of \cite{MeSj:82} combined with either the paramatrix argument of \cite{Va:75} or Lax--Phillips theory \cite{LaPh:89},
 or (ii) the defect-measure argument of  \cite{Bu:02}\footnote{The arguments in \cite{Bu:02} actually require that, additionally, $\partial\domain_-$ has no points where the tangent vector makes infinite-order contact with $\partial\domain_-$.}
proves the estimate that, given $k_0>0$ and $R>0$,
\beq\label{eq:res1}
\N{\gu}_{L^2(\domain_R)} + k \N{u}_{L^2(\domain_R)} \leq C_1(A, n, \domain_-, R, k_0) \N{f}_{L^2(\domain_+)}
\eeq
for all $k\geq k_0$, where $\domain_R:= \domain_+\cap B_R(\bze)$, 
and $C_1(A, n, \domain_-,R, k_0)$ is some (unknown) function of $A, n, \domain_-, R,$ and $k_0$, but is independent of $k$.
Without the nontrapping assumption, and assuming $n\equiv 1$ (but still $A,\domain_-\in C^\infty$), the Carleman-estimate argument of \cite{Bu:98} proves the estimate 
\beq\label{eq:res2}
\N{\gu}_{L^2(\domain_R)} + k \N{u}_{L^2(\domain_R)} \leq C_2(A, n, \domain_-, R, k_0) \exp\Big(k\,
 C_3(A, n, \domain_-, R, k_0) 
\Big)\N{f}_{L^2(\domain_+)},
\eeq
and this estimate has recently been obtained in \cite{Sh:17} for the case when $\domain_-=\emptyset$, $A\equiv I$, and $n$ is Lipschitz.
The estimates \eqref{eq:res1}, \eqref{eq:res2} are known as \emph{cut-off resolvent estimates}, with, e.g., \eqref{eq:res1}, often written as 
\beq\label{eq:res3}
\N{R_\chi(k)}_{L^2(\domain_+)\rightarrow L^2(\domain_+)} \leq \frac{C_1}{k}, \quad\N{R_\chi(k)}_{L^2(\domain_+)\rightarrow H^1 (\domain_+)} \leq C_4,
\eeq
where $R_\chi(k):= \chi (\nabla\cdot(A\nabla) +k^2n)^{-1}\chi$, for $\chi$ a $C^\infty$ function with compact support, is the \emph{cut-off resolvent} with zero Dirichlet boundary conditions on $\partial\domain_-$ and the Sommerfeld radiation condition at infinity, and $C_1, C_4$ are independent of $k$ (but dependent on $A,n, \domain_-, R,$ and $k_0$). 

The situation when $A$ and $n$ can jump is more complicated, since rays can be trapped by jumps of a particular sign. For the full-space problem (i.e.~no Dirichlet obstacle) and $A$ and $n$ jumping across a single common $C^\infty$ interface, \cite{Be:03a} proved that the bound \eqref{eq:res2} holds regardless of the sign of the jumps. The situation when the common interface is $C^\infty$ and has strictly positive curvature is well-understood thanks to the work of \cite{CaPoVo:99, PoVo:99, PoVo:99a}: when the sign of the jumps is such that the problem is trapping, there exists a sequence of resonances super-algebraically close to the real axis and the resolvent grows super-algebraically \cite{PoVo:99}, but when the jumps are nontrapping the estimate \eqref{eq:res3} holds \cite{CaPoVo:99}. 
Furthermore, when the interface is Lipschitz and star-shaped, the resolvent estimate \eqref{eq:res3} was recently proved for certain nontrapping jumps in \cite{MoSp:17}.

\paragraph{Informal discussion of the main results, and their novelty.}

Our focus is on determining conditions on $A$ and $n$ for which we can prove that the nontrapping resolvent estimate \eqref{eq:res1} holds, with 
the constant $C_1(A, n, \domain_-, R,k_0)$ given explicitly in terms of $A$, $n$, $R$, $k_0$, and geometric parameters describing $\domain_-$. 

We prove that the nontrapping resolvent estimate \eqref{eq:res1} holds for the EDP when $\domain_-$ is Lipschitz and star-shaped and one of the following three conditions holds
\beq\label{eq:A1a}
A(\bx) - (\bx\cdot\nabla)A(\bx) \geq \mu_1 \quad \tand \quad n(\bx)+ \bx\cdot\nabla n(\bx) \geq \mu_2, 
\eeq
\beq\label{eq:A2a}
2A(\bx) - (\bx\cdot\nabla)A(\bx) \geq \mu_3 \quad \tand \quad n(\bx)=1,
\eeq
\beq\label{eq:n3a}
A(\bx)=I \quad\tand \quad2n(\bx)+ \bx\cdot\nabla n(\bx) \geq \mu_4, 
\eeq
for almost every $\bx\in \domain_+$, where $\mu_j>0,\, j=1,\ldots,4$, and where the inequalities for $A$ are understood in the sense of quadratic forms.
For example, under \eqref{eq:A1a} we prove in Theorem \ref{thm:EDP1} below
that 
\beq\label{eq:bound1}
\mu_1 \N{\gu}^2_{L^2(\domain_R)} + \mu_2 k^2 \N{u}^2_{L^2(\domain_R)} \leq C_1 \N{f}^2_{L^2(\domain_+)},
\eeq
for all $k>0$, where
\beq\label{eq:C1}
C_1 := 4\left[\frac{R^2}{\mu_1} + \frac{1}{\mu_2}\left(R+ \frac{d-1}{2k}\right)^2\right];
\eeq
observe that if we impose the condition that $k\geq k_0$ for some $k_0>0$, then $C_1$ can be bounded above, independently of $k$, for $k\geq k_0$, and the bound \eqref{eq:bound1} is therefore of the form \eqref{eq:res1}.

We prove these results using the vector-field/commutator argument of Morawetz, with the commutation relations expressed as the identities in \S\ref{sec:4}. Recall that these identities were used in  \cite{MoLu:68, Mo:75} with the vector field $\bx$ to prove the resolvent estimate \eqref{eq:res1} when $\domain_-$ is smooth and star-shaped, $A\equiv I$, and $n\equiv 1$; we use the same vector field here, hence the star-shaped restriction on $\domain_-$ and the appearance of the vector field $\bx$ in \eqref{eq:A1a}-\eqref{eq:n3a}. One could conceivably generalise these results to wider classes of domains using the vector fields in \cite[Section 4]{MoRaSt:77} and \cite{BlKa:77}, but we do not explore this here.

We first prove these bounds under the assumption that both $A$ and $n$ are Lipschitz (Theorem \ref{thm:EDP1}), but we then use approximation and density arguments to prove them 
when $A, n\in L^\infty$, with \eqref{eq:A1a}--\eqref{eq:n3a} understood in a distributional sense (which then allows both $A$ and $n$ to be discontinuous); see  Condition \ref{cond:1altL} and Theorem \ref{thm:EDP2}. These arguments were inspired by analogous results in \cite{Th:06} in the setting of rough-surface scattering and $A\equiv I$.

The main novelty of these bounds is as follows.
\bit
\item To our knowledge, the bound \eqref{eq:bound1} is the first a priori bound on the solution of the EDP, with explicit dependence of the constant on $A$ and $n$, in a case when both $A$ and $n$ vary.
\item Our bound for $A, n \in L^\infty$ in Theorem \ref{thm:EDP2} implies new results about the well-posedness of the EDP in 3-d.
Indeed, as stated above, the best well-posedness results for this problem are for piecewise-Lipschitz $A$ via the unique continuation principle, the argument of  \cite{BaCaTs:12}, and Fredholm theory; we show in \S\ref{sec:EDPresults1} and \S\ref{sec:previous} below how Theorem \ref{thm:EDP2} proves well-posedness of the EDP for certain $A$ not covered by the argument in \cite{BaCaTs:12}.
\item Our bound for $A,n \in L^\infty$ in Theorem \ref{thm:EDP2} contains the first resolvent estimate for the transmission problem involving penetrable obstacles that are star-shaped, but not Lipschitz (and thus their boundaries can have, e.g., cusps). By  \cite[Lemma 2.3]{Vo:99}, this estimate then implies that the cut-off resolvent has a resonance-free strip beneath the real axis. Until recently the only results about resonances of the transmission problem were for $C^\infty$ convex obstacles with strictly positive curvature \cite{CaPoVo:99, PoVo:99, PoVo:99a, CaPoVo:01, Ga:15}; the recent paper \cite{MoSp:17} obtained results for Lipschitz star-shaped obstacles, and the present paper removes the Lipschitz assumption.
\eit
In addition, we explicitly show in \S\ref{sec:6} both (i) the relationship of the condition \eqref{eq:n3a} to the nontrapping of rays (partly using results in \cite{Ra:71}), and (ii) why one expects the condition \eqref{eq:n3a} to arise when using Morawetz's identities. 
The point (i) is implicit in results about geometric controllability of the wave equation (see, e.g., \cite{BaLeRa:92}, \cite{Ya:99}, \cite[Theorem 2.8]{AlCa:18})
and the reasons for (ii) are essentially known in the semiclassical-analysis community, but we could not find either (i) or (ii) explicitly stated in the literature, and neither appears to be known
by a large number of users of Morawetz's identities in the applied-analysis/numerical-analysis communities. 

\paragraph{Outline of the paper.}
In \S \ref{sec:2} we formulate the EDP and state our main results.
In \S \ref{sec:3} and \S\ref{sec:4} we collect preliminary results needed for the proofs, with \S\ref{sec:4} containing Morawetz's identities and associated results. 
In \S\ref{sec:5} and \S\ref{sec:5a} we prove the main results.
In \S\ref{sec:6} we explain the relationship of the condition \eqref{eq:n3a} to the nontrapping of rays, and why one expects this condition to arise when using Morawetz's identities. 
In Appendix \ref{sec:TEDP} we outline how our results can be extended to the TEDP.

\section{Formulation of the problem and statement of the main results}\label{sec:2}

\subsection{Formulation of the problem and geometric definitions}

\textbf{Notation:} $L^\infty(\domain)$ denotes complex-valued $L^\infty$ functions on a Lipschitz open set $\domain$. 
When the range of the functions is not $\Com$, it will be given in the second argument; e.g. 
$L^\infty(\domain,\Rea^{d\times d})$ denotes the space of $d\times d$ matrices with each entry a real-valued $L^\infty$ function on $\domain$.
We write $A\subset \subset B$ iff $A$ is \emph{compactly contained} in $B$ (i.e. $A$ is a compact subset of the open set $B$). 
We use $\gamma$ to denote the trace operator $H^1(\domain)\rightarrow H^{1/2}(\partial \domain)$ and use $\dist(\cdot,\cdot)$ to denote the distance function.

\begin{definition}[Exterior Dirichlet Problem (EDP)]\label{def:EDP}
Let $\domain_-$ be a bounded Lipschitz open set such that the open complement $\domain_+:= \Rea^d\setminus \overline{\domain_-}$ is connected. Let $\Gamma_D:= \partial \domain_-$. 
Given 
\bit
\item $f\in L^2(\domain_+)$ with $\supp f$ a compact subset of $\Rea^d$, 
\item $g_D\in H^1(\Gamma_D)$,
\item $n\in L^\infty(\domain_+,\Rea)$ such that $1-n $ has compact support and 
\beq\label{eq:nlimitsEDP}
0<\varmin \leq \var(\bx)\leq\varmax<\infty\,\, \text{ for almost every } \bx \in \domain_+,
\eeq
\item $A \in L^\infty (\domain_+ , \Rea^{d\times d})$ such that $I -A$ has compact support, $A$ is symmetric, and there exist $0<A_{\min}\leq A_{\max}<\infty$ such that
\beq\label{eq:AellEDP}
 A_{\min} |\bxi|^2\leq\big(A(\bx) \bxi\big) \cdot\overline{ \bxi}  \leq A_{\max}|\bxi|^2 \quad\text{ for almost every }\bx \in \domain_+ \text{ and for all } \bxi\in \Com^d,
\eeq
\eit
we say $u\in H^1_{\rm{loc}}(\domain_+)$ satisfies the \emph{exterior Dirichlet problem} if 
\beq\label{eq:EDP1}
\opL_{A,n} u:= \nabla\cdot(A \gu ) + k^2 n u = -f \quad \tin \domain_+,
\eeq
\beqs
\gamma u =g_D \quad\ton \Gamma_D,
\eeqs
and $u$ satisfies the Sommerfeld radiation condition 
\beq\label{eq:src}
\pdiff{u}{r}(\bx) - \ri k u(\bx) = o \left( \frac{1}{r^{(d-1)/2}}\right)
\eeq
as $r:= |\bx|\tendi$, uniformly in $\hatx:= \bx/r$.
\end{definition}

\noi Some remarks: 

(i) The equation \eqref{eq:EDP1} is understood in the following weak sense:
\beq\label{eq:weak}
\int_{\domain_+} 
\gu \cdot(A\overline{\nabla\phi}) 
- k^2 n u\overline{\phi} 
=\int_{\domain_+} f\, \overline{\phi} \quad\tfa \phi \in C^\infty_{0}(\domain_+),
\eeq
where $C^\infty_{0}(\domain_+):= \{ \phi \in C^\infty(\domain_+) : \supp \,\phi \text{ is a compact subset of } \domain_+\}.$

(ii) We can legitimately impose the radiation condition \eqref{eq:src} on the function $u \in H^1_{\loc}(\domain_+)$ since 
 $u$ satisfies the equation $\Delta u +k^2 u=0$ outside a ball of finite radius, and then  $u$ is $C^\infty$ outside this ball by elliptic regularity.

(iii) One usually prescribes $g_D\in H^{1/2}(\Gamma)$, but the Morawetz identities that we use to obtain a bound on the solution require that $g_D \in H^1(\Gamma)$.

(iv) If  $\domain_-=\emptyset$, then $\domain_+=\Rea^d$ and the EDP becomes the full-space problem.

\paragraph{Geometric definitions:}
Let $B_{a}(\bx_0):=\{\bx\in\Rea^d,|\bx-\bx_0|_2<a\}$, where $|\cdot|_2$ denotes the vector 2-norm (from here on we drop the subscript $2$). When $\bx_0=\bze$ we write $B_a$ for $B_a(\bx_0)$. We let $\domain_R:= \domain_+ \cap B_R$.

We now define the 
notions of \emph{star-shaped} and \emph{star-shaped with respect to a ball}.

\begin{definition}
(i) $\domaingen$ is \emph{star-shaped with respect to the point $\bx_0$} if, whenever $\bx \in \domaingen$, the segment $[\bx_0,\bx]\subset \domaingen$.

\noi (ii) $\domaingen$ is \emph{star-shaped with respect to the ball $B_{a}(\bx_0)$} if it is star-shaped with respect to every point in $B_{a}(\bx_0)$.

\end{definition}

These definitions make sense even for non-Lipschitz $\domaingen$, but when $\domaingen$ is Lipschitz one can characterise star-shapedness with respect to a point or ball in terms of 
$(\bx-\bx_0)\cdot\bnu(\bx)$ for $\bx \in \partial\domaingen$, where $\bnu(\bx)$ is the outward-pointing unit normal vector at $\bx\in \partial\domaingen$.

\begin{lemma}\textbf{\emph{(\cite[Lemma 5.4.1]{Mo:11})}}\label{lem:star}
(i) If $\domaingen$ is Lipschitz, then it is star-shaped with respect to $\bx_0$ if and only if $(\bx-\bx_0)\cdot\bnu(\bx)\geq 0$ for all $\bx \in\partial\domaingen$ for which 
$\bnu(\bx)$ is defined.

\noi (ii) If $\domaingen$ is Lipschitz, then $\domaingen$ is star-shaped with respect to $B_{a}(\bx_0)$ if and only if $(\bx-\bx_0) \cdot \bnu(\bx) \geq {a}$ for all  $\bx \in \partial \domaingen$ for which $\bnu(\bx)$ is defined.
\end{lemma}

\subsection{The main results (with $A$ and $n$ satisfying \eqref{eq:A1a})}
\label{sec:EDPresults1}

As mentioned in \S\ref{sec:intro}, there are three sets of conditions on $A$ and $n$,  \eqref{eq:A1a}-\eqref{eq:n3a}, under which we prove results. For clarity of exposition, in this subsection we describe the main results when $A$ and $n$ satisfy \eqref{eq:A1a}; in \S\ref{sec:EDPresults2} we then describe the results when $A$ and $n$ satisfy either \eqref{eq:A2a} or \eqref{eq:n3a}.

Our first result is the nontrapping resolvent estimate when $A$ and $n$ satisfy \eqref{eq:A1a} and are both Lipschitz.

\begin{condition}[$A$ and $n$ both Lipschitz, $g_D\equiv 0$, $\domain_-$ star-shaped]\label{cond:1}
$d=2,3$, $\domain_-$ is star-shaped with respect to the origin, $A\in C^{0,1}(\overline{\domain_+}, \Rea^{d\times d})$, $n\in C^{0,1}(\overline{\domain_+}, \Rea)$, $g_D\equiv 0$, and there exist $\mu_1, \mu_2>0$ such that
\beq\label{eq:A1}
A(\bx) - (\bx\cdot\nabla)A(\bx) \geq \mu_1,\,\, \text{ in the sense of quadratic forms, for almost every }\bx\in \domain_+,
\eeq
and 
\beq\label{eq:n1}
n(\bx)+ \bx\cdot\nabla n(\bx) \geq \mu_2 \quad\text{ for almost every }\bx\in \domain_+.
\eeq
\end{condition}

Recall that if $\domain$ is a bounded Lipschitz open set then $C^{0,1}(\overline{\domain})=W^{1,\infty}(\domain)$ (see, e.g., \cite[\S4.2.3, Theorem 5]{EvGa:92}), and so the conditions \eqref{eq:A1} and \eqref{eq:n1} make sense for Lipschitz $A$ and $n$.

\begin{theorem}[Bounds on the EDP under Condition \ref{cond:1}]\label{thm:EDP1}
If $\domain_-, A$, $n$, $f$, and $g_D$ satisfy the requirements in the definition of the EDP (Definition \ref{def:EDP}), along with the requirements in Condition \ref{cond:1}, then the solution of the EDP exists and is unique. Furthermore, given $R>0$ with $\supp(I-A), $ $\supp(1-n),$ and $\supp \,f$ all compactly contained in $\domain_R$, then \eqref{eq:bound1} holds for all $k>0$ where $C_1$ is given by \eqref{eq:C1}.
\end{theorem}

When $A$ and $n$ are $C^1$, the conditions \eqref{eq:A1} and \eqref{eq:n1} can be rewritten as
\beqs
 -r^2 \pdiff{}{r} \left(\frac{A}{r}\right)\geq \mu_1  \quad\tand\quad\pdiff{}{r}\big( rn \big) \geq \mu_2,
\eeqs
suggesting that the most general $L^\infty$ $A$ and $n$ for which we can prove a bound are 
\beqs
A(\bx) = \mu_1 I - r \widetilde{\Pi}(\bx) \quad\tand\quad n(\bx) = \mu_2 + \frac{\widetilde{\pi}(\bx)}{r},
\eeqs
where both $\widetilde{\Pi}\in L^\infty(\domain_{+},\Rea^{d\times d})$ and
 $\widetilde{\pi}\in L^\infty(\domain_{+},\Rea)$ are monotonically non-decreasing in the radial direction.
To minimise technicalities arising from the factors of $r$ (and the singularity of $1/r$ at the origin), we prove a bound under the following slightly-more-restrictive conditions.

\begin{condition}[An analogue of Condition \ref{cond:1} with $A$ and $n$ in $L^\infty(\Oe)$]\label{cond:1altL}
$d=2,3$,
$\domain_-$ is star-shaped with respect to the origin, $g_D\equiv 0$, $A$ and $n$ satisfy the requirements in the Definition of the EDP (Definition \ref{def:EDP}), 
\beqs
A(\bx) = A_{\max} I  - \Pi(\bx)\quad\text{ for almost every }\, \bx\in \domain_+,
\eeqs
where $\Pi\in L^\infty(\domain_{+},\Rea^{d\times d})$ 
is monotonically non-decreasing in the radial direction in the sense of quadratic forms, i.e. for all $h\geq 0$,
\beqs
\essinf_{\bx\in \domain} \big[ \Pi(\bx + h\vecer )- \Pi(\bx)\big] \geq 0
\eeqs
in the sense of quadratic forms, where $\vecer$ is the unit vector in the radial direction, and 
\beqs
n(\bx) = \varmin  + \pi(\bx)\quad\text{ for almost every }\, \bx\in \domain_+,
\eeqs
where $\pi\in L^\infty(\domain_{+},\Rea)$ 
is monotonically non-decreasing in the radial direction, i.e. for all $h\geq 0$,
\beqs
\essinf_{\bx\in \domain} \big[ \pi(\bx + h\vecer )- \pi(\bx)\big] \geq 0.
\eeqs
\end{condition}

\begin{theorem}[Bounds on the EDP under Condition \ref{cond:1altL}]\label{thm:EDP2}
If $\domain_-, A,$ and $n$ satisfy the requirements in the definition of the EDP (Definition \ref{def:EDP}), along with the requirements in Condition \ref{cond:1altL}, then the solution of the EDP exists and is unique. Furthermore, given $R>0$ with $\supp(I-A), $ $\supp(1-n),$ and $\supp\, f$ all compactly contained in $\domain_R$, then 
\eqref{eq:bound1} holds for all $k>0$, with $\mu_1=A_{\min}$ and $\mu_2= \varmin$.
\end{theorem}

We now highlight three particular situations in which Condition \ref{cond:1altL} is satisfied with $A$ and $n$ piecewise constant.
In these examples, we use $\boldsymbol{1}_G$ to denote the indicator function of a set $G$.

\begin{figure}
  \centering
  \subfloat[The subdomains on which $A$ and $n$ are constant in Condition \ref{cond:simple}]
  {
\scalebox{1}{
\begin{tikzpicture}[scale=3,inner sep=0.5mm]

\draw[domain = 0:6.2832,smooth,variable=\t,samples=400]
plot ({\t r}:{(0.4-sin(\t r))*1.35});

\node at (0,-0.31) [circle,draw,fill=black] {};

\node at (0,-0.5) {{\large $\Omega_1$}}; 
\node at (-0.6,-0.7) {{\large $\Omega_2$}};
\node at (0,0.2) {{\large $\Omega_3$}};
\end{tikzpicture}
    }
    \label{fig:1}
  }
  \quad
  \subfloat[The subdomains on which $A$ and $n$ are constant in Condition \ref{cond:Yves}.]
  {
  \scalebox{1}{

\begin{tikzpicture}[scale=3,inner sep=0.5mm]

\fill[fill=black] (0,0) circle [radius =1];
\fill[fill=white] (0,0) circle [radius =11/12];

\foreach \x in {1,2,...,11}
\draw (0,0) circle [radius = (\x/(\x+1))];

\draw (0,0) circle [radius = 1];

\node at (0,0) [circle,draw,fill=black] {};

\end{tikzpicture}

    }
    \label{fig:2}
  }

  \subfloat[An admissible $\Omega_T$ in Condition \ref{cond:T} (shape taken from \cite{We:17}).]
  {
  \scalebox{1}{

\begin{tikzpicture}[scale=3,inner sep=0.5mm]

\draw[domain = 1.56:1.58,smooth,variable=\t,samples=100] 
plot ({\t r}:{(2 - 2 * sin(\t r) + sin(\t r)*sqrt(abs(cos(\t r)))/(sin(\t r)+1.4))/2}); 

\draw[domain = 4.65:4.75,smooth,variable=\t,samples=100] 
plot ({\t r}:{(2 - 2 * sin(\t r) + sin(\t r)*sqrt(abs(cos(\t r)))/(sin(\t r)+1.4))/2}); 

\draw[domain = 0:1.56,smooth,variable=\t,samples=200] 
plot ({\t r}:{(2 - 2 * sin(\t r) + sin(\t r)*sqrt(abs(cos(\t r)))/(sin(\t r)+1.4))/2}); 

\draw[domain = 1.58:4.66,smooth,variable=\t,samples=200] 
plot ({\t r}:{(2 - 2 * sin(\t r) + sin(\t r)*sqrt(abs(cos(\t r)))/(sin(\t r)+1.4))/2}); 

\draw[domain = 4.749:6.29,smooth,variable=\t,samples=200] 
plot ({\t r}:{(2 - 2 * sin(\t r) + sin(\t r)*sqrt(abs(cos(\t r)))/(sin(\t r)+1.4))/2}); 

\node at (0,-0.31) [circle,draw,fill=black] {};

\node at (-0.4,-0.6) {{\large $\Omega_T$}};


\end{tikzpicture}

    }
    \label{fig:3}
  }
  
  \caption{Examples illustrating Condition \ref{cond:1altL}; in all three, the black dot marks the origin.}
\end{figure}
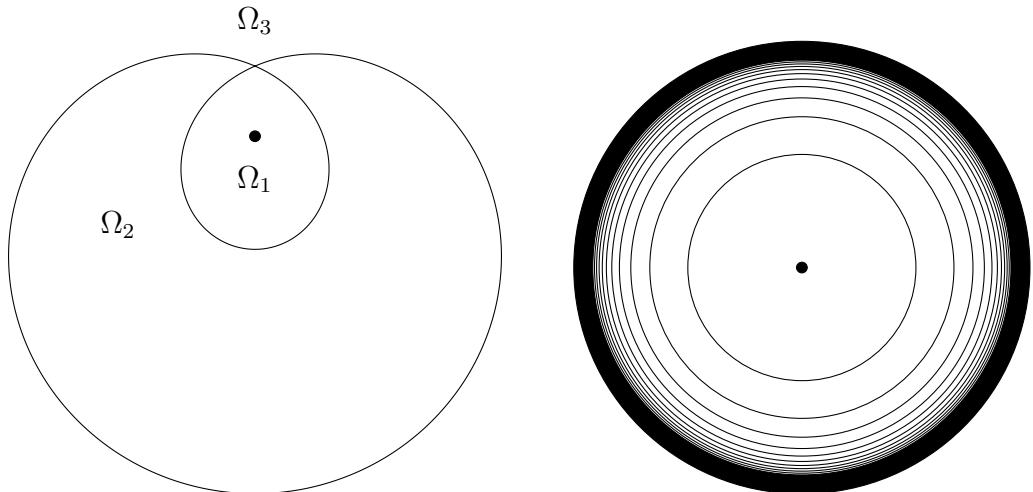
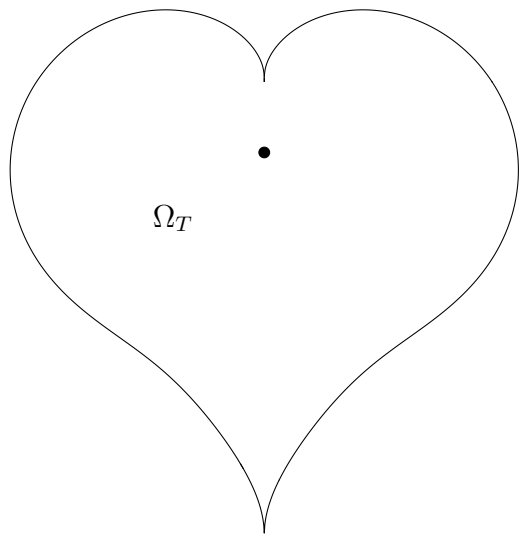

\begin{condition}[First particular case of Condition \ref{cond:1altL}]\label{cond:simple}
$\Oi=\emptyset$, $\Omega_1,\, \Omega_2,$ and $\Omega_3$ are as in Figure \ref{fig:1}, and 
\beqs
A(\bx) = \left(\sum_{j=1}^3 a_j \boldsymbol{1}_{\Omega_j}\right)I, \quad n(\bx) =\sum_{j=1}^3 n_j \boldsymbol{1}_{\Omega_j}
\eeqs
where $a_1\geq a_2\geq a_3=1$ and $0<n_1\leq n_2 \leq n_3=1$.
\end{condition}

\begin{condition}[Second particular case of Condition \ref{cond:1altL}]\label{cond:Yves}
$\Oi=\emptyset$, 
\beqs
A(\bx) = 
\left( \boldsymbol{1}_{\Rea^d\setminus B_1(\bze)}(\bx)+
\sum_{j=1}^\infty a_j \boldsymbol{1}_{G_j}(\bx)
\right)I 
\eeqs
and 
\beqs
n(\bx) = 
\boldsymbol{1}_{\Rea^d\setminus B_1(\bze)}(\bx)+
\sum_{j=1}^\infty n_j \boldsymbol{1}_{G_j}(\bx),
\eeqs
where 
\beqs
A_{\max}= a_1\geq \ldots \geq a_j \geq a_{j+1}\geq \ldots \geq 1=A_{\min},
\eeqs
\beqs
n_{\min}= n_1\leq \ldots \leq n_j \leq n_{j+1}\leq \ldots \leq 1=n_{\max},
\eeqs
and the sets $G_j$ are defined by 
\beqs
G_j:=B_{j/(j+1)}(\bze) \setminus B_{(j-1)/j}(\bze).
\eeqs
\end{condition}

\begin{condition}\textbf{\emph{(Particular case of Condition \ref{cond:1altL} corresponding to the transmission problem)}}\label{cond:T}
$\Oi= \emptyset$, 
\beq\label{eq:AnT}
A(\bx) = 
\left(a_i\boldsymbol{1}_{\Omega_T}(\bx)+ a_o\boldsymbol{1}_{\Rea^d\setminus\overline{\Omega_T}}(\bx) \right)I 
\quad\tand\quad n(\bx)
= n_i\boldsymbol{1}_{\Omega_T}(\bx) +n_o\boldsymbol{1}_{\Rea^d\setminus\overline{\Omega_T}}(\bx),
\eeq
where $a_o, a_i, n_o,$ and $n_i$ are positive real numbers satisfying
\beq\label{eq:AnT2}
1=a_o \leq a_i \quad\tand \quad 1=n_o \geq n_i,
\eeq
and $\Omega_T$ is a $C^0$ bounded open set that is star-shaped with respect to a point.
\end{condition}

When $A$ and $n$ satisfy Condition \ref{cond:T} and $\Omega_T$ is additionally Lipschitz, the Helmholtz EDP of Definition \ref{def:EDP} reduces to the \emph{Helmholtz transmission problem} (see, e.g., \cite[Lemma 4.19]{Mc:00}).

\begin{definition}[The Helmholtz transmission problem]\label{def:HTP}
Let $k>0$ 
and let $n_i,n_o,a_i,a_o$ be positive real numbers.
Let $\Omega_T$ be a bounded connected Lipschitz open set.
Let $f_i\in L^2(\Oin)$, $f_o\in L^2(\Oout)$, and assume that $f_o$ has compact support.
The Helmholtz transmission problem is: find $u_i \in H^1(\Omega_T)$ and $u_0\in H^1_{\loc}(\Rea^d\setminus\overline{\Omega_T})$ such that
\begin{align}
\begin{aligned}
a_i\Delta u_i+k^2 n_i u_i &=f_i &&\text{ in }\Oin,\\
a_o\Delta u_o+k^2 n_o u_o &=f_o &&\text{ in } \Oout,\\
u_o&= u_i&&\text{ on } \partial\Omega_T,\\
a_o\pdiff{u_o}{\nu}&=a_i\pdiff{u_i}{\nu}&&\text{ on } \partial\Omega_T,
\end{aligned}
\label{eq:BVP}
\end{align}
and $u_o$ satisfies the Sommerfeld radiation condition \eqref{eq:src} with $k$ replaced by $k\sqrt{n_o/a_o}$.
\end{definition}

Since the definition of the conormal derivative requires that $\Omega_T$ be Lipschitz (see, e.g., \cite[Lemma 4.3]{Mc:00}), the transmission problem for non-Lipschitz $\Omega_T$ must be understood as the EDP with $A$ and $n$ given by \eqref{eq:AnT}.

Theorem \ref{thm:EDP2} can be used to obtain a result about the resonances of the EDP under Condition \ref{cond:1altL};
and thus in particular about resonances of the transmission problem corresponding to Condition \ref{cond:T}.
 To state this result, we introduce the following notation. Let $R(k)$ denote the solution operator of the EDP of Definition \ref{def:EDP} when $g_D\equiv 0$; i.e.~$R(k): f \mapsto u$.
Although $R(k)$ depends 	also on $A$ and $n$, in what follows we consider these fixed and consider $k$ as variable.
Given  $\chi \in C^\infty_{0}(\Rea^d)$ such that $\chi\equiv 1$ in a neighbourhood of $\Oi$, define
the \emph{cut-off resolvent} 
\beqs
R_\chi(k):= \chi R(k)\chi; 
\eeqs
then $R_\chi(k): L^2(\Omega_+)\rightarrow H^1(\Omega_+)$ for $k\in \Rea\setminus\{0\}$. The \emph{resonances} of the EDP are then defined to be the poles of the meromorphic continuation of $R_\chi(k)$ into $\Im k<0$.

\begin{corollary}[Resonance-free strip beneath the real-$k$ axis]\label{cor:res1}
The operator family $R_\chi(k)$ is well-defined and holomorphic for $\Im k>0$.
If $\domain_-, A,$ and $n$ satisfy the requirements in the definition of the EDP (Definition \ref{def:EDP}), along with the requirements in Condition \ref{cond:1altL}, 
then there exist $\cC_j>0$, $j=1,2,3$ (independent of $k$ but dependent on $A$, $n$, and $\Oi$) such that $R_\chi(k)$
extends from the upper-half plane to a holomorphic operator family on 
$|\Re k|\geq \cC_1  , \Im k\geq -\cC_2$
satisfying the  estimate 
\beq\label{eq:55}
\N{R_\chi(k)}_{L^2(\Oe) \rightarrow L^2(\Oe)} \leq\frac{\cC_3}{|k|}
\eeq 
in this region.
\end{corollary}

This corollary follows from the bound of Theorem \ref{thm:EDP2} using the result \cite[Lemma 2.3]{Vo:99}. 
Recall that \cite[Lemma 2.3]{Vo:99} takes a resolvent estimate for $k$ real (such as \eqref{eq:res3}/\eqref{eq:bound1}) and converts it into a resolvent estimate in a strip beneath the real-$k$ axis (such as \eqref{eq:55}).
In principle, one could go into the details of  \cite[Lemma 2.3]{Vo:99} and make the width of the strip ($\cC_2$ in Corollary \ref{cor:res1}) explicit in the constant from the bound for $k$ real. 
Since  Theorem \ref{thm:EDP2}  gives an explicit expression for that constant, 
we would then have an explicit lower bound for the width of the strip.

\bre[Bounds for non-Lipschitz $\domain_-$]\label{rem:CWM}
In formulating the EDP we assumed that $\Gamma_D:=\partial \domain_-$ was Lipschitz. 
We see below that our proofs of the bounds in Theorems \ref{thm:EDP1} and \ref{thm:EDP2} also hold in the non-Lipschitz case; i.e.~with the only assumptions on $\Oi$ that it is a $C^0$ bounded open set that is star-shaped with respect to the origin (the BVP is then understand as the variational problem \eqref{eq:EDPvar} below). These arguments were inspired by the arguments in \cite{ChMo:08}, where the bound \eqref{eq:bound1} was proved for such $\Oi$ when $A\equiv I$ and $n\equiv 1$.
\ere

\bre[$H^2$-regularity]\label{rem:H2}
The standard finite-element-error analysis of the Helmholtz equation is most comfortable when the solution of the adjoint BVP is in $H^2$ ($H^2_{\rm{loc}}$ for exterior problems), with a bound on the $H^2$-norm of the adjoint BVP playing a key role (for the reasons behind this see, e.g., \cite[Theorem 6.32 and Remarks 26 and 31]{Sp:15} and the references therein).
When $\partial\domain_- \in C^{1,1}$, $A\in C^{0,1}(\overline{\Oe},\Rea^{d\times d})$, and $g_D\in H^{3/2}(\Gamma_D)$, the solution of the EDP (Definition \ref{def:EDP}) is in $H^2_{\rm{loc}}(\Oe)$, and 
the bound in Theorem \ref{thm:EDP1} can then be combined with 
\cite[Theorem 4.18 (i)]{Mc:00} or \cite[Theorems 2.3.3.2 and 2.4.2.5]{Gr:85} to produce a bound on the $H^2$-norm. (When $A\not\equiv I$, however, 
\cite[Theorem 4.18 (i)]{Mc:00} and \cite[Theorem 2.3.3.2]{Gr:85} do not state explicitly how the constants in the bound depend on $A$, although this dependence could be in principle be determined.)
\ere 

\bre[$A$ and $n$ perturbations]\label{rem:pert}
The case when $A$ and $n$ are $L^\infty$ perturbations of $A_0$ and $n_0$ satisfying Condition \ref{cond:1}/\ref{cond:1altL} is investigated more in \cite[Remark 1.14]{PeSp:18}, but we mention here the particular case when $n=n_0+ \eta$. Writing the PDE as
\beq\label{eq:pert}
\nabla \cdot(A\gu) +k^2 n_0 u = -f - k^2 \eta u
\eeq
and applying Theorem \ref{thm:EDP1}/\ref{thm:EDP2}, we see that if $k\|\eta\|_{L^\infty(\domain_R)}$ is sufficiently small, then 
bounds similar to those in Theorem \ref{thm:EDP1}/\ref{thm:EDP2} hold, since 
one can absorb the contribution from the $k^2 \eta u$ term in \eqref{eq:pert} into the $\mu_2 k^2\|u\|^2_{L^2(\OR)}$ appearing on the left-hand side of \eqref{eq:bound1}.
\ere

In Theorems \ref{thm:EDP1} and \ref{thm:EDP2} the right-hand side of the PDE, $f$, is in $L^2(\domain_+)$ with compact support.
A standard argument appearing in, e.g., \cite[Text between Lemmas 3.3 and 3.4]{ChMo:08}, \cite[Proof of Corollary 1.10]{BaSpWu:16} can then be used to prove bounds for $f \in (H^1_{0, D}(\domain_R))'$  where $H^1_{0,D}(\domain_R)$ is the  space in which the variational formulation of the EDP with $g_D\equiv 0$ is formulated (i.e.~$H^1$ with a zero Dirichlet boundary condition on $\Gamma_D$; see \eqref{eq:spaceEDP} below). A bound with data in $(H^1)'$ is then equivalent to a bound on the inf-sup constant (see, e.g., \cite[Theorem 2.1.44]{SaSc:11}).
In stating this corollary, we use the norm
\beq\label{eq:1knorm}
\N{v}^2_{\HoDkk} := \N{\gv}^2_{L^2(\domain_R)} + k^2 \N{v}^2_{L^2(\domain_R)} \quad \tfor v \in H^1_{0,D}(\domain_R).
\eeq

\begin{corollary}[Bound on the inf-sup constant of the EDP under Condition \ref{cond:1}/\ref{cond:1altL}]
\label{cor:H1}
Let $\Omega_-, A$, and $n$ satisfy the requirements in the definition of the EDP (Definition \ref{def:EDP}), along with the requirements in Condition \ref{cond:1}/\ref{cond:1altL}.
Given $F\in (H^1_{0,D}(\domain_R))'$ let $u$ be the solution of the variational formulation of the EDP \eqref{eq:EDPvar} below. Then $u$ exists, is unique, and satisfies the bound
\beq\label{eq:H1bound}
\N{u}_{\HoDkk} \leq \frac{1}{\min(A_{\min},\varmin)}\left( 1 + 2\sqrt{\frac{C_1}{\min(\mu_1,\mu_2)}}\varmax  k\right) \N{F}_{(\HoDkk)'}
\eeq
for all $k>0$, where $C_1$ is given by \eqref{eq:C1}. Thus, the inf-sup constant of the sesquilinear form $a(\cdot,\cdot)$ defined by \eqref{eq:EDPa} satisfies the lower bound
\beqs
\inf_{u\in H^1_{0,D}(\domain)\setminus\{0\}} \sup_{v\in H^1_{0,D}(\domain)\setminus\{0\}} \frac{\big|a(u,v)\big|}{\N{u}_{\HoDkk}\N{v}_{\HoDkk}} \geq 
\min(A_{\min},\varmin)\left(1 + 2\sqrt{\frac{C_1}{\min(\mu_1,\mu_2)}}\varmax  k\right)^{-1}.
\eeqs
\end{corollary}
 
\subsection{The analogues of the results in \S\ref{sec:EDPresults1} when $A$ and $n$ satisfy \eqref{eq:A2a} or \eqref{eq:n3a}}\label{sec:EDPresults2}

\begin{condition}[$A$ Lipschitz, $n\equiv 1$, $g_D\equiv 0$, $\domain_-$ star-shaped]\label{cond:2}
$d=2,3$, $\domain_-$ is star-shaped  with respect to the origin, $A\in C^{0,1}(\overline{\domain_+}, \Rea^{d\times d})$, $n\equiv 1$, $g_D\equiv 0$, and there exists $\mu_3>0$ such that
\beq\label{eq:A2}
2A(\bx) - (\bx\cdot\nabla)A(\bx) \geq \mu_3,\,\, \text{ in the sense of quadratic forms, for almost every }\bx\in \domain_+.
\eeq
\end{condition}

\begin{condition}[$n$ Lipschitz, $A\equiv I$, $g_D\not\equiv 0$, $\domain_-$ star-shaped w.r.t.~a ball]\label{cond:3b}
$d=2,3$, $\domain_-$ is star-shaped with respect to a ball centred at the origin, $A\equiv I$, $n\in C^{0,1}(\overline{\domain_+})$ and there exists $\mu_4>0$ such that 
\beq\label{eq:n3}
2n(\bx)+ \bx\cdot\nabla n(\bx) \geq \mu_4 \quad\text{ for almost every }\bx\in \domain_+.
\eeq
\end{condition}

\noi Observe that the conditions on $A$ and $n$ varying on their own, i.e. \eqref{eq:A2} and \eqref{eq:n3}, are less restrictive than the conditions on $A$ and $n$ varying together, i.e.~\eqref{eq:A1} and \eqref{eq:n1}.

\begin{theorem}[Bounds on the EDP under Conditions \ref{cond:2} and \ref{cond:3b}]\label{thm:EDP3}

\

\noi (i) If $\domain_-, A$, $n$, $f$, and $g_D$ satisfy the requirements in the definition of the EDP (Definition \ref{def:EDP}), along with the requirements in Condition \ref{cond:2}, then the solution of the EDP exists and is unique. Furthermore, given $R>0$ with both $\supp(I-A)$ and $\supp\, f$ compactly contained in $\domain_R$,
\beq\label{eq:bound2}
\mu_3 \left(\N{\gu}^2_{L^2(\domain_R)} + k^2 \N{u}^2_{L^2(\domain_R)} \right)\leq C_2 \N{f}^2_{L^2(\domain_+)},
\eeq
for all $k>0$, where 
\beqs
C_2 = \frac{4}{\mu_3}\left( R^2 +\left(R+1+ \frac{d}{2k}\right)^2\right)\left(1 + 4A_{\max} + \frac{24(A_{\max})^2}{k^2}\right)^2
+ \frac{4\mu_3}{k^2}.
\eeqs

\

\noi (ii) If $\domain_-, A$, $n$, $f$, and $g_D$ satisfy the requirements in the definition of the EDP (Definition \ref{def:EDP}) along with the requirements in Condition \ref{cond:3b}, then the solution of the EDP exists and is unique. 
Let $L_D:= \max_{\bx\in\Gamma_D}|\bx|$, and let $aL_D$ be the radius of the ball with respect to which $\domain_-$ is star-shaped.
Given $R>0$ with both $\supp(1-n)$ and $\supp\, f$ compactly contained in $\domain_R$, 
\begin{align}\nonumber
\mu_4\left(\N{\gu}^2_{L^2(\domain_R)} +\right.&\left.  k^2 \N{u}^2_{L^2(\domain_R)}\right)+ \frac{aL_D}{2} \N{\dudnu}^2_{L^2(\Gamma_D)} \\
&\leq C_3 \N{f}^2_{L^2(\domain_+)}  + C_4  \N{\nabla_{\Gamma_D}g_D}^2_{L^2(\Gamma_D)} + C_5 k^2 \N{g_D}^2_{L^2(\Gamma_D)}\label{eq:bound3b}
\end{align}
for all $k\geq \sqrt{3/8}R^{-1}$, 
where 
\beq\label{eq:C3}
C_3 
= \frac{4}{\mu_4}\left( R^2+\left(2R+ \frac{d-2}{2k}\right)^2\right)\left(1 + \frac{3}{2}\varmax\right)^2
+ \frac{2\mu_4}{k^2 \varmax}.
\eeq
\beqs
C_4 
= 2\left(1 + \frac{3}{2}\varmax \right)L_D\left(1 + \frac{4}{a}\right),
\eeqs
and 
\beqs
C_5
= 2\left(1 +\frac{3}{2}\varmax \right)\frac{4}{aL_D}\left(2R+ \frac{d-2}{2k }\right)^2+ \frac{2\mu_4^2}{aL_D k^2}.
\eeqs

\noi Observe that, given $k_0>0$, each of $C_j, j=2,\ldots, 5$ can be bounded above, independently of $k$, for $k\geq k_0$.
\end{theorem}

Analogues of Remarks \ref{rem:CWM}, \ref{rem:H2}, and \ref{rem:pert}, and Corollary \ref{cor:H1} hold for the results of Theorem \ref{thm:EDP3}. The only exception is that the assumption that $\Oi$ is Lipschitz cannot be removed (as described in Remark \ref{rem:CWM}) from Part (ii) of Theorem \ref{thm:EDP3}, since 
 the bound \eqref{eq:bound3b} involves the normal derivative  on the boundary, and 
one needs the domain to be Lipschitz for the normal derivative to be well-defined. 

\bre[Star-shaped vs.~star-shaped with respect to a ball]
In Conditions \ref{cond:1} and \ref{cond:2} we assumed that $g_D\equiv 0$ and $\domain_-$ is star-shaped,
and we obtained bounds under these conditions in Theorems \ref{thm:EDP1} and \ref{thm:EDP3}(i) (respectively) using Morawetz identities.
Morawetz identities can be used to obtain bounds when $g_D\not\equiv 0$ if $\domain_-$ is star-shaped with respect to a ball, but when $A\not\equiv I$ the calculations involving the surface gradient are more involved than the case $A\equiv I$ (compare \eqref{eq:nT1} and \eqref{eq:7a} below); we therefore exclude the case $g_D\not\equiv 0$ from Theorems \ref{thm:EDP1} and 
\ref{thm:EDP3}(i) for brevity. 

On the other hand, Condition \ref{cond:3b} (where $A\equiv I$) allows $g_D\not\equiv 0$ and $\domain_-$ that are star-shaped with respect to a ball. Recall that $g_D\not\equiv 0$ is needed to cover scattering of a plane-wave; indeed, in this case, $u$ in Definition \ref{def:EDP} is the scattered field and $g_D$ equals minus the trace of the incident plane wave (which is not zero); see, e.g., \cite[Page 107]{ChGrLaSp:12}.
Note that if $\domain_-$ is only star-shaped and additionally $g_D\equiv 0$, then the bound \eqref{eq:bound3b} holds with all the norms over $\Gamma_D$ removed.
\ere

\bre[Bounds on the full-space problem]\label{rem:trans}
Taking the obstacle $\domain_-$ to be the empty set, changing $\domain_R$ to $B_R$, and ignoring all the norms on $\Gamma_D$, the results of Theorem \ref{thm:EDP3} apply to the full-space problem.
\ere

\bre[Dimensional arguments]
From dimensional arguments, we expect the factors $C_1, C_2$, and $C_3$ to have the dimension of $($length$)^2$, and $C_4$ and $C_5$ to have dimension of length (see, e.g., \cite[Remarks 3.6 and 3.8]{MoSp:14}). This is the case provided that we (a) interpret the $R+1$ in $C_2$ as having the dimension of length, and (b) interpret the $(A_{\max})^2/k^2$ in $C_2$ as having an extra $($length$)^{-2}$ factor (and thus being dimensionless). 
These apparent discrepancies arise because we choose a distance to be equal to one in the proof of the bound \eqref{eq:bound2}: we apply the Morawetz identity in a ball of radius $R+1$, as opposed to $R+b$ for some arbitrary $b$ with dimensions of length. Thus in (a) above we have $R+1$ instead of $R+b$ and in (b) we have $k^2= ((R+1)-R)^2 k^2$ instead of $b^2 k^2= ((R+b)-R)^2 k^2$.
\ere 
 
\subsection{Discussion of previous work}\label{sec:previous}

This discussion focuses on the exterior Dirichlet and full-space problems; 
nevertheless, there is a substantial literature on well-posedness and $k$-explicit bounds for
the Helmholtz equation when \emph{either} the obstacle $\Oi$ is unbounded \emph{or} at least one of $\supp (I-A)$, $\supp(1-n)$, and $\supp\, f$ is not compact (e.g.~scattering by infinite rough surfaces and/or infinite rough layers); see \cite[Remark 3.9]{MoSp:17} and the references therein for an overview of results about these problems.

\paragraph{Existence and uniqueness results.}
As discussed in \S\ref{sec:intro}, the unique continuation principle (UCP) gives uniqueness of the solution of \eqref{eq:1} when $A$ is Lipschitz and $n \in L^\infty$.

The Baire-category argument of \cite{BaCaTs:12} 
uses the fact that a UCP is known for the time-harmonic Maxwell system with Lipschitz coefficients \cite{NgWa:12} to prove uniqueness of the solution of the time-harmonic Maxwell problem posed in $\Rea^3$ with \emph{piecewise}-Lipschitz coefficients, provided that the subdomains on which the coefficients are defined 
satisfy \cite[Assumption 1]{BaCaTs:12}.
This assumption allows a large class of subdomains (including any bounded finite collection), but does not allow the subdomain boundaries to concentrate from below on a surface in $\Oe$, and thus the subdomains in Figure \ref{fig:2} are ruled out (see \cite[Figure 1]{BaCaTs:12}).

The argument in \cite{BaCaTs:12} equally applies to the Helmholtz equation with 
piecewise-Lipschitz $A$ 
and piecewise-$L^\infty$ $n$, proving uniqueness (and hence, by Fredholm theory, well-posedness) of the BVPs corresponding to Condition \ref{cond:simple} (see Figure \ref{fig:1}) and Condition \ref{cond:T} (see Figure \ref{fig:3}), but \emph{not} of the BVP corresponding to Condition \ref{cond:Yves} (see Figure \ref{fig:2}). 

\paragraph{Bounds on the EDP/full-space problem when $A$ and $n$ are continuous.}

The paper \cite{Bl:73} considers $A\in C^1$, $n\equiv 1$, $A\rightarrow I$ uniformly as $r\tendi$, \cite{BlKa:77} considers  $A\equiv I$, $n\in C^1$, $n\rightarrow 1$ as $r\tendi$, and both prove weighted estimates rather than cut-off resolvent estimates. 
The paper \cite{Bl:73} uses Morawetz identities with the vector field $\bx$, and hence requires $\domain_-$ to be star-shaped with respect to ball.
The paper \cite{BlKa:77} considers a more general class of domains than star-shaped \cite[Definition 2.1]{BlKa:77} by modifying vector field $\bx$ in a neighbourhood of the obstacle (see also \cite{BlKa:74}).
The results of \cite{Bl:73} are obtained under inequalities about norms of $A$ in $\domain_+$ that essentially guarantee that \eqref{eq:A2a} holds; the results of \cite{Bl:73} are obtained under a condition \cite[Equation 4.3]{Bl:73} similar to \eqref{eq:n3a}, but different since the vector field is no longer $\bx$.

The paper \cite{PeVe:99} considers the full-space problem, but with $\supp(1-n)$ not compact. They use Morawetz identities with the vector field $\bx$, and effectively prove a cut-off resolvent estimate under a condition \cite[Equation 1.8]{PeVe:99}, similar to \eqref{eq:n3a}, that limits the decrease of $n$ in the radial direction.

\paragraph{Bounds on the EDP/full-space problem when $A$ and $n$ are not continuous (transmission problems).}

The nontrapping resolvent estimate \eqref{eq:res1} for the full-space problem when $A$ and $n$ have a single common nontrapping jump across a $C^\infty$ interface with strictly positive curvature was obtained in \cite{CaPoVo:99} (following earlier work in \cite{PoVo:99a})\footnote{In fact, \cite{CaPoVo:99} also considers the analogous nontrapping transmission problem for the EDP when the impenetrable obstacle $\domain_-$ is nontrapping.}. 
With the transmission problem written as in Definition \ref{def:HTP}, the nontrapping case is when 
\beq\label{eq:ntT}
\frac{n_i}{n_o}\leq \frac{a_i}{a_o},
\eeq
and the trapping case is when this inequality is reversed. In the trapping case, when the interface $\partial \Omega_T$ is $C^\infty$ with strictly positive curvature, \cite{PoVo:99} proved that there is a sequence of resonances super-algebraically close to the real axis. The worse-case bound of exponential growth in $k$ \eqref{eq:res3} for general $C^\infty$ $\partial \Omega_T$ was then proved in \cite{Be:03a}.

In \cite{MoSp:17}, a resolvent estimate was proved for the nontrapping problem where $A$ and $n$ have one jump across a star-shaped Lipschitz $\partial \Omega_T$, and 
satisfy the slightly-stronger condition than \eqref{eq:ntT} that \eqref{eq:AnT2} holds, i.e.~that
\beqs
\frac{n_i}{n_o}\leq 1\leq \frac{a_i}{a_o};
\eeqs
the corresponding resonance-free region was then deduced (as in Corollary \ref{cor:res1} above) from the results of \cite{Vo:99}. Theorem \ref{thm:EDP2} and Corollary \ref{cor:res1} therefore contain the generalisations of the results of \cite{MoSp:17} to non-Lipschitz (but still star-shaped) $\partial \Omega_T$ \footnote{Although \cite{MoSp:17} considers more-general transmission conditions than in \eqref{eq:BVP}, namely $u_o =u_i + g_D$ and $a_o (\partial  u_o/\partial \nu) = a_i (\partial u_i/\partial \nu)  +g_N$ for arbitrary $g_D\in H^1(\partial \Omega_T)$ and $g_N\in L^2(\partial \Omega_T)$.}.

In the case when
$A\equiv I$ and $n$ has one jump on a sphere (i.e.~transmission through a penetrable ball), a priori estimates of Sobolev norms of arbitrary order on spherical surfaces in $\domain_R$ were obtained in \cite{Ca:12, CaLePa:12}, with the spherical surfaces needing to be a sufficient distance from the obstacle in the trapping case. 
Resolvent estimates when both $A$ and $n$ are discontinuous across a sphere, but $n$ is complex were obtained in \cite{NgVo:12}.

In terms of the methods used, \cite{CaPoVo:99, PoVo:99a, PoVo:99} use microlocal analysis and propagation of singularities, \cite{NgVo:12, MoSp:17} use Morawetz identities with the vector field $\bx$, and \cite{Ca:12, CaLePa:12} use separation of variables and results about the asymptotics of Bessel and Hankel functions.

\paragraph{Local energy decay of EDP for the wave equation.}
The relationship between a resolvent estimate on the time-harmonic problem 
and local-energy decay for solutions of the corresponding wave equation is well-understood in scattering theory (see, e.g, \cite[Theorem 1.1]{Vo:99}). We therefore mention briefly the following results on local-energy decay for the EDP for the wave equation \eqref{eq:2}, which are essentially equivalent to resolvent estimates on \eqref{eq:1}. Again these use Morawetz identities with the vector field $\bx$, and therefore require $\domain_-$ to be star-shaped. 
The report \cite{Li:64} considers $A\equiv I$, $n\in C^1$, and radial, and obtains local energy decay under essentially the condition \eqref{eq:n3a}. The paper
\cite{Za:66} considers $A, n\in C^1$, with $A$ scalar, and obtains local energy decay under essentially under the condition \eqref{eq:A1a}; the paper \cite{BlKa:88} is the analogue of \cite{Za:66} except with matrix-valued $A$.
For more recent extensions and generalisations of these arguments and results, see \cite{MeTa:09} and the references therein.

\section{Preliminary results and inequalities}\label{sec:3}

\subsection{Background theory of the EDP}

We now give the variational formulation of the EDP with zero Dirichlet data.
This formulation is based on Green's identity; for a proof of this identity, see, e.g., \cite[Lemma 4.3]{Mc:00}.

\ble[Green's identity and conormal derivative]\label{lem:Green}
Let $\domaingen$ be a bounded Lipschitz open set. If $A\in 
L^\infty(\domaingen,\Rea^{d\times d})$, $u\in H^1(\domaingen)$, and $\nabla\cdot (A \nabla u) \in L^2(\domaingen)$ (understood as in \eqref{eq:weak}) then there exists a uniquely defined $\phi \in H^{-1/2}(\partial \domaingen)$ such that
\beq\label{eq:Green}
\langle \phi, \gamma v\rangle_{\partial \domaingen} = \int_\domaingen \vb\, \nabla\cdot (A \gu) + \int_\domaingen 
(A\gu)\cdot \gvb
\quad \tfa v \in H^1(\domaingen),
\eeq
where $\langle \cdot,\cdot\rangle_{\partial \domaingen}$ denotes the duality pairing on $\partial \domaingen$. 
Furthermore, if $u\in H^2(\domaingen)$ and $A\in C^{0,1}(\overline{\domaingen},\Rea^{d\times d})$ then $\phi = \bnu \cdot \gamma (A \gu)$ and thus we denote $\phi$ by $\dudnuAw$.
\ele

Recall the notation that $B_R:=\{\bx: |\bx|<R\}$ and let $\Gamma_R := \partial B_R =\{\bx: |\bx|=R\}$. 
Define $T_R: H^{1/2}(\Gamma_R) \rightarrow H^{-1/2}(\Gamma_R)$ to be the Dirichlet-to-Neumann map for the equation $\Delta u+k^2 u=0$ posed in the exterior of $B_R$ with the Sommerfeld radiation condition \eqref{eq:src}. The definition of $T_R$ in terms of Hankel functions and polar coordinates (when $d=2$)/spherical polar coordinates (when $d=3$) is given in, e.g., \cite[Equations 3.5 and 3.6]{ChMo:08} \cite[\S2.6.3]{Ne:01}, \cite[Equations 3.7 and 3.10]{MeSa:10}. Two key properties of $T_R$ that we use below are contained in the following lemma.

\ble[Two key properties of $T_R$]

\

\noi (i)
\beq\label{eq:TR}
\Re \big( - \langle T_R \phi,\phi\rangle_{\Gamma_R}\big) \geq 0 \quad\tfa \phi \in H^{1/2}(\Gamma_R),
\eeq
where $\langle\cdot,\cdot\rangle_{\GR}$ denotes the duality pairing on $\GR$.

\noi (ii) There exists $C>0$, independent of $k$, such that
\beq\label{eq:TR2}
\big|\big\langle T_R(\gamma u), \gamma v\rangle_{\Gamma_R}\big\rangle\big| \leq C \N{u}_{H^1_k (\domain_R)}  \N{v}_{H^1_k (\domain_R)}  \quad \tfa u, v \in H^1(\domain_R).
\eeq
\ele

\bpf[References for the proof]
(i) is proved in \cite[Corollary 3.1]{ChMo:08} or \cite[Theorem 2.6.4]{Ne:01}, and (ii) is proved in \cite[Lemma 3.3]{MeSa:10}.
\epf

\ble[Variational formulation of EDP with $g_D\equiv 0$]\label{lem:EDP}
With $\domain_+, f, n$, and $A$ as in Definition \ref{def:EDP},
choose $R>0$ such that $\supp (I-A)$, $\supp(1-n)$, and $\supp\, f$ are all compactly contained in $B_R$. 
Let $\domain_R:= \domain_+\cap B_R$ and let 
\beq\label{eq:spaceEDP}
H_{0,D}^1(\domain_R):= \big\{ v\in H^1(\domain_R) : \gamma v=0 \ton \Gamma_D\big\}.
\eeq
The variational formulation of the EDP  of Definition \ref{def:EDP} with $g_D= 0$ is:
\beq\label{eq:EDPvar}
\text{ find } u \in H^1_{0,D}(\domain_R) \tst \quad a(u,v)=F(v) \quad \tfa v\in H^1_{0,D}(\domain_R),
\eeq
where
\beq\label{eq:EDPa}
a(u,v):= \int_{\domain_R} 
\Big((A \gu)\cdot\gvb
 - k^2 nu\vb\Big) - \langle T_R \gamma u,\gamma v\rangle_{\Gamma_R}\quad\tand\quad
F(v):= \int_{\domain_R} f\, \vb.
\eeq
\ele

\bpf[Proof of Lemma \ref{lem:EDP}]
If $u\in H^1_{\text{loc}}(\domain_+)$ satisfies the EDP of Definition \ref{def:EDP} with $g_D\equiv 0$ then its restriction to $\domain_R$ 
is in $H^1_{0,D}(\domain_R)$. Then, applying Green's identity \eqref{eq:Green} in $\domain_R$ and using the fact that $u$ satisfies $\Delta u +k^2 u=0$ in a neighbourhood of $\Gamma_R$ , we find that $u|_{\domain_R}$ satisfies the variational problem \eqref{eq:TEDPvar}.

Conversely, given $u\in H^1_{0,D}(\domain_R)$ satisfying the variational problem \eqref{eq:EDPvar}, by letting $v\in C^\infty_0(\domain_R)$ 
in \eqref{eq:EDPvar}, we see that $u$ satisfies the Helmholtz equation in $\domain_R$ in the sense of \eqref{eq:weak}.
With $\domain_R^c:= \Rea^d \setminus\overline{\domain_R}$ and $h_R:= \gamma u$ on $\Gamma_R$, we extend $u$ by setting $u|_{\domain_R^c}$ to be the solution of the Dirichlet problem for the homogeneous Helmholtz equation in $\domain_R^c$ satisfying the Sommerfeld radiation condition \eqref{eq:src}, with Dirichlet data on $\Gamma_R$ equal to $h_R$. 
Using the variational problem \eqref{eq:EDPvar}, Green's identity, and the definition of $T_R$, one can show that the Neumann traces on either side of $\Gamma_R$ of the extended function are equal.
The fact that both the Dirichlet and Neumann traces are continuous across $\Gamma_R$ then implies that the extended function satisfies $\Delta u+k^2u=0$ in a neighbourhood of $\Gamma_R$ and thus is $C^\infty$ in this neighbourhood. The extended function is therefore in $H^1_{\text{loc}}(\domain_+)$ and satisfies the EDP.
\epf

\ble[Continuity of the sesquilinear form]\label{lem:cont}
The sesquilinear form $a(\cdot,\cdot)$ of the EDP defined by \eqref{eq:EDPa} is continuous on $H^1(\domain_R)$. 
\ele

\bpf[Proof of Lemma \ref{lem:cont}]
This follows from the Cauchy-Schwarz inequality and the inequalites \eqref{eq:nlimitsEDP}, \eqref{eq:AellEDP}, \eqref{eq:TR2}.
\epf

\ble[Existence and uniqueness from an a priori bound]\label{lem:Fred}
If, under the assumption of existence, one has a bound on the solution of the EDP in terms of the data, then the solution exists and is unique.
\ele

\bpf
We first consider the case when $g_D\equiv 0$. 
The sesquilinear forms of both the EDP and TEDP with $g_D\equiv 0$ satisfy G\aa rding inequalities. Indeed, in the case of the EDP, \eqref{eq:TR} and \eqref{eq:AellEDP} imply that, for any $v\in H^1_{0,D}(\domain_R)$, 
\begin{align*}
\Re a(v,v) &\geq A_{\min} \N{\gv}^2_{L^2(\domain_R)} - k^2 \varmax \N{v}^2_{L^2(\domain_R)},\\
&\geq A_{\min} \N{v}^2_{H^1_k(\domain_R)} - k^2 (\varmax + A_{\min})\N{v}^2_{L^2(\domain_R)},
\end{align*}
where $\|v\|_{H^1_k(\domain_R)}$ is defined by \eqref{eq:1knorm}.
A bound on the solution, under the assumption of existence, shows that the solution of the BVP (if it exists) is unique. Since the sesquilinear form is continuous and satisfies a G\aa rding inequality, Fredholm theory implies that the solution of the variational problem exists and is unique; see, e.g., \cite[Theorem 2.34]{Mc:00}, \cite[\S6.2.8]{Ev:98}, \cite[Theorem 6.31]{Sp:15}.

When $g_D\not\equiv 0$, uniqueness of the solution still follows from the a priori bound. We recall that there exists $\widetilde{\gamma}: H^{1/2}(\Gamma_D)\rightarrow H^1_{\loc}(\domain_+)$ such that $\gamma \widetilde{\gamma} \phi =\phi$ for all $\phi \in H^{1/2}(\Gamma_D)$; see, e.g, \cite[Theorem 3.37]{Mc:00}. Given $u$ satisfying the EDP with $g_D\not\equiv 0$, we have that $u- \widetilde{\gamma}g_D$ satisfies the EDP with zero Dirichlet data and with suitably modified $f$ and $g$ (the compact support of $\widetilde{\gamma}g_D$ ensures that $u- \widetilde{\gamma}g_D$ satisfies the radiation condition). Existence of the solution of the EDP then follows from the case $g_D\equiv 0$. 
\epf

\

We now recall a regularity result due to Ne\v{c}as.

\begin{theorem}\textbf{\emph{(Ne\v{c}as' regularity result \cite[\S5.1.2, 5.2.1]{Ne:67}, \cite[Theorem 4.24]{Mc:00})}}\label{thm:Necas}
Let $\domaingen$ be a bounded Lipschitz open set, let $A\in C^{0,1}(\overline{\domaingen}, \Rea^{d\times d})$ with $A_{ij}=A_{ji}$ and let $A$ satisfy the inequality \eqref{eq:AellEDP} with $\domain_+$ replaced by $\domain$ and $A_{\min}>0$. 
If $u\in H^1(\domain)$ with $\nabla \cdot(A\gu)\in L^2(\domaingen)$ (understood in the sense of \eqref{eq:weak}), then $\dudnuAw\in L^2(\partial \domaingen)$ if and only if $\gamma u \in H^1(\partial \domaingen)$.
\end{theorem}

This result is proved using the analogue of the identity \eqref{A2} below with the vector field $\bx$ replaced by a general vector field.

We now apply this result to the solution of the EDP. It is convenient to introduce the following notation:
for $\domaingen$ a bounded Lipschitz open set and $A\in C^{0,1}(\overline{\domaingen}, \Rea^{d\times d})$ with $A_{ij}=A_{ji}$, define the space
\beqs
V(\domaingen):= \bigg\{ 
v\in H^1(\domaingen) :\, \, \nabla\cdot(A\gv)\in L^2(\domaingen), \,\, \dvdnuAw\in L^2(\partial \domaingen), \,\,\gamma v \in H^1(\partial \domaingen)
\bigg\}.
\eeqs
Observe that Theorem \ref{thm:Necas} implies that either of the conditions $\dvdnuAw\in L^2(\partial \domaingen)$ and $\gamma v \in H^1(\partial \domaingen)$ can be left out of the definition of $V(\domaingen)$ (since one implies the other).

\begin{corollary}[Ne\v{c}as' regularity result applied to the EDP]\label{cor:Necas}
Let $R$ be such that $\overline{\domain_-}\subset\subset \domain_R$.
If $A\in C^{0,1}(\overline{\domain_+},\Rea^{d\times d})$, $A_{ij}=A_{ji}$, and $A$ satisfies the inequality \eqref{eq:AellEDP} with $A_{\min}>0$, then the solution of the EDP of Definition \ref{def:EDP} is in $V(\domain_R)$.
\end{corollary}

\bpf[Proof of Corollary \ref{cor:Necas}]
The facts that $f\in L^2(\domain_R)$ and $u\in H^1(\domain_R)$ imply that $\nabla\cdot(A\gu)\in L^2(\domain_R)$, and
the fact that $g_D \in H^1(\Gamma_D)$ implies that $\gamma u \in H^1(\Gamma_D)$. 
Interior $H^2$-regularity of the operator
$\opL_{A,n}$ on compact subsets of $\Omega_R$ \cite[Theorem 4.16]{Mc:00}, \cite[\S6.3.1]{Ev:98} 
and the fact that $f\in L^2(\domain_+)$ with compact support imply that $\gamma u \in H^1(\Gamma_R)$ and $\dudnuAw \in L^2(\Gamma_R)$. 
Applying the Ne\v{c}as result in $\domain_R$ then gives that $\dudnuAw\in L^2(\Gamma_D)$ and we are done.
\epf

\

When $A\equiv I$, we have the following density result for the space $V(\Omega)$.

\ble{\bf (\cite[Lemmas 2 and 3]{CoDa:98})}\label{lem:density}
 If $\Omega$ is Lipschitz and $A\equiv I$, then $\DOmegabar:= \big\{ u|_\domaingen : u \in C^\infty(\Rea^d)\}$ is dense in $V(\domaingen)$.
\ele

This result, which relies on results in \cite{JeKe:81} and \cite{Da:80}, allows us to easily apply the integrated Morawetz identity, \eqref{eq:morid1int}, below to the solution of the EDP when $A\equiv I$.
Lemma \ref{lem:density} can be generalised to the case $A\in C^1$ using \cite[Proposition 7.4]{MiTa:99}. However, when $A\not\equiv I$ we use simpler techniques that bypass this density issue by first proving the a priori bounds for smooth star-shaped domains, and then using the following approximation result.

\ble\textbf{\emph{(Approximation of star-shaped domains \cite[Lemma 3.7]{ChMo:08})}}\label{lem:CWM}
Let $S:=\{ \bx \in \Rea^d : |\bx|=1\}$.
If $\Oi$ is $C^0$ and starshaped then, for every $\phi\in C^\infty_{\rm comp}(\Oe)$ and $R> \sup_{\bx \in \Gamma}|\bx|$, there exists $f\in C^\infty(S,\Rea)$ with $\min_{\hatx\in S}f(\widehat{\bx})>0$ such that
\beqs
\supp \,\phi \subset \Omega'_+ := \Big\{ s f(\widehat{\bx}) \widehat{\bx} \in \Rea^d \, :\, \widehat{\bx} \in S, \, s>1\Big\},
\eeqs
and $\Rea^d\setminus \overline{B_R} \subset \Omega'_+ \subset \Oe$.
\ele

\subsection{Preliminary inequalities}

We repeatedly use both the Cauchy-Schwarz inequality and the inequality
\beq\label{eq:Cauchy}
2 a b \leq \eps a^2 + \frac{b^2}{\eps} \quad \text{ for } a, b, \,\text{and } \eps>0
\eeq
(following \cite[\S B.2]{Ev:98} we refer to \eqref{eq:Cauchy} as the Cauchy inequality).

We now use Green's identity to bound the $L^2$ norm of $\gu$ in terms of the $L^2$ norm of $u$ (and vice versa) along with norms of the data and traces of $u$. This is a well-known method: see, e.g., 
\cite[Theorem I.1]{Mo:75} and \cite[Lemma 2.2]{Sp:14}, and note that the idea is similar to that of Caccioppoli inequalities in the Calculus of Variations.

\begin{lemma}\textbf{\emph{(Bounding the $H^1$ semi-norm of $u$ via the $L^2$ norm and $L^2$ norm of $f$)}}\label{lem:control} 
Assume there exists a solution to the EDP of Definition \ref{def:EDP}. 

\

\noi 
(i) Let $R>0$ be such that $\overline{\domain_-}$, $\supp(I-A)$, $\supp(1-n)$, and $\supp f$ are all compactly contained in $\domain_R$. Then, for all $k>0$,
\beq\label{eq:GreenEDP1}
A_{\min}\N{\gu}_{L^2(\domain_R)}^2 \leq \frac{3}{2}\varmax k^2\N{u}^2_{L^2(\domain_{R})} + \frac{1}{2\varmax k^{2}}\N{f}^2_{L^2(\domain_+)}+ \N{\gamma u}_{L^2(\Gamma_D)}\N{\dudnuA}_{L^2(\Gamma_D)},
\eeq

\noi 
(ii) Let $R>0$ be such that $\overline{\domain_-} \subset\subset \Omega_{R}$. Then, for all $k>0$,
\beq\label{eq:GreenEDP2}
\varmin k^2\N{u}^2_{L^2(\domain_R)} \leq 4 \left(A_{\max}+ \frac{6(A_{\max})^2}{\varmin k^2 }\right)\N{\gu}_{L^2(\domain_{R+1})}^2 +\frac{2}{\varmin k^{2}}\N{f}^2_{L^2(\domain_+)}+ 4\N{\gamma u}_{L^2(\Gamma_D)}\N{\dudnuA}_{L^2(\Gamma_D)}.
\eeq
\end{lemma}

\bpf
(i) Applying Green's identity \eqref{eq:Green} with  $\domain=\domain_{R}$, with $u$ the solution of the EDP, and with $v= u$, we obtain
\beq\label{eq:hot1new}
-\left\langle \dudnuA, \gamma u\right\rangle_{\Gamma_D} + \int_{\Gamma_R} \overline{u} \pdiff{u}{r}= 
\int_{\domain_{R}}
(A \gu)\cdot\overline{\gu} - k^2 n \nus - \ub f,
\eeq
where we have used the fact that $u\in C^\infty$ in a neighbourhood of $\Gamma_R$ (by elliptic regularity) to write the duality pairing on $\Gamma_R$ as an integral.
The key point now is that the inequality \eqref{eq:TR}/\eqref{2.12} involving the term on $\Gamma_R$ in \eqref{eq:hot1new} allows us to obtain an upper bound on $\int_{\domain_{R}}
(A \gu)\cdot\overline{\gu}$. Indeed, taking the real part of \eqref{eq:hot1new}, using the inequality \eqref{eq:TR}/\eqref{2.12}, and then the Cauchy-Schwarz inequality and the inequalities on $A$ and $n$ \eqref{eq:AellEDP} and \eqref{eq:nlimitsEDP}, we obtain that
\beqs
A_{\min}\N{\gu}^2_{L^2(\domain_R)} \leq  \varmax k^2\N{u}^2_{L^2(\domain_R)} + \N{u}_{L^2(\domain_R)}\N{f}_{L^2(\domain_R)} + \N{\gamma u}_{L^2(\Gamma_D)} \N{\dudnuA}_{L^2(\Gamma_D)}.
\eeqs
The result \eqref{eq:GreenEDP1} then follows from using the Cauchy inequality \eqref{eq:Cauchy} on $ \|u\|_{L^2(\domain_R)}\|f\|_{L^2(\domain_R)}$.

(ii) The sign property of the inequality \eqref{eq:TR}/\eqref{2.12} does not allow us to obtain an upper bound on $k^2 \int_{\domain_R} n\nus$ via the argument in Part (i). Instead we apply Green's identity in $\domain_{R+1}$ with $u$ the solution of the EDP and $v= \chi u$, where
$\chi(\bx):= \chi(r)$ is such that $\chi\equiv 1$ on $[0,R]$, $\chi(R+1)=0$, and $\chi(r)= F(R+1-r)$ for $r \in [R,R+1]$, where $F(t):=t^2(3-2t)$. Observe that $F$ increases from $0$ to $1$ as $t$ increases from $0$ to $1$, and thus $\chi$ decreases from $1$ to $0$ as $r$ increases from $R$ to $R+1$. 
This particular choice of $F$
is motivated by the fact that there exists an $M>0$ such that 
\beq\label{eq:hot2}
\frac{(F'(t))^2}{F(t)}\leq M \quad \tfa 0\leq t\leq 1;
\eeq
in fact, one can easily verify that this last inequality holds with $M=12$.

Applying Green's identity \eqref{eq:Green} as described above we obtain
\beq\label{eq:hot1}
-\left\langle \dudnuA, \gamma u\right\rangle_{\Gamma_D}= 
\int_{\domain_{R+1}}
\chi (A \gu)\cdot\overline{\gu} + (A\gu)\cdot(\overline{u} \,\nabla \chi) - k^2 n\chi\nus - \chi\ub f
\eeq
(where we use the convention on $\Gamma_D$ that the normal points \emph{out} of $\domain_-$ and thus \emph{into} $\domaingen_{R+1}$); observe that, since $\chi(R+1)=0$, there is no contribution from $\Gamma_R$, and thus we have avoided the issue with the sign in the inequality \eqref{eq:TR}/\eqref{2.12}.
Now, by the Cauchy-Schwarz and Cauchy inequalities
\beq\label{eq:hot1a}
\left| \int_{\domain_{R+1}} (A\gu)\cdot (\overline{u} \,\nabla\chi )\right|\leq \frac{\eps}{2} \int_{\domain_{R+1}}
\chi |u|^2
 + \frac{1}{2\eps} \int_{\domain_{R+1}} \frac{|A\gu|^2 |\nabla \chi|^2}{\chi}.
\eeq
Then, using the second inequality in \eqref{eq:AellEDP}, the inequality \eqref{eq:hot2} with $M=12$, and choosing  $\eps= n_{\min} k^2$ we obtain that
\beqs
\left| \int_{\domain_{R+1}} (A\gu)\cdot (\overline{u} \, \nabla\chi ) \right| \leq \frac{ n_{\min}k^2}{2} \int_{\domain_{R+1}} \chi |u|^2 + \frac{6(A_{\max})^2}{n_{\min} k^2}\int_{\domain_{R+1}}|\gu|^2.
\eeqs
Using this last inequality in \eqref{eq:hot1}, we find that
\begin{align*}
\frac{n_{\min}k^2}{2} \int_{\domain_{R+1}} \chi |u|^2 \leq &\left( A_{\max} + \frac{6 (A_{\max})^2}{n_{\min}k^2}\right) \N{\gu}^2_{L^2(\domain_{R+1})}
 +\int_{\domain_{R+1}}\chi \ub f +\N{\gamma u}_{L^2(\Gamma_D)} \N{\dudnuA}_{L^2(\Gamma_D)}.
\end{align*}
Using the Cauchy inequality \eqref{eq:Cauchy} again on the term $\int_{\domain_{R+1}}\chi \ub f$ with weight $\eps= n_{\min} k^2$, we obtain \eqref{eq:GreenEDP2}.
\epf

\bre[Dimensions of the factors in \eqref{eq:GreenEDP1} and \eqref{eq:GreenEDP2}]
The dimensions of the factors in front of the norms in \eqref{eq:GreenEDP1} and \eqref{eq:GreenEDP2} are as expected 
apart from
\beqs
\left(A_{\max}+ \frac{6(A_{\max})^2}{\varmin k^2}\right);
\eeqs
this expression should be non-dimensional, but instead the second term has dimension $(\text{length})^2$. This discrepancy is because there is the factor $1=((R+1)-R)^2$ (the distance between $B_{R+1}$ and $B_R$ squared) multiplying the $k^2$s, providing the missing $(\text{length})^{-2}$.
\ere

\section{Morawetz-type identities and associated results}\label{sec:4}

\subsection{Morawetz-type identity for the operator $\opL_{A,n}$}

When writing these identities, it is convenient to use the notation that $\langle \textbf{a},\textbf{b}\rangle :=\sum_{j=1}^d a_j \overline{b_j}$ for $\textbf{a}, \textbf{b}\in \Com^d$.

\ble[Morawetz-type identity]
Let $\domaingen\subset \Rea^d$, $d\geq 2$. Let $v\in C^2(\domaingen)$, $A\in C^1(\domaingen, \Rea^{d\times d})$ with $A_{ij}=A_{ji}$, $n\in C^1(\domaingen,\Rea)$, and $\alpha,\beta\in C^1(\domaingen,\Rea)$.
Let 
\beqs
\opL_{A,n} v:= \nabla\cdot (A \gv) + k^2 n\, v,
\eeqs
and let
\beq\label{eq:cM}
\cM v:= \bx\cdot \gv - \ri k \beta v + \alpha v.
\eeq
Then 
\begin{align}\nonumber 
2 \Re \big(\overline{\cM v } \,\opL_{A,n} v \big) = &\, \nabla \cdot \bigg[ 2 \Re\big(\overline{\cM v}\, A \gv\big) + \bx\Big(k^2n \nvs - 
\langle A\gv, \gv\rangle
\Big)\bigg] 
- (2\alpha -d +2) 
\langle A\gv, \gv\rangle 
\\
&\qquad+ 
\big\langle \big((\bx\cdot\nabla)A\big)\gv,\gv\big\rangle
 - \big( (d-2\alpha) n + \bx\cdot \nabla n\big) k^2 \nvs \nonumber\\
&\qquad -2 \Re\big(\ri k \vb
\langle A\gv, \nabla \beta\rangle
\big)-2\Re \big(\vb \langle A \gv, \nabla\alpha\rangle\big).
\label{eq:morid1}
\end{align}
\ele

\noi Observe that the term  
$\langle (\bx\cdot\nabla)A\gv,\gv\rangle$ equals $x_i (\partial_i A_{jl})\partial_l v \, \overline{\partial_j v}$ under the summation convention (note that, here and in the rest of the paper, we use the convention that repeated indices are summed over, but all indices are lowered).

\

\bpf
Splitting $\cM v$ up into its component parts, we see that the identity \eqref{eq:morid1} is the sum of the following three identities:
\begin{align}\nonumber
2 \Re \big(\bx\cdot \overline{\nabla v} \,\opL_{A,n} v \big) &= \nabla \cdot \bigg[ 2 \Re\big( \bx\cdot \overline{\nabla v}\,A\nabla v \big) + \bx\big(k^2n\, \nvs - \langle A\gv,\gv\rangle\big) \bigg]+ 
(d-2)\langle A\gv,\gv\rangle \\
&\qquad+ \big\langle \big((\bx\cdot\nabla)A\big)\gv,\gv\big\rangle - (dn  +\bx\cdot\nabla n)k^2 \nvs,
\label{A}
\end{align}
\beq\label{B}
2 \Re \big( \ri k \beta \overline{v} \,\opL_{A,n} v\big) = \nabla \cdot \big[ 2 \Re \big(\ri k \beta \overline{v}\, A\nabla v\big)\big]-2 \Re\big(\ri k \vb\langle A\gv, \nabla \beta\rangle\big) ,
\eeq
and
\beq\label{C}
2 \Re \big(\alpha \overline{v} \,\opL_{A,n} v \big) = \nabla \cdot \big[ 2 \Re (\alpha \overline{v}\, A\nabla v )\big] + 2\alpha k^2 n\nvs
-2\alpha\langle A\gv,\gv\rangle- 2\Re \big(\vb \langle A \gv, \nabla\alpha\rangle\big).
\eeq
To prove \eqref{B} and \eqref{C}, expand the divergences on the right-hand sides (remembering that $\alpha$ and $\beta$ are real and that $A$ is symmetric, so $\langle A \bxi,\bxi\rangle$ is real for any $\bxi\in \Com^d$).

The basic ingredient of \eqref{A} is the identity
\beq\label{basic}
(\overline{\bx\cdot \nabla v})\nabla\cdot(A\gv) = \nabla \cdot \big[(\bx\cdot \overline{\nabla v})A\nabla v\big] - 
\langle A\gv,\gv\rangle- \big((\bx\cdot\nabla)\gvb\big)\cdot A\gv.
\eeq
To prove this, expand the divergence on the right-hand side and use the fact that the second derivatives of $v$ commute.
We would like each term on the right-hand side of \eqref{basic} to either be single-signed or be the divergence of something.
To deal with the final term we use the identity
\beq\label{Melenktrick}
2 \Re \big( (\bx\cdot \nabla)\gvb \cdot A\gv\big)= \nabla \cdot \big[ \bx\langle A\gv,\gv\rangle \big] - d \langle A\gv,\gv\rangle -
\big\langle \big((\bx\cdot\nabla)A\big)\gv,\gv\big\rangle, 
\eeq
which can be proved by expanding the divergence on the right-hand side and using the fact that $A$ is symmetric.
Therefore, taking twice the real part of \eqref{basic} and using \eqref{Melenktrick} yields
\begin{align}\nonumber
2 \Re \big((\bx\cdot \gv)\nabla\cdot(A\gv)\big) =&  \nabla \cdot \bigg[ 2 \Re\big( (\bx\cdot \overline{\nabla v})\,A\nabla v \big) -\bx\langle A\gv,\gv\rangle\bigg] + (d-2)\langle A\gv,\gv\rangle \\&\qquad\qquad\qquad+ \big\langle \big((\bx\cdot\nabla)A\big)\gv,\gv\big\rangle.  \label{A2}
\end{align}
Now add $k^2$ times
\beqs
2 \Re \left[ (\bx\cdot \overline{\nabla v})nv  \right]= \nabla \cdot \big[ \bx n \nvs \big] -dn\nvs - \bx\cdot\nabla n \nvs
\eeqs
(which is the analogue of \eqref{Melenktrick} with the vector $\nabla v$ replaced by the scalar $v$ and the matrix $A$ replaced by the scalar $n$) to \eqref{A2} to obtain \eqref{A}. 
\epf

\

Our next goal is to prove an integrated version of the identity \eqref{eq:morid1} over a domain using the divergence theorem (Lemma \ref{lem:morid1int} below).
Before doing this, we need to introduce some notation regarding tangential differential operators (we mainly follow the notation in \cite[\S2]{PaWe:58}).

Let $\domaingen$ be a bounded Lipschitz open set with outward pointing unit normal vector $\bnu$. 
Recall that the surface gradient $\nT :H^1(\pD)\rightarrow L^2(\pD,\Com^{d})$ is such that, if $u$ is $C^1$ in a neighbourhood of the boundary $\pD$, then 
\beq\label{eq:nT1}
\gu = \nT u + \bnu \dudnu;
\eeq
for an explicit expression for $\nT u$ in terms of a parametrisation of the boundary, see, e.g., \cite[Page 276]{ChGrLaSp:12}. 
We now defined an operator analogous to $\nT$ when the normal derivative $\dudnuw$ is replaced by the conormal derivative $\dudnuAw$.

Given $A\in C^{0,1}(\overline{\domaingen}, \Rea^{d\times d})$ with $A_{ij}=A_{ji}$, and $u$ a $C^1$ function in a neighbourhood of $\pD$ let $\bT(u)$ be the differential operator defined by
\beq\label{eq:7}
T_j(u) := A_{ij} \left( \partial_i u - \dudnuA \frac{\nu_i}{\nu}\right)
\eeq
where 
\beq\label{eq:nu}
\nu:= A_{ij} \nu_i\nu_j.
\eeq
These definitions and the fact that $\partial u/\partial \nu_A= \nu_j A_{ij}\partial_i u$ imply that
$T_j(u) \nu_j=0$, so $\bT(u)$ is a tangential differential operator. 
Observe that \eqref{eq:7} implies that
\beq\label{eq:7a}
\gu = A^{-1}\bT(u) + \frac{\bnu}{\nu}\dudnuA,
\eeq
and so when $A\equiv I$, $\bT(u)=\nT u$ (compare \eqref{eq:7a} to \eqref{eq:nT1}); $\bT(u)$ can therefore be understood as the surface gradient in the metric induced by $A$.

Just as there exists an expression for $\nT u$  in terms of a parametrisation of $\pD$, in principle one can find an expression for $\bT(u)$ in terms of parametrisation of $\pD$ (e.g., using material in \cite[\S2]{PaWe:58}). We do not need this explicit form of $\bT(u)$  in what follows, but we use the fact that we can define $\bT$ as a mapping from $H^1(\pD)$ to $L^2(\pD,\Com^d)$.

\begin{lemma}[Integrated form of the Morawetz-type identity \eqref{eq:morid1}]
\label{lem:morid1int}
Let $\domaingen$ be a bounded Lipschitz open set, with boundary $\pD$ and outward-pointing unit normal vector $\bnu$. Let $\gamma$ denote the trace map and 
$\partial/\partial \nu_A$ the conormal derivative (defined by Lemma \ref{lem:Green}).
If $v \in \DOmegabar:= \big\{ u|_\domaingen : u \in C^\infty(\Rea^d)\}$,
$A\in C^{0,1}(\overline{\domaingen}, \Rea^{d\times d})$ with $A_{ij}=A_{ji}$, $n\in C^{0,1}(\overline{\domaingen},\Rea)$, and $\alpha,\beta\in C^{0,1}(\overline{\domaingen},\Rea)$, then
\begin{align}\nonumber 
&\int_\domain 2 \Re \big(\overline{\cM v } \,\opL_{A,n} v \big) 
+ (2\alpha -d +2) \big\langle A\gv, \gv\big\rangle -\big\langle \big((\bx\cdot\nabla)A\big)\gv,\gv\big\rangle +\big( (d-2\alpha) n + \bx\cdot \nabla n\big) k^2 \nvs
\\ &\hspace{5cm}+2 \Re\big(\ri k \vb\big\langle A\gv, \nabla \beta\big\rangle\big)+2\Re \big(\vb \big\langle A \gv, \nabla\alpha\big\rangle\big) \nonumber
\\
&=\int_{\pD}(\bx\cdot\bnu)\left(\frac{1}{\nu}\left|\dvdnuA\right|^2 - \big\langle A^{-1}\bT(\gamma v), \bT(\gamma v)\big\rangle + k^2 n |\gamma v|^2\right) \nonumber
\\ & \hspace{5cm}+ 2\Re\left(\Big( \big\langle \bx, A^{-1}\bT(\gamma v)\big\rangle + \ri k \beta \overline{\gamma v} + \alpha \overline{\gamma v}\Big)\dvdnuA\right)
\label{eq:morid1int}
\end{align}
where $\nu$ is defined by  \eqref{eq:nu} and $\bT(\cdot)$ is defined in terms of a parametrisation of the boundary as discussed above.
\ele

\begin{proof}[Proof of Lemma \ref{lem:morid1int}]
Recall that the divergence theorem $\int_\domain \nabla \cdot \textbf{F}
= \int_{\pD} \textbf{F} \cdot \bnu$
is valid when $\textbf{F} \in C^1(\overline{\domaingen}, \Com^d)$ \cite[Theorem 3.34]{Mc:00}, and thus for $\bF\in H^1(\domaingen,\Com^d)$ by the density of $C^1(\overline{\domaingen})$ in $H^1(\domain)$ \cite[Theorem 3.29]{Mc:00} and the continuity of trace operator from $H^1(\domain)$ to $H^{1/2}(\pD)$ \cite[Theorem 3.37]{Mc:00}.
Recall that the product of an $H^1(\Omega)$ function and a $C^{0,1}(\overline{\Omega})$ function is in $H^1(\Omega)$,
and the usual product rule for differentiation holds for such functions. This result implies that
\beqs
\bF=2 \Re\big(\overline{\cM v}\, A \gv\big) + \bx\big(k^2n \nvs - \big\langle A\gv, \gv\big\rangle\big)
\eeqs
is in $H^1(D, \Com^d)$ 
and then \eqref{eq:morid1} implies that $\nabla\cdot\bF$ is given by the integrand on the left-hand side of \eqref{eq:morid1int}.
To complete the proof, we need to show that $\bF\cdot\bnu$ equals the integrand on the right-hand side of \eqref{eq:morid1int}.
Since $v\in \DOmegabar$ and $A$ is Lipschitz, $(A\cdot \gv)\cdot \bnu= \dvdnuAw$, and we see that we only need to show that
\begin{align}\nonumber
2 \Re\left((\bx\cdot\gvb)\dvdnuA\right) - (\bx\cdot\bnu)\big\langle A\gv, \gv\big\rangle= &(\bx\cdot\bnu)\left(\frac{1}{\nu}\left|\dvdnuA\right|^2 - \big\langle A^{-1}\bT(v), \bT(v)\big\rangle\right)\\
&\qquad\qquad
+ 2\Re\left(\big\langle \bx, A^{-1}\bT(v)\big\rangle\dvdnuA\right),\label{eq:9}
\end{align}
with $\bT(\cdot)$ given by \eqref{eq:7}.
By multiplying \eqref{eq:7} by $\overline{\partial_j v}$, we find that 
\beq\label{eq:10}
A_{ij} \partial_i v \overline{\partial_j v} = T_j(v) \overline{\partial_j v} + \left| \dvdnuA\right|^2 \frac{1}{\nu},
\eeq
and by taking the complex conjugate of \eqref{eq:7} and multiplying by $(A^{-1})_{jk}T_k(v)$ we find that
\begin{align}\nonumber
(A^{-1})_{jk}T_k(v)\overline{T_j(v)}&= (A^{-1})_{jk}A_{ij}T_k(v)\left( \overline{\partial_i v} - \overline{\dvdnuA}\frac{\nu_i}{\nu}\right),\\
&=T_k(v)\overline{\partial_k v}\label{eq:11}
\end{align}
(since $T_k(v)\nu_k=0$).
Putting \eqref{eq:10} and \eqref{eq:11} together we get
\beq\label{eq:12}
A_{ij}\partial_i v \overline{\partial_j v} = (A^{-1})_{jk}T_k(v) \overline{T_j(v)} + \left| \dvdnuA\right|^2 \frac{1}{\nu}.
\eeq
Next we use \eqref{eq:7a} to show that 
\begin{align}
x_i \overline{\partial_i v}\dvdnuA &=\left( x_i (A^{-1})_{ji} \overline{T_j(u)} + \overline{\dvdnuA}x_i \frac{\nu_i}{\nu}\right) \dvdnuA=  x_i (A^{-1})_{ji} \overline{T_j(v)} \dvdnuA+\frac{ (\bx\cdot\bnu)}{\nu} \left|\dvdnuA\right|^2 \label{eq:13}.
\end{align}
Using \eqref{eq:12} and \eqref{eq:13} in the left-hand side of \eqref{eq:9} (and recalling that $A$, and hence also $A^{-1}$, is symmetric) we see that \eqref{eq:9} holds and the proof is complete.
\epf

\bre[Bibliographic remarks on Morawetz-type identities for $\opL_{A,n}$]
The idea of multiplying second-order PDEs with first-order expressions has been used by many authors; multiplying $\Delta v$ by a derivative of $v$ goes back to Rellich \cite{Re:40, Re:43}, and multiplying $\nabla\cdot(A\gv)$ by a derivative of $v$ goes back to H\"ormander \cite{Ho:53} and Payne and Weinberger \cite{PaWe:58} (e.g., the identity \eqref{eq:morid1} with $n$, $\alpha$, and $\beta$ all equal zero appears as \cite[Equation 2.4]{PaWe:58}).

In the context of the Helmholtz equation, the identity \eqref{eq:morid1} with $A\equiv I$, $n\equiv 1$, $\bx$ replaced by a general vector field, and $\alpha$ and $\beta$ replaced by general scalar fields was the heart of Morawetz's paper \cite{Mo:75} (following the earlier work by Morawetz and Ludwig \cite{MoLu:68}-- see Lemma \ref{lem:ML} below).
The identity with $A$ variable and $n\equiv 1$ was used by Bloom in \cite{Bl:73}, and the identity \eqref{eq:morid1} with $A\equiv I$ and variable $n$ was used in Bloom and Kazarinoff in \cite{BlKa:77}. 
\ere

\subsection{The Morawetz-Ludwig identity for the operator $\opL v:= (\Delta +k^2 )v$}

\ble\textbf{\emph{(Morawetz-Ludwig identity, \cite[Equation 1.2]{MoLu:68})}}\label{lem:ML}
 Let $v \in C^2(\domaingen)$ for some $\domaingen\subset \Rea^d$, $d\geq 2$. Let $\opL v:=(\Delta +k^2)v$ and let
 \beqs
\cM_\alpha v :=	 r\left(v_r -\ri k v + \frac{\alpha}{r}v\right),
\eeqs
where $\alpha \in \Rea$ and $v_r=\bx\cdot \gv/r$. Then 
\bal\nonumber
2\Re( \overline{\cM_\alpha v} \opL v) =  &\,\nabla \cdot \bigg[2\Re \left(\overline{\cM_{\alpha} v} \gv\right)+ \left(k^2\nvs - \ngvs \right)\bx\bigg] \\&+ \big(2\alpha -(d-1)\big)\big(k^2 \nvs - \ngvs\big) - \big(\ngvs -\nvrs\big)- \big| v_r -\ri k v\big|^2.
\label{eq:ml2d}
\end{align}
\ele

The Morawetz-Ludwig identity is a particular example of the identity \eqref{eq:morid1} with $A\equiv I, n\equiv 1, \beta=r,$ and $\alpha$ a constant, and some further manipulation of the non-divergence terms (using the fact that $\bx= \beta\nabla\beta$). For a proof, see \cite{MoLu:68}, \cite[Proof of Lemma 2.2]{SpChGrSm:11}, or \cite[Proof of Lemma 2.3]{SpKaSm:15}.

The Morawetz-Ludwig identity \eqref{eq:ml2d} has two key properties. With this identity written as $\nabla\cdot\bQ(v)=P(v)$, the key properties are:
\ben
\item 
If $u$ is a solution of $\opL u=0$ in $\Rea^d\setminus \overline{B_R}$ satisfying the Sommerfeld radiation condition \eqref{eq:src}, then
\beq\label{eq:N0}
\int_\GR \bQ(u) \cdot\hatx
 \tendo \quad \tas R\tendi,
\eeq
(independent of the value of $\alpha$ in the multiplier $\cM_\alpha u$); see \cite[Proof of Lemma 5]{MoLu:68}, \cite[Lemma 2.4]{SpChGrSm:11}. 
\item 
If $\opL u=0$ and $2\alpha= (d-1)$,  then 
\beq\label{eq:N00}
P(u) \geq 0.
\eeq
\een
Since $I-A$ and $1-n$ both have compact support, $\opL_{A,n }=\opL$ outside a sufficiently large ball.
The two properties above of the Morawetz-Ludwig identity  mean that if the multiplier that we use on the operator $\opL_{A,n}$ is equal to $\cM_{(d-1)/2}$ outside a large ball, then there is no contribution from infinity. A convenient way to encode this information is in the following lemma (which first appeared in \cite[Lemma 2.1]{ChMo:08}).

\ble[First inequality on $\GammaR$ used to deal with the contribution from infinity]\label{lem:2.1}
Let $u$ be a solution of the homogeneous Helmholtz equation in $\Rea^d\setminus \overline{B_{R_0}}$, $d=2,3$, for some $R_0>0$, satisfying the Sommerfeld radiation condition \eqref{eq:src}.
Let $\alpha,\beta\in \Rea$ with 
$\beta\geq R$ and $2\alpha\geq d-1$. Then, for $R>R_0$, 
\beq\label{eq:2.1}
\int_{\GammaR} R\left( \left|\pdiff{u}{r}\right|^2 - |\nabla_S u|^2 + k^2 |u|^2\right)   - 2 k \beta\, \Im \int_{\GammaR} \bar{u} \pdiff{u}{r} + 2\alpha\Re \int_{\GammaR}\bar{u}\pdiff{u}{r} \leq 0,
\eeq
where $\nabla_S$ is the surface gradient on $r=R$.
\ele

We have purposely chosen the two constants in \eqref{eq:2.1} to be $\beta$ and $\alpha$, emphasising the fact that the left-hand side of \eqref{eq:2.1} is 
$\int_\GR \bQ(u) \cdot\hatx$ with $\bQ(u)$ arising from the multiplier $\cM u:= \bx\cdot \gu - \ri k\beta u +\alpha u$.

\

\bpf[Proof of Lemma \ref{lem:2.1}] 
We first show that it is sufficient to prove \eqref{eq:2.1} with $\beta=R$ and $2\alpha=d-1$. Indeed, the two inequalities 
\beq\label{2.12}
\Re \int_{\GammaR} \bar{u}\pdiff{u}{r} \,\rd s \leq 0\quad\tand\quad \Im \int_{\GammaR} \bar{u}\pdiff{u}{r} \,\rd s \geq 0 
\eeq
can be proved using the explicit expression for the solution of the Helmholtz equation in the exterior of a ball (i.e. an expansion in either trigonometric polynomials, for $d=2$, or spherical harmonics, for $d=3$, with coefficients given in terms of Bessel and Hankel functions) and then proving bounds on the particular combinations of Bessel and Hankel functions; see \cite[Theorem 2.6.4, p.97]{Ne:01} or \cite[Lemma 2.1]{ChMo:08} (observe that the first inequality in \eqref{2.12} is equivalent to the inequality \eqref{eq:TR}).
If we establish \eqref{eq:2.1} with $\beta=R$ and $2\alpha=d-1$, then the two inequalities in \eqref{2.12} give the result for $\beta\geq R$ and $2\alpha\geq d-1$.

We now integrate \eqref{eq:ml2d} with $v=u$ and $2\alpha =d-1$ over $B_{R_1}\setminus B_R$, use the divergence theorem, and then let $R_1\tendi$
(note that using the divergence theorem is allowed since $u$ is $C^\infty$ by elliptic regularity).
The first key property of the Morawetz-Ludwig identity \eqref{eq:N0} implies that the surface integral on $\vert \bx \vert ={R_1}$ tends to zero as $R_1\tendi$ \cite[Lemma 2.4]{SpChGrSm:11}. Then, recalling that $\hatx:= \bx/r$, and
using the decomposition $\gv = \nabla_S v + \widehat{\bx}v_r$ on the integral over $\Gamma_R$ (or equivalently using the right-hand side of \eqref{eq:morid1int} with $A\equiv I, n\equiv 1, \beta= r, 2\alpha=d-1, \pD=\Gamma_R$, and $\bnu=\widehat{\bx}$), we obtain that
\begin{align*}
\int_\GR \bQ(u) \cdot\hatx &=
\int_{\GammaR} R\left( \left|\pdiff{u}{r}\right|^2  - |\nabla_S u|^2+ k^2 |u|^2\right)   - 2 k R\, \Im \int_{\GammaR} \bar{u} \pdiff{u}{r} + (d-1)\Re \int_{\GammaR}\bar{u}\pdiff{u}{r}  \nonumber \\
& =  -\int_{\Rea^d\setminus B_R}\left(\big(\ngus -\nurs\big)+ \left|
 u_r - \ri k u
 \right|^2\right)\leq 0
 \end{align*}
 (where this last inequality is the second key property \eqref{eq:N00}); i.e. we have established \eqref{eq:2.1} with $\beta=R$ and $2\alpha=d-1$ and we are done.
 \epf

\

To prove the bound \eqref{eq:bound3b} in Part (ii) of Theorem \ref{thm:EDP3}, we actually need the inequality \eqref{eq:2.1} with $2\alpha=d-2$. We now prove this result using Lemma \ref{lem:2.1}  and the multiplier 
$\cM u:= \bx\cdot \gu - \ri k\beta u +\alpha u$ with $\alpha$ variable.
The price we pay is that $\beta$ must be larger ($\geq 2R$ instead of $\geq R$), and the inequality only holds for $kR$ sufficiently large.

\ble[Second inequality on $\GammaR$ used to deal with the contribution from infinity]\label{lem:2.1mod}
Let $u$ be a solution of the homogeneous Helmholtz equation in $\Rea^d\setminus B_{R_0}$, for some $R_0>0$, satisfying the Sommerfeld radiation condition. 
Let $\alpha,\beta\in \Rea$ with 
$\beta\geq 2R$ and $2\alpha\geq d-2$. If $kR \geq \sqrt{3/8}$, then \eqref{eq:2.1} holds.
\ele

\bpf[Proof of Lemma \ref{lem:2.1mod}]
From the two inequalities \eqref{2.12}, it is sufficient to prove the result with $\beta=2R$ and $2\alpha=d-2$; i.e. to prove that 
\beq\label{eq:a0}
\int_{\GammaR} R\left( \left|\pdiff{u}{r}\right|^2 - |\nabla_S u|^2+ k^2 |u|^2 \right)   - 2 k (2R) \Im \int_{\GammaR} \bar{u} \pdiff{u}{r} + (d-2)\Re \int_{\GammaR}\bar{u}\pdiff{u}{r} \leq 0.
\eeq
The overall idea of the proof is the following: we apply the identity arising from the multiplier
\beq\label{eq:2.1mult}
\cM u :=\bx\cdot \gu - \ri (2R)k u + \frac{d-2+\chi}{2}u
\eeq
in $B_{2R}\setminus B_R$, where $\chi=\chi(r)$ is a $C^1$ function such that $\chi(R)=0$ and $\chi(2R)=1$; i.e. $\chi$ brings $2\alpha$ up from $d-2$ on $\Gamma_R$ to $d-1$ on $\Gamma_{2R}$. The contribution from $\Gamma_{2R}$ is dealt with using Lemma \ref{lem:2.1}, i.e. the Morawetz-Ludwig identity (and this is why we chose $\beta=2R$ in the multiplier), and the contribution from $B_{2R}\setminus B_R$ can be controlled by choosing $\chi$ appropriately and making $R$ sufficiently large.

We chose $R>R_0$ (so that $\opL u=0$ in $\Rea^d\setminus \overline{B_R}$) and write the identity arising from the multiplier \eqref{eq:2.1mult} (i.e. \eqref{eq:morid1} with $A\equiv I, n\equiv 1, 2\alpha = d-2 + \chi$, and $\beta=2R$) as $\nabla\cdot\bQ(u) = P(u)$. By applying the divergence theorem (justified since $u\in C^\infty$, by elliptic regularity, and $\chi \in C^1[R,2R]$) we have
\beq\label{eq:a1}
\int_{\Gamma_R}\bQ(u)\cdot\widehat{\bx} + \int_{B_{2R}\setminus B_R} P(u) =  \int_{\Gamma_{2R}}\bQ(u)\cdot\widehat{\bx}.
\eeq
The right-hand side of \eqref{eq:a1} is given by the right-hand side of \eqref{eq:morid1int} with $A,n,\alpha$, and $\beta$ as above, $\pD=\Gamma_{2R}$ and $\bnu=\widehat{\bx}$, and thus equals
\beqs
\int_{\Gamma_{2R}} 2R\left( \left|\pdiff{u}{r}\right|^2  - |\nabla_S u|^2+ k^2 |u|^2\right)   - 2 k (2R) \Im \int_{\Gamma_{2R}} \bar{u} \pdiff{u}{r} + (d-1)\Re \int_{\Gamma_{2R}}\bar{u}\pdiff{u}{r},
\eeqs
which is $\leq 0$ by \eqref{eq:2.1}.
We therefore have that
\beq\label{eq:a1a}
\int_{\Gamma_R}\bQ(u)\cdot\widehat{\bx} + \int_{B_{2R}\setminus B_R} P(u) \leq 0.
\eeq
The definition of the multiplier \eqref{eq:2.1mult} and the explicit expression for $\bQ(u)\cdot\widehat{\bx}$ contained in the right-hand side of \eqref{eq:morid1int} imply that $\int_{\Gamma_R}\bQ(u)\cdot\widehat{\bx}$ equals the left-hand side of \eqref{eq:a0}. Therefore, from \eqref{eq:a1a} we see that to prove \eqref{eq:a0} we only need to show that, once $R$ is sufficiently large,
\beq\label{eq:a2}
 \int_{B_{2R}\setminus B_R} P(u)\geq 0.
\eeq
Using the definition of the multiplier \eqref{eq:2.1mult}, we find 
\begin{align}
\int_{B_{2R}\setminus B_R} P(u) 
= \int_{B_{2R}\setminus B_R} \chi \ngus + \big(2-\chi\big) k^2 \nus +   \Re\big(\ub\, \gu\cdot\nabla\chi\big).\label{eq:a3}
\end{align}
By the Cauchy-Schwarz and Cauchy inequalities, and the fact that $|\nabla \chi|=|\chi'|$, we have
\begin{align*}
\left|
\int_{B_{2R}\setminus B_R}2 \Re\big(\ub \,\gu\cdot\nabla\chi\big)
\right| \leq  \int_{B_{2R}\setminus B_R} |u| |\gu||\chi'|&\leq   \int_{B_{2R}\setminus B_R}\chi^{1/2} |\gu| \, |u|\frac{|\chi'|}{\chi^{1/2}},\\
&\leq\frac{\eps}{2} \int_{B_{2R}\setminus B_R}\chi \ngus + \frac{1}{2\eps}\int_{B_{2R}\setminus B_R}\frac{|\chi'|^2}{\chi}\nus,
\end{align*}
for any $\eps>0$. Letting $\eps=2$ and using the resulting inequality in \eqref{eq:a3}, we have
\beqs
\int_{B_{2R}\setminus B_R} P(u)  \geq  \int_{B_{2R}\setminus B_R}\left( (2-\chi) k^2 -  \frac{|\chi'|^2}{4\chi}\right)\nus.
\eeqs
The only requirements on $\chi$ we have imposed so far are that $\chi \in C^1[R,2R]$ with $\chi(R)=0$ and $\chi(2R)=1$.
We now assume that $0\leq \chi(r)\leq 1$ for $r\in [R,2R]$.
To obtain the inequality \eqref{eq:a2} (and hence \eqref{eq:a0}, and hence the result) we need to choose $\chi$ such that 
\beq\label{eq:a3a}
\frac{|\chi'(r)|^2}{\chi(r)\big(2-\chi(r)\big)} \leq 
4k^2 \quad\tfor R\leq r\leq 2R.
\eeq
We let $\chi(r):= F((r-R)/2R)$ where $F(t) = t^2(3-2t)$. Observe that $F$ increases from $0$ to $1$ as $t$ increases from $0$ to $1$, and thus 
 $\chi$ increases from $0$ to $1$ as $r$ increases from $R$ to $2R$.
As in the proof of Lemma \ref{lem:control}, this choice of $F$ is motivated by the fact that there exists an $M>0$ such that \eqref{eq:hot2} holds (and thus the left-hand side of \eqref{eq:a3a} will be bounded).
One can easily verify that 
\beqs
\frac{(F'(t))^2}{F(t)\big(2-F(t)\big)}\leq 6 \quad\tfor 0\leq t\leq 1,
\eeqs
and then the chain rule implies that \eqref{eq:a3a} is satisfied (and hence \eqref{eq:2.1} holds) if $kR \geq \sqrt{3/8}$.
\epf

\section{Proofs of the main results in \S\ref{sec:EDPresults1}}
\label{sec:5}

\bpf[Proof of Theorem \ref{thm:EDP1} under the additional assumption that $\Omega_-$ is $C^{1,1}$]
By Lemma \ref{lem:Fred}, 
it is sufficient to prove the bound \eqref{eq:bound1} under the assumption of existence.

Since $\Omega$ is $C^{1,1}$ and $\gamma u =0$ on $\Gamma_D$, elliptic regularity implies that $u\in H^2(\Omega_R)$ for every $R>0$; see, e.g., \cite[Theorem 4.18]{Mc:00}. Since $\cD(\overline{\Omega_R})$ is dense in $H^2(\Omega_R)$ \cite[Page 77]{Mc:00} and the integrated identity \eqref{eq:morid1int} is continuous in $v$ with respect to the topology of $H^2(\domaingen)$ (this requires the inequality $\|\bT(\gamma v)\|_{L^2(\pD)}\leq C \|\gamma v\|_{H^1(\pD)}$, which holds since $\bT(u)$ is tangential; see \cite[Lemma 4.23(i)]{Mc:00}),  \eqref{eq:morid1int} holds for $\Omega= \Omega_R$ and $v=u$. Then we let 
$\beta=R$ and $2\alpha=d-1$, we use the fact that $\gamma u=0$ on $\Gamma_D$ to simplify the terms on $\Gamma_D$, and we use the facts that $A\equiv I$ and $n\equiv 1$ in a neighbourhood of $\Gamma_R$ to simplify the terms on $\Gamma_R$. The result is
\begin{align}\nonumber 
&\int_{\domain_R}  \big\langle \big(A-(\bx\cdot\nabla)A\big)\gu, \gu\big\rangle+\big( n + \bx\cdot \nabla n\big) k^2 \nus
+ \int_{\Gamma_D} \frac{\bx\cdot\bnu}{\nu}\left|\dudnuA\right|^2 \\
& = 
2 \Re \int_{\domain_R}\overline{\cM u } \,f+\int_{\GammaR} R\left( \left|\pdiff{u}{r}\right|^2 - |\nabla_S u|^2+ k^2 |u|^2 \right)   - 2 k R\, \Im \int_{\GammaR} \bar{u} \pdiff{u}{r} + (d-1)\Re \int_{\GammaR}\bar{u}\pdiff{u}{r},  \label{eq:m1}
\end{align}
where $\nabla_S$ is the surface gradient on $\Gamma_R$, $\nu$ is defined by \eqref{eq:nu} above, and the normal vector $\bnu$ on $\Gamma_D$ is taken to point \emph{out} of $\domain_-$ (and thus \emph{into} $\domain_+$).

Our choices of $\beta$ and $2\alpha$ ensure that the combination of the integrals over $\Gamma_R$ in \eqref{eq:m1} is $\leq 0$ by Lemma \ref{lem:2.1}, and the fact that $\domain_-$ is star-shaped (along with Part (i) of Lemma \ref{lem:star}) implies that the integral over $\Gamma_D$ is $\geq 0$. Using the conditions on $A$ and $n$ \eqref{eq:A1} and \eqref{eq:n1}, the definition of the multiplier $\cM u$ \eqref{eq:cM}, and the Cauchy-Schwarz and Cauchy inequalities, we find that
\begin{align*}\nonumber
\mu_1 &\N{\gu}^2_{L^2(\domain_R)} + \mu_2 k^2 \N{u}^2_{L^2(\domain_R)} \leq 2 \Re \int_{\domain_R} \left( \bx\cdot \overline{\gu} + \ri k R \ub + \frac{d-1}{2} \ub\right) f,\\ \nonumber
&\quad\qquad\leq 2 \left( R \N{\gu}_{L^2(\domain_R)} + \left( R + \frac{d-1}{2k}\right)k \N{u}_{L^2(\domain_R)}\right) \N{f}_{L^2(\domain_+)},\\
&\quad\qquad\leq \eps_1 \N{\gu}^2_{L^2(\domain_R)} + \eps_2 k^2 \N{u}^2_{L^2(\domain_R)} + \left(\frac{1}{\eps_1} R^2 + \frac{1}{\eps_2}\left(R + \frac{d-1}{2k}\right)^2\right)\N{f}^2_{L^2(\domain_+)}
\end{align*}
for any $\eps_1, \eps_2>0$. Choosing $\eps_1=\mu_1/2$ and $\eps_2=\mu_2/2$ we get the result \eqref{eq:bound1}.
\epf

\

To reduce the smoothness of $\Omega$ from $C^{1,1}$ to Lipschitz (and hence prove Theorem \ref{thm:EDP1}), we need the following lemma; this lemma is also the main ingredient in the proof of Corollary \ref{cor:H1}.

\begin{lemma}[Bound for $f\in L^2$ implies bound for $F\in (H^1)'$]
\label{lem:H1}
Assume that $\Omega_-, A$, and $n$ are such that, given $f\in L^2(\domain_+)$ with compact support, the solution of the EDP of Definition \ref{def:EDP} with $g_D\equiv 0$ exists, is unique, and satisfies the bound 
\beq\label{eq:bound_lem1}
\mu_1\N{\gu}^2_{L^2(\Omega_R)}+ \mu_2 k^2 \N{u}^2_{L^2(\Omega_R)}\leq C(A,n,\Omega_-,R, k) \N{f}^2_{L^2(\Omega_+)} \quad \tfa k>0,
\eeq
and for some $\mu_1, \mu_2, C>0$.
Given $F\in (H^1_{0,D}(\domain_R))'$, let $\tilde{u}$ satisfy the variational formulation of the EDP \eqref{eq:EDPvar} with $g_D\equiv 0$. 
Then $\tilde{u}$ exists, is unique, and satisfies the bound \eqref{eq:H1bound} with $u$ replaced by $\tilde{u}$ and $C_1$ replaced by $C$.
\end{lemma}

\bpf[Proof of Lemma \ref{lem:H1}]
Existence and uniqueness of $\tilde{u}$ follows from 
Lemma \ref{lem:Fred}.
Define 
\beqs
a_+(u,v):= \int_{\domain_R} \Big(
(A\gu)\cdot\gvb
+k^2 \var u \vb\Big)-\big\langle T_R \gamma u, \gamma v\big\rangle_{\GR},
\eeqs
define $u_+ \in \HoDk$ to be the solution of the variational problem $a_+(u_+,v)=F(v)$ for all $v\in \HoDk$, 
and define $w\in \HoDk$ to be the solution of the variational problem 
$a(w,v)=2k^2 \int_{\domain_R} n u_+\vb$ for all $v\in \HoDk$. 
The whole point of these definitions is that $\tilde{u}=u_++w$, $a_+(\cdot,\cdot)$ is coercive, and $w$ satisfies an EDP with data in $L^2(\Omega_R)$. Indeed,
\beqs
\Re a_+(v,v) \geq A_{\min}\N{\gv}^2_{L^2(\domain_R)} + k^2 \varmin\N{v}^2_{L^2(\domain_R)} \geq \min(A_{\min},\varmin)\N{v}^2_{\HoDkk},
\eeqs
and so by the Lax--Milgram theorem 
\beq\label{eq:LM1}
\N{u_+}_{\HoDkk}\leq \big(\min(A_{\min},\varmin)\big)^{-1}\N{F}_{(\HoDk)'}.
\eeq
Combining the bounds \eqref{eq:bound_lem1} and \eqref{eq:LM1}, we have
\begin{align*}
\N{w}_\HoDkk^2 \leq \frac{C}{\min(\mu_1,\mu_2)} (2k^2)^2 \N{n u_+}^2_{L^2(\Omega_R)} &\leq \frac{C}{\min(\mu_1,\mu_2)} 4 (\varmax)^2 k^2 \N{u_+}^2_\HoDkk \\
& \leq \frac{4 C (\varmax)^2 }{\min(\mu_1,\mu_2) (\min(A_{\min},\varmin))^2}k^2\N{F}^2_{(\HoDkk)'}.
\end{align*}
The bound \eqref{eq:H1bound} with $u$ replaced by $\tilde{u}$ and $C_1$ replaced by $C$ follows from using this last bound along with \eqref{eq:LM1} and fact that $\tilde{u}=u_+ +w$.
\epf

\

\bpf[Proof of Theorem \ref{thm:EDP1}]
By Lemma \ref{lem:Fred}, 
we only need to show that, 
under the assumption of existence of a solution, the bound \eqref{eq:bound1} holds.

Let $u$ solve the variational problem \eqref{eq:EDPvar} with $\Oi$ Lipschitz and star-shaped (actually, we only need $\Oi$ to be $C^0$ and star-shaped to apply Lemma \ref{lem:CWM}). 
Let $W:= \{ \phi|_{\Omega_R} : \phi \in C_{0}^\infty(\Oe)\}$. In this proof only, we use the notation that $\|u\|^2_\mu := \mu_1 \|\gu\|_{L^2(\Omega_R)}^2+ \mu_2 k^2 \|u\|^2_{L^2(\Omega_R)}$. Since $\|\cdot\|_\mu$ is equivalent to $\|\cdot\|_{H^1_{k}(\OR)}$, and $W$ is dense in $H^1_{0,D}(\OR)$ \cite[Page 77]{Mc:00}, given $\eps>0$, there exists a $u_\eps \in W$ such that $\|u-u_\eps\|_\mu <\eps$. Let $\phi_\eps \in C_{0}^\infty(\Oe)$ be such that $u_\eps= \phi_\eps$ on $\Omega_R$. 
Applying Lemma \ref{lem:CWM} to $\phi_\eps$, we have that there exists $f_\eps \in C^\infty(S,\Rea)$ with $\min_{\hatx \in S}f(\hatx)>0$, $\supp\, \phi_\eps \subset \Omega^\eps_+$ and $\Rea^d\setminus \overline{B_R} \subset\Omega^\eps_+ \subset \Omega_+$, where
\beqs
\Omega^\eps_+ :=\Big\{ s f_\epsilon(\hatx)\hatx \in \Rea^d \, :\, \hatx \in S, \, s>1\Big\}
\eeqs
is $C^\infty$ and star-shaped. 
Let $\Omega^\eps_R:= \Omega^\eps_+ \cap B_R$; with this definition $u_\eps \in H_{0,D}^1(\Omega^\eps_R)$.
Given $v_\eps \in H^1_{0,D}(\Omega^\eps_R)$, let $v$ denote its extension by zero from $\Omega^\eps_R$ to $\Omega_R$; then $v\in H_{0,D}^1(\Omega_R)$ and we can therefore regard $H_{0,D}^1(\Omega^\eps_R)$ as a subspace of $H^1_{0,D}(\Omega_R)$. 

Let $a_\eps(u,v)$ be defined by \eqref{eq:EDPa} with $\Omega_R$ replaced by $\Omega^\eps_R$ and 
observe that, for all $v\in H_{0,D}^1(\Omega_R^\eps)$,
\beqs
a_\eps(u_\eps,v) = a(u_\eps,v) = F(v) -a(u-u_\eps,v).
\eeqs
Let $u', u''\in H_{0,D}^1(\Omega^\eps_R)$ be defined as the solutions of the variational problems
\beqs
a_\eps(u',v)= F(v) \quad \tand \quad a_\eps(u'', v) = -a(u-u_\eps, v) \quad \tfa v\in H^1_{0,D}(\Omega_R^\eps).
\eeqs
The key features of these definitions are that: (i) $u_\eps= u' + u''$, (ii) both $u'$ and $u''$ satisfy EDPs with the obstacle $C^\infty$ and star-shaped, and with $A$ and $n$ satisfying the analogue of Condition \ref{cond:1} with $\Oe$ replaced by its subset $\Omega_+^\eps$, (iii) $u'$ satisfies the EDP with the same right-hand side as $u$, and (iv)
$u''$ satisfies the EDP with right-hand side depending on $u-u_\eps$, which can be made arbitrarily small.

Using (a) our proof of Theorem \ref{thm:EDP1} with $C^{1,1}$ and star-shaped $\Oi$, (b) Lemma \ref{lem:H1} and (c) the definition of $u_\eps$, we have that both $u'$ and $u''$ exist, are unique, and satisfy 
\beq\label{eq:u'u''bounds}
\N{u'}_\mu \leq \sqrt{C_1}\N{f}_{L^2(\Omega_+)} \quad \tand \quad \N{u''}_{H^1_k(\Omega_R)} \leq C\, \widetilde{C_c} \N{u-u_\eps}_\mu \leq C\, \widetilde{C_c} \eps.
\eeq
respectively, where $C_1$ is given by \eqref{eq:C1}, $C$ is the constant on the right-hand side of \eqref{eq:H1bound}, and  $\widetilde{C_c}$ is the continuity constant of the sesquilinear form $a(\cdot,\cdot)$ in the norm $\|\cdot\|_\mu$; the exact forms of $C$ and $\widetilde{C_c}$ are not important in what follows, the key point is that they are both independent of $\eps$. Since 
\beqs
\N{u}_\mu \leq \N{u_\eps}_\mu + \eps \leq \N{u'}_\mu + \N{u''}_\mu + \eps,
\eeqs
the result that $\|u\|_\mu \leq \sqrt{C_1}\N{f}_{L^2(\Omega_+)}$ follow from the bounds \eqref{eq:u'u''bounds}, the fact that $\|\cdot\|_\mu$ is equivalent to $\|\cdot\|_{H^1_k(\Omega_R)}$, and the fact that $\eps$ was arbitrary; the proof is therefore complete.
\epf

\

\bpf[Proof of Theorem \ref{thm:EDP2}]
By Lemma \ref{lem:Fred}, we only need to show that, 
under the assumption of existence of a solution, the bound \eqref{eq:bound1} holds.

The first step is to approximate $A$ by $A_\delta \in C^\infty(\domain_+, \Rea^{d\times d})$ 
and $n$ by $n_\delta \in C^\infty(\domain_+, \Rea)$
such that 
\beqs
\N{A_\delta-A}_{L^2(\domain_R)}\,\tand\,\N{n_\delta -n}_{L^2(\domain_R)}  \tendo  \quad\tas\quad \delta \tendo,
\eeqs
$A_\delta$ satisfies \eqref{eq:A1} with $\mu_1=A_{\min}$, and $n_\delta$ satisfies \eqref{eq:n1} with $\mu_2=\varmin$.
We first show how to use mollifiers to construct such an $n_\delta$ (the construction of $A_\delta$ is analogous); to do this when $\Oi$ is nonempty, we need to define $n$ on a slightly larger set than $\domain_+$. 
For $\bx \in \overline{\domain_-}$, we let $\pi(\bx):=0$, so that $n(\bx)= \varmin$ for $\bx \in \overline{\domain_-}$; it follows immediately that the extended $\pi(\bx)$ is monotonically non-decreasing in the radial direction on all of $\Rea^d$.
Let $\psi\in C_{0}^\infty(\Rea^d)$ be defined by 
\beqs
\psi(\bx) := 
\left\{
\begin{array}{ll}
C \exp\left( \frac{1}{|\bx|^2 -1}\right) & \text{ if } |\bx| <1,\\
0 & \text{ if } |\bx|\geq 1,
\end{array}
\right.
\eeqs
where $C$ is chosen so that $\int_{\Rea^d} \psi(\bx) \rd \bx =1$.
Define $\psi_\delta(\bx) := \psi(\bx/\delta)/\delta^d$, so that $\psi_\delta(\bx)=0$ if $|\bx|>\delta$ and $\int_{\Rea^d} \psi_\delta(\bx) \rd \bx =1$.
Define $n_\delta$ by 
\beqs
n_\delta (\bx) := (n * \psi_\delta)(\bx)=\varmin+ (\pi * \psi_\delta)(\bx)= \varmin+ \int_{|\by|<\delta}
\pi(\bx-\by)\psi_\delta(\by)\, \rd \by.
\eeqs
Standard properties of mollifiers (see, e.g., \cite[\S C.4 Theorem 6]{Ev:98}) imply that 
$\|n_\delta- n\|_{L^2(\domain_R)}\tendo$ as $\delta\tendo$, 
and also that $\varmin\leq n_\delta(\bx)\leq \varmax$ for all $\bx \in \Rea^d$.
To show that $n_\delta$ satisfies \eqref{eq:n1}, we first observe that if $n$ is $C^1$ then \eqref{eq:n1} is equivalent to
\beqs
n + r\pdiff{n}{r} \geq \mu_2 
\eeqs
(where the derivatives are standard derivatives, as opposed to weak derivatives).
Thus it is sufficient to show that 
\beq\label{eq:Tot1}
n_\delta + r\pdiff{}{r}\left( \pi *\psi_\delta\right)\geq \varmin.
\eeq
Since $n_\delta \geq \varmin$,  \eqref{eq:Tot1} holds if $\partial\left( \pi *\psi_\delta\right)/\partial r\geq 0$.
However,
\beq\label{eq:Tot2}
(\pi * \psi_\delta)(\bx + h\vecer ) - (\pi * \psi_\delta)(\bx) = \int_{|\by|<\delta} \Big[\pi(\bx -\by+ h\vecer ) - \pi(\bx-\by)\Big]\psi(\by)\, \rd \by
\eeq
which is $\geq 0$ for all $h\geq 0$, since $\pi$ is monotonically nondecreasing in the radial direction; hence $\partial (\pi * \psi_\delta)/\partial r\geq 0$.

We now define $A_\delta$ in an analogous way: $\Pi$ is extended so that $\Pi(\bx):=0$ inside the scatterer, and 
\beqs
A_\delta (\bx) :=(A * \psi_\delta)(\bx)=A_{\max} I -(\Pi * \psi_\delta)(\bx)
\eeqs
(where the convolution is understood element-wise). Again, standard properties of mollifiers imply  that $\|A_\delta-A\|_{L^2(\domain_R)} \tendo$ as $\delta\tendo$ 
and $A_\delta(\bx)\geq A_{\min}$ for all $\bx \in \Rea^d$ (in the sense of quadratic forms).
Finally, \eqref{eq:Tot2} with $\pi$ replaced by $\Pi$ shows that 
$\partial (\Pi * \psi_\delta)/\partial r\geq 0$, and so $A_\delta$ satisfies \eqref{eq:A1} with $\mu_1= A_{\min}$.

Having achieved the desired approximations of $A$ and $n$, we now use them to prove that the bound \eqref{eq:bound1} holds under the assumption of existence. This proof follows a similar format to the proof of Theorem \ref{thm:EDP1}: we approximate $u$ by a smooth function $u_\eps$, and then write $u_\eps$ as $u' + u''$, for suitably chosen $u'$ and $u''$.
As in the proof Theorem \ref{thm:EDP1}, let  $W:= \{ \phi|_{\Omega_R} : \phi \in C_{0}^\infty(\Oe)\}$. In this proof only, we use the notation that $\vertiii{u}^2 := A_{\min} \|\gu\|_{L^2(\Omega_R)}^2+ n_{\min} k^2 \|u\|^2_{L^2(\Omega_R)}$. Since $\vertiii{\,\cdot\,}$ is equivalent to $\|\cdot\|_{H^1_{k}(\OR)}$, and $W$ is dense in $H^1_{0,D}(\OR)$ \cite[Page 77]{Mc:00}, given $\eps>0$, there exists a $u_\eps \in W$ such $\vertiii{u-u_\eps}<\eps$, and then 
\beqs
a(u_\eps, v)= F(v) - a(u-u_\eps,v) \quad\tfa v\in \HoDk.
\eeqs
Define $a_\delta(u,v)$ by \eqref{eq:EDPa} with $A$ and $n$ replaced by $A_\delta$ and $n_\delta$ respectively; then
\begin{align*}
a_\delta(w,v)
=a(w,v) + k^2 \int_\OR (n-n_\delta) w\, \vb - \int_\OR 
\big((A-A_\delta)\gw \big)\cdot\gvb
\end{align*}
for all $w,v\in H^1_{0,D}(\OR)$. We then have that 
\beqs
a_\delta(u_\eps,v)= F(v)  - a(u-u_\eps,v)
+ k^2 \int_\OR (n-n_\delta) u_\eps\, \vb
- \int_\OR 
\big((A-A_\delta) \gu_\eps\big)\cdot \gvb
\eeqs
for all $v\in H^1_{0,D}(\OR)$. Given $u$ and $u_\eps$, let $u', u'' \in H^1_{0,D}(\OR)$ be the solutions of the variational problems
\beqs
a_\delta(u', v) = F(v) 
\quad\tfa v\in \HoDk
\eeqs
and 
\beqs
a_\delta(u'', v) = -a(u-u_\eps,v)
+ k^2 \int_\OR (n-n_\delta) u_\eps\, \vb
- \int_\OR 
\big((A-A_\delta) \gu_\eps\big)\cdot \gvb
 \quad\tfa v\in \HoDk
\eeqs
(it is straightforward to check that the right-hand side of this last equation is a well-defined functional on $\HoDk$).
Since $A_\delta$ and $n_\delta$ satisfy Condition \ref{cond:1}, Theorem \ref{thm:EDP1} implies that $u'$ exists and is unique and Corollary \ref{cor:H1}
implies that $u''$ exists and is unique.

The key features of these definitions are that (i) $u_\eps= u' + u''$, (ii) both $u'$ and $u''$ satisfy EDPs with coefficients $A_\delta$ and $n_\delta$, (iii) $u'$ satisfies the EDP with the same right-hand side as $u$, and (iv)
$u''$ satisfies the EDP with each term on the right-hand side depending on one of $u-u_\eps$, $n-n_\delta$, and $A-A_\delta$, each of which can be made arbitrarily small.

Indeed, since $A_\delta$ and $n_\delta$ satisfy Condition \ref{cond:1},
 by Theorem \ref{thm:EDP1},
the bound \eqref{eq:bound1} holds 
with $u$ replaced by $u'$, $\mu_1=A_{\min}$, and $\mu_2=\varmin$; i.e. 
\beq\label{eq:star1}
\vertiii{u'}\leq \sqrt{C_1}\N{f}_{L^2(\Omega_+)}.
\eeq
Furthermore, by Corollary \ref{cor:H1} and the definition of the norm on $(\HoDk)'$, 
\beq\label{eq:102}
\N{u''}_{\HoDkk} \leq \widetilde{C} \left(C_c \N{u-u_\eps}_{\HoDkk} +k \N{(n-n_\delta)u_\eps}_{L^2(\domain)} + \N{(A-A_\delta)\gu_\eps}_{L^2(\domain)}\right),
\eeq
where $\widetilde{C}$ is the constant on the right-hand side of \eqref{eq:H1bound}
and $C_c$ is the continuity constant of the sesquilinear form $a(\cdot,\cdot)$ in the $\|\cdot\|_{H^1_k(\Omega_R)}$ norm;
the exact forms of $\widetilde{C}$ and $C_c$ are not important in what follows, the key point is that they are both independent of $\eps$ and $\delta$.
Now, since $u_\eps\in C^\infty(\overline{\OR})$,
\beqs
\N{(\var-\var_\delta)u_\eps}^2_{L^2(\OR)} = \int_\OR |\var-\var_\delta|^2 |u_\eps|^2 \leq \N{u_\eps}^2_{L^\infty(\OR)} \N{\var-\var_\delta}_{L^2(\OR)}.
\eeqs
An analogous inequality holds for $A-A_\delta$, and thus, 
since both $\|A-A_\delta\|_{L^2(\OR)}$ and $\|\var-\var_\delta\|_{L^2(\OR)} \tendo$ as $\delta\tendo$, given $\eps>0$ there exists a $\delta>0$ such that $k\N{(\var-\var_\delta)u_\eps}_{L^2(\OR)} <\eps$
and $\|(A-A_\delta)\gu_\eps\|_{L^2(\OR)} <\eps$. Therefore, using these inequalities in \eqref{eq:102} we have 
\beq\label{eq:starstar1}
\N{u''}_{\HoDkk} \leq \widetilde{C} (C_c+2)\eps.
\eeq
The result $\vertiii{u}\leq \sqrt{C_1}\|f\|_{L^2(\Omega_+)}$ then follows from combining the inequalities
\beqs
\vertiii{u}\leq \vertiii{u_\eps}+\eps \leq \vertiii{u'} + \vertiii{u''}+ \eps,
\eeqs
with the bounds \eqref{eq:star1} and \eqref{eq:starstar1} (using the fact that $\vertiii{\,\cdot\,}$ is equivalent to $\|\cdot\|_{H^1_{k}(\OR)}$),  and recalling that $\eps$ was arbitrary.
\epf

\

\bpf[Proof of Corollary \ref{cor:res1}]
The result \cite[Lemma 2.3]{Vo:99} implies that the assertion will hold if 
(i)  $R_\chi(k)$ is holomorphic for  $\Im k> 0$, (ii)  $R_\chi(k)$ is well-defined for $k\in \Rea\setminus\{0\}$ and 
there exist $C_4>0$ and $k_0>0$ such that
\beq\label{eq:pert1}
\N{R_\chi(k)}_{L^2(\Oe)\rightarrow L^2(\Oe)}\leq \frac{C_4}{k} \quad \tfa k\geq k_0.
\eeq
We note that the set-up in \cite{Vo:99} concerns scattering by a bounded, $C^\infty$ obstacle, where $A$ and $n$ are piecewise $C^\infty$; nevertheless, the particular result \cite[Lemma 2.3]{Vo:99} assumes only that the differential operator equals $\cL:= \Delta +k^2$ outside a large ball, and nothing about the scatterer or smoothness of the coefficients, and thus is applicable here.

By Theorem \ref{thm:EDP2}, $R_\chi(k)$ is well-defined for $k\in \Rea\setminus\{0\}$ and the bound on the real axis \eqref{eq:pert1} holds, and thus we only need to show that $R_\chi(k)$ is well-defined and holomorphic for $\Im k>0$.

When $\Im k>0$, an a priori bound on the solution of the EDP (with $g_D\equiv 0$) can be found using Green's identity (see, e.g., \cite[Lemma 3.3]{BaSpWu:16}, \cite[Theorem 2.7]{GaGrSp:15}). By Lemma \ref{lem:Fred}, the operator family $R_\chi(k)$ is therefore well-defined for $\Im k>0$.
Analyticity follows by applying the Cauchy--Riemann
operator $\partial/\partial \overline{k}$ to the variational problem \eqref{eq:EDPvar}; $\partial u/\partial \overline{k}$ then satisfies the EDP with $f=0$ and $g_D\equiv 0$, and is therefore zero by the uniqueness results.
\epf

\

\bpf[Proof of Corollary \ref{cor:H1}]
Existence, uniqueness, and the bound \eqref{eq:H1bound} follow from Lemma \ref{lem:H1}. The result about the inf-sup constant then follows from, e.g., \cite[Theorem 2.1.44]{SaSc:11}.
\epf

\section{Proofs of the results in \S\ref{sec:EDPresults2}}\label{sec:5a}

The proof of Theorem \ref{thm:EDP3} follows the proof of Theorem \ref{thm:EDP1} closely; the main difference is that, in the bound from the Morawetz identity, we only obtain \emph{either} $\|\gu\|^2_{L^2}$ (Part (i)) \emph{or} $\|u\|_{L^2}^2$ (Part (ii)) on the left-hand side (as opposed to the full weighted-$H^1$ norm), and we use the inequalities \eqref{eq:GreenEDP2} and \eqref{eq:GreenEDP1}  to put the missing part of the weighted-$H^1$ norm back in.

\

\bpf[Proof of Theorem \ref{thm:EDP3}]
(i) As in the proof of Theorem \ref{thm:EDP1}, we first assume that $\domain$ is $C^{1,1}$. The steps to then use the bound in this case to prove the bound in the case when $\domain$ is Lipschitz are identical to those in Theorem \ref{thm:EDP1}, and so we omit them.

We follow the proof of Theorem \ref{thm:EDP1}, 
but this time we replace $R$ by $R+1$ and also change $\alpha$; i.e. we 
 apply the Morawetz identity in $\domain_{R+1}$ with $v=u$ (justified as in the proof of Theorem \ref{thm:EDP1} since we're assuming $\Omega_-$ is $C^{1,1}$), $\beta=R+1$, and $2\alpha=d$. Recalling that we're assuming that $n\equiv 1$, and that $\supp\, f\subset \domain_R$, we find that the analogue of \eqref{eq:m1} is 
\begin{align}\nonumber 
&\int_{\domain_{R+1}} \big\langle \big(2A-(\bx\cdot\nabla)A\big)\gu, \gu\big\rangle
+ \int_{\Gamma_D} \frac{\bx\cdot\bnu}{\nu}\left|\dudnuA\right|^2  = 
2 \Re \int_{\domain_R}\overline{\cM u } \,f\\
&\qquad+\int_{\Gamma_{R+1}} (R+1)\left( \left|\pdiff{u}{r}\right|^2 + k^2 |u|^2 - |\nabla_S u|^2\right)   - 2 k (R+1)\, \Im \int_{\Gamma_{R+1}} \bar{u} \pdiff{u}{r} + d\,\Re \int_{\Gamma_{R+1}}\bar{u}\pdiff{u}{r}.  \label{eq:m1a}
\end{align}
The combination of the integrals over $\Gamma_{R+1}$ is $\leq 0$ by Lemma \ref{lem:2.1}, and the fact that $\domain_-$ is star-shaped implies that the integral over $\Gamma_D$ is $\geq 0$. Using the condition on $A$ \eqref{eq:A2} and the Cauchy-Schwarz inequality we have 
\begin{align}
\mu_3 \N{\gu}^2_{L^2(\domain_{R+1})} 
\leq 2 \left( R \N{\gu}_{L^2(\domain_R)} + \left( R +1+ \frac{d}{2k}\right)k \N{u}_{L^2(\domain_R)}\right) \N{f}_{L^2(\domain_+)}.\label{eq:m2}
\end{align}
Our choice of $\alpha$ has given us no $k^2\|u\|^2_{L^2(\domain_R)}$ on the left-hand side, but we reintroduce this term using the inequality \eqref{eq:GreenEDP2}. Indeed, recalling that $n\equiv 1$, and combining \eqref{eq:m2} with \eqref{eq:GreenEDP2} we have
\beq\label{eq:m3}
k^2  \N{u}^2_{L^2(\domain_R)} \leq \frac{2\cC}{\mu_3} \left( R \N{\gu}_{L^2(\domain_R)} + \left( R +1+ \frac{d}{2k}\right)k \N{u}_{L^2(\domain_R)}\right) \N{f}_{L^2(\domain_+)} + \frac{2}{k^2}\N{f}^2_{L^2(\domain_+)},
\eeq
where 
\beqs
\cC:= 4\left(A_{\max}+ \frac{6(A_{\max})^2}{k^2}\right).
\eeqs
Therefore, multiplying \eqref{eq:m3} by $\mu_3$ and combining with \eqref{eq:m2}, we obtain
\begin{align*}
&\mu_3 \left(\N{\gu}^2_{L^2(\domain_R)}+ k^2 \N{u}^2_{L^2(\domain_R)}\right)\\
&\leq 2(1+ \cC) \left( R \N{\gu}_{L^2(\domain_R)} + \left( R+1 + \frac{d}{2k}\right)k \N{u}_{L^2(\domain_R)}\right) \N{f}_{L^2(\domain_+)}
+ \frac{2\mu_3}{k^2}\N{f}^2_{L^2(\domain_+)},\\
&\leq \eps_1 \N{\gu}^2_{L^2(\domain_R)}+ \eps_2 k^2 \N{u}^2_{L^2(\domain_R)}
+ \left(\frac{(1+\cC)^2}{\eps_1} R^2 + \frac{(1+\cC)^2 }{\eps_2} \left(R+ 1+\frac{d}{2k}\right)^2 + \frac{2\mu_3}{k^2 }\right)\N{f}^2_{L^2(\domain_+)}
\end{align*}
for any $\eps_1, \eps_2>0$.
Choosing $\eps_1=\eps_2=\mu_3/2$, we obtain the result \eqref{eq:bound2}.

(ii) 
We apply the Morawetz identity in $\domain_{R}$ with $v=u$ but this time with $\beta=2R$ and $2\alpha=d-2$. 
In this case, however, since $A\equiv I$, our use of the identity is allowed for Lipschitz $\Omega_-$ by the density result of Lemma \ref{lem:density}. Indeed, Corollary \ref{cor:Necas} implies that the solution of the EDP is in $V(\domain_R)$ for every $R>0$. The integrated identity \eqref{eq:morid1int} holds with $D=\domain_{R}$ and $v=u$ by (i) the density of $\cD(\overline{\domain_{R}})$ in $V(\domain_{R})$ (Lemma \ref{lem:density}), and (ii) the fact that \eqref{eq:morid1int} is continuous in $v$ with respect to the topology of $V(\domain)$ (using, as in Part (i), the fact that  $\|\bT(\gamma v)\|_{L^2(\pD)}\leq C \|\gamma v\|_{H^1(\pD)}$  \cite[Lemma 4.23(i)]{Mc:00})).
Recalling that we're assuming that $A\equiv  I$ and $g_D\not\equiv 0$, the analogue of \eqref{eq:m1a} is now 
\begin{align}\nonumber 
&\int_{\domain_{R}}  \big(2n + \bx\cdot \nabla n\big)k^2 \mus\, \\ \nonumber
&+\int_{\Gamma_D} (\bx\cdot\bnu)\left(\left|\dudnu\right|^2 - |\nabla_{\Gamma_D} g_D|^2 + k^2 n |g_D|^2 \right) + 2 \Re \left( \left(\bx\cdot \overline{\nabla_{\Gamma_D}g_D} + \ri k (2R) \overline{g_D} + \frac{(d-2)}{2} \overline{g_D}\right)\dudnu\right)\\
&=2 \Re \int_{\domain_R}\overline{\cM u } \,f +\int_{\Gamma_{R}} R\left( \left|\pdiff{u}{r}\right|^2  - |\nabla_S u|^2+ k^2 |u|^2\right)   
- 2 k (2R)\, \Im \int_{\Gamma_{R}} \bar{u} \pdiff{u}{r} 
+ (d-2)\Re \int_{\Gamma_{R}}\bar{u}\pdiff{u}{r}.  \label{eq:m6}
\end{align}
If $kR\geq \sqrt{3/8}$ then the combination of integrals over $\Gamma_{R}$ is $\leq 0$ by Lemma \ref{lem:2.1mod} (note that Lemma \ref{lem:2.1} is not applicable since $2\alpha=d-2$).
Using the condition on $n$ \eqref{eq:n3} and the Cauchy-Schwarz inequality we have 
\begin{align}\nonumber
\mu_4 k^2 &\N{u}^2_{L^2(\domain_{R+1})} +\int_{\Gamma_D}(\bx\cdot\bnu)\left|\dudnu\right|^2 \\ \nonumber
&\leq
2 \left( R \N{\gu}_{L^2(\domain_R)} + \left( 2R + \frac{d-2}{2k}\right)k \N{u}_{L^2(\domain_R)}\right) \N{f}_{L^2(\domain_+)}
+\int_{\Gamma_D}(\bx\cdot\bnu)\left|\nabla_{\Gamma_D}g_D\right|^2 \\
&\qquad- 2 \Re \int_{\Gamma_D} \left(\bx\cdot \overline{\nabla_{\Gamma_D}g_D} + \ri k \big(2R\big) \overline{g_D} + \frac{(d-2)}{2} \overline{g_D}\right)\dudnu
\label{eq:m4alt}
\end{align}
(where we have used the fact that $\varmin>0$ in neglecting the $|g_D|^2$ term on $\Gamma_D$).
Using the inequality 
$|\bx|\leq L_D$ on $\Gamma_D$, along with the Cauchy-Schwarz and Cauchy inequalities, we have
\begin{align}\nonumber
&\hspace{-0.25cm}- 2 \Re \int_{\Gamma_D} \left(\bx\cdot \overline{\nabla_{\Gamma_D}g_D} + \ri k \big(2R\big) \overline{g_D} + \frac{(d-2)}{2} \overline{g_D}\right)\dudnu\\ \nonumber
&\hspace{0.5cm}\leq 2 \left( L_D \N{\nabla_{\Gamma_D} g_D}_{L^2(\Gamma_D)} + \left(2R + \frac{d-2}{2k}\right) k \N{g_D}_{L^2(\Gamma_D)}\right) \N{\dudnu}_{L^2(\Gamma_D)},\\
&\hspace{0.5cm}\leq  (\eps_1 + \eps_2) \N{\dudnu}^2_{L^2(\Gamma_D)} + \frac{1}{\eps_1}L_D^2 \N{\nabla_{\Gamma_D} g_D}^2_{\LtGD} + \frac{1}{\eps_2} \left(2R + \frac{d-2}{2k}\right)^2 k^2 \N{g_D}^2_\LtGD,\label{eq:m5alt}
\end{align}
for any $\eps_1, \eps_2>0$. 
Combining \eqref{eq:m5alt} with \eqref{eq:m4alt}, using the inequality $\bx\cdot\bnu\geq aL_D$ on $\Gamma_D$ in the left-hand side, and choosing $\eps_1=\eps_2=aL_D/4$, we obtain
\begin{align}\nonumber
\mu_4 k^2 \N{u}^2_{L^2(\domain_{R+1})} &+\frac{aL_D}{2}\N{\dudnu}^2_\LtGD  \\ \nonumber
\leq
&\,2 \left( R \N{\gu}_{L^2(\domain_R)} + \left( 2R + \frac{d-2}{2k}\right)k \N{u}_{L^2(\domain_R)}\right) \N{f}_{L^2(\domain_+)}\\
&\,+ L_D \left(1 + \frac{4}{a}\right) \N{\nabla_{\Gamma_D}g_D}^2_\LtGD + \frac{4}{aL_D} \left(2R + \frac{d-2}{2k}\right)^2 k^2 \N{g_D}^2_\LtGD.\label{eq:m5alt2}
\end{align}
Our choice of $\alpha$ has given us no $\|\gu\|^2_{L^2(\domain_R)}$ on the left-hand side, but we reintroduce this term using the inequality \eqref{eq:GreenEDP1}. Indeed, recalling that $A\equiv I$, and combining \eqref{eq:m5alt2} with \eqref{eq:GreenEDP1} we have
\begin{align}\nonumber
\N{\gu}^2_{L^2(\domain_R)} \leq &\,\, 
\frac{3 n_{\max}}{2\mu_4} \left[2 \bigg( R \N{\gu}_{L^2(\domain_R)} + \left( 2(R+1) + \frac{d-2}{2k}\right)k \N{u}_{L^2(\domain_R)}\right) \N{f}_{L^2(\domain_+)}\\ \nonumber
&\hspace{1.5cm}+ L_D \left(1 + \frac{4}{a}\right) \N{\nabla_{\Gamma_D}g_D}^2_\LtGD + \frac{4}{aL_D} \left(2R + \frac{d-2}{2k}\right)^2 k^2 \N{g_D}^2_\LtGD\bigg]\\
&\hspace{4cm}+ \frac{1}{2k^2 \varmax}\N{f}^2_{L^2(\domain_+)} + \N{g_D}_\LtGD\N{\dudnu}_{\LtGD}.\label{eq:m5alt3}
\end{align}
Combining \eqref{eq:m5alt2} and \eqref{eq:m5alt3}, we have
\begin{align*}
&\mu_4 \left(\N{\gu}^2_{L^2(\domain_R)}+k^2 \N{u}^2_{L^2(\domain_R)} \right)+\frac{aL_D}{2}\N{\dudnu}^2_\LtGD  \\ \nonumber
&\leq \left(1+\frac{3}{2}\varmax\right)\left[2 \bigg( R \N{\gu}_{L^2(\domain_R)} + \left( 2R + \frac{d-2}{2k}\right)k \N{u}_{L^2(\domain_R)}\right) \N{f}_{L^2(\domain_+)}\\
&\hspace{1.75cm}+ L_D \left(1 + \frac{4}{a}\right) \N{\nabla_{\Gamma_D}g_D}^2_\LtGD + \frac{4}{aL_D} \left(2R + \frac{d-2}{2k}\right)^2 k^2 \N{g_D}^2_\LtGD\bigg]\\&
\hspace{6cm}+ \frac{\mu_4}{2k^2 \varmax}\N{f}^2_{L^2(\domain_+)}+ \mu_4 \N{g_D}_\LtGD\N{\dudnu}_{\LtGD}.
\end{align*}
The result \eqref{eq:bound3b} then follows from using the Cauchy inequality on the right-hand side.

Although the bound \eqref{eq:bound3b} is only valid for $k \geq \sqrt{3/8}R^{-1}$, to obtain existence and uniqueness for all $k>0$ in these cases, we can use the unique continuation principle; see, e.g., \cite[Theorem 2.5]{GrSa:18} and the references in \S\ref{sec:intro}.
\epf

\section{Geometric interpretation of the condition $2n + \bx\cdot\nabla n>0$}\label{sec:6}

Before stating the main results of this section (Lemma \ref{lem:R1} and Theorem \ref{thm:Ralston}) we need to recall the definitions of the \emph{bicharacteristics} and \emph{rays} of the wave equation (see, e.g., \cite{Ho:83}, \cite{Ho:85}, \cite[Chapter VI]{CoHi:62}, and the expository accounts in \cite[\S8.7]{Str:93}, \cite{Ta:78}).
In this section, we only consider the full-space problem for the wave equation and the Helmholtz equation, i.e.~we have $\bx\in \Rea^d$, and we only consider the case when $A\equiv I$ (but see Remark \ref{rem:AnotI} for some comments on the case $A\not\equiv I$).

\begin{definition}[Bicharacteristics of the wave equation]
The \emph{bicharacteristics} of the wave equation
\beq\label{eq:R1}
\Delta U - n \pdiff{^2 U}{t^2 } =0, \quad \bx\in \Rea^d,\,\, t\in \Rea, \,\,n\in C^1(\Rea^d),
\eeq
are the solutions $(\bx(s), t(s), \bxi(s), \tau(s))$, 
with $\bx, \bxi\in \Rea^d$, $t, \tau\in \Rea$, 
of the Hamiltonian system
\begin{subequations}
\begin{align*}
\dot{x}_i &= \pdiff{}{\xi_i}p\big(x, t, \xi ,\tau\big), \quad \dot{t} = \pdiff{}{\tau}p\big(x, t, \xi ,\tau\big),\\
\dot{\xi}_i &= -\pdiff{}{x_i}p\big(x, t, \xi ,\tau\big),\quad
\dot{\tau}_i = -\pdiff{}{t}p\big(x, t, \xi ,\tau\big),
\end{align*}
\end{subequations}
where dot denotes differentiation 
with respect to $s$ and 
\beq\label{eq:p}
p\big(\bx, t, \bxi ,\tau\big):= |\bxi|^2 - n(\bx)\tau^2 
\eeq
(i.e., $p$ is the principal symbol of the wave equation \eqref{eq:R1}). That is, the bicharacteristics 
are the solutions of the system of equations
\begin{subequations}\label{eq:64}
\begin{align}\label{eq:64a}
\dot{x}_i &= 2 \xi_i,  \qquad \dot{t}\, = -2 n\tau, \\ \label{eq:64b}
\dot{\xi}_i &= \left(\pdiff{n}{x_i}(\bx)\right)\tau^2, \quad \dot{\tau}\, =0.
\end{align}
\end{subequations}
The \emph{null bicharacteristics} are the bicharacteristics for which $p=0$.
\end{definition}

\bre[Using $t$ instead of $s$ as the parameter along the bicharacteristics]\label{rem:bichar}
Since $\dot{\tau}=0$ in \eqref{eq:64b}, $\tau$ is a constant function of $s$.
Without loss of generality, assume that $\tau$ is negative. Then, since $n>0$, \eqref{eq:64a} implies that $\dot{t}>0$ and thus $t$ is an increasing function of $s$.
 We can therefore use $t$ instead of $s$ as the variable along the bicharacteristics.
\ere

\bre[Bicharacteristics and singularities]
The significance of bicharacteristics in the study of the wave equation stems from the result of 
\cite[\S VI]{DuHo:72} that (in the absence of boundaries) singularities of pseudodifferential operators (understood in terms of the \emph{wavefront set}) travel along null bicharacteristics; see, e.g., \cite[Chapter 24]{Ho:85}, \cite[\S12.3]{Zw:12}.
\ere

\begin{definition}[Rays]\label{def:ray}
The \emph{rays} of the wave equation are given by $\bx(s)$, where $(\bx(s), t(s), \bxi(s), \tau(s))$ is a null bicharacteristic of the wave equation; i.e.,
the rays are the projection of null bicharacteristics in the $\bx$ variables.
\end{definition}

\ble[Geometric interpretation of the condition $2n + \bx\cdot\nabla n>0$]
\label{lem:R1}
If $\bx(s)$ is a ray of the wave equation (defined by Definition \ref{def:ray}), then
\beq\label{eq:R9}
\half\diff{^2}{s^2} \left(\half \big|\bx(s)\big|^2\right) = 2 n\big(\bx(s)\big) + \bx(s)\cdot \nabla n\big(\bx(s)\big).
\eeq
\ele

\noi We postpone the proof of Lemma \ref{lem:R1} to \S\ref{sec:Rproof} (after we have recapped the necessary results from semiclassical analysis in \S\ref{sec:RSCA}), and  we now use this result to relate the condition $2n + \bx\cdot\nabla n>0$ to $n$ being trapping/nontrapping.

\begin{definition}[Trapping/nontrapping $n$]
Given $n\in C^1(\Rea^d)$ with $\supp(1-n)$ compact, consider the bicharacteristics defined by \eqref{eq:64} with $\tau=-1$ and thought of as a function of $t$ by Remark \ref{rem:bichar}. We say that $n$ is \emph{nontrapping} if, given $B_R(\bze)\supset \supp\, n$, there exists a $T(R)>0$ such that all bicharacteristics with $|\bx(t_0)|<R$ satisfy $|\bx(t)|>R$ for all $t\geq T(R)$; i.e.~all rays starting inside $B_R(\bze)$ at time $t_0$ have left $B_R(\bze)$ by time $T$. 
We say that $n$ is \emph{trapping} if $n$ is not nontrapping.
\end{definition}

\begin{theorem}[Relationship of the condition $2n + \bx\cdot \nabla n>0$ to trapping/nontrapping]\label{thm:Ralston}
Let $n\in C^1(\Rea^d)$ with $\supp(1-n)$ compact.

(i) If there exists a $\mu>0$ such that $2n + \bx\cdot \nabla n \geq \mu$ for all $\bx\in \Rea^d$, then $n$ is nontrapping.

(ii) (Due to \cite{Ra:71}.) If $d=3$, $n(\bx)= n(r)$ (with $r:=|\bx|$) and 
there exists an $r_1$ such that
\beq\label{eq:Ralston}
2 n(r_1) + r_1 \pdiff{n}{r}(r_1)<0,
\eeq
then $n$ is trapping and there exists a sequence of resonances exponentially close to the real-$k$ axis.
\end{theorem}

\bpf
(i) Since $\nmin>0$ and $\rd t/\rd s= 2n \geq 2\nmin$ it is sufficient to show that, with the bicharacteristics thought of as a function of $s$, there exists an 
$S(R)>0$ such that all bicharacteristics with $|\bx(s_0)|<R$ satisfy $|\bx(s)|>R$ for all $s\geq S(R)$. Without loss of generality, we take $s_0=0$.

If $f\in C^2(\Rea)$ with $(\rd^2 f/\rd s^2)(s)\geq \mu>0$ for all $s$, then
\beqs
\diff{f}{s}(s)-\diff{f}{s}(0) = \int_{0}^s \diff{^2 f}{s^2}(s) \, \rd s \geq \mu s, \quad\text{and similarly}\quad
f(s) \geq \frac{\mu s^2}{2} + \left(\diff{f}{s}(0)\right)s + f(0).
\eeqs
Thus, given $f(0)$ and $(\rd f/\rd s)(0)$ with $f(0)<R^2/2$ and $|(\rd f/\rd s)(0)|<\infty$, there exists an $S(R)$ such that $f(s) > R^2/2$ for all $s\geq S(R)$. Applying this result to $f(s) = |\bx(s)|^2/2$, and using Lemma \ref{lem:R1}, we obtain the result in (i).

(ii) With $n=c^{-2}$, the condition $2n + \bx\cdot\nabla n>0$ becomes $\bx\cdot \nabla c  <c$, which is $r c'(r)<c$ if $c$ is radial. This condition can then be rewritten as 
\beq\label{eq:cond_c}
\diff{}{r} \left(\frac{r}{c(r)}\right)>0 \quad\text{ or } \quad \diff{}{r} \left(\frac{c(r)}{r}\right)<0. 
\eeq
Since there exists an $r_1$ such that \eqref{eq:Ralston} holds, $(\rd/\rd r)( c(r)/r) |_{r=r_1}>0$ and then, since $c\in C^1$, there 
exist $r_0$ and $r_2$ with $r_0<r_1<r_2$ and
\beq\label{eq:Ralston1}
\frac{c(r_0)}{r_0}<\frac{c(r_2)}{r_2};
\eeq
without loss of generality, we can take $r_2$ such that $c(r)/r$ is monotonically decreasing for $r\geq r_2$ (observe that the assumption $\supp(1-n)$ is compact implies that $c\equiv 1$ for sufficiently large $r$).
Ralston showed in \cite{Ra:71} that if (a) \eqref{eq:Ralston1} holds, and (b) one has energy decay for solutions of the wave equation in $r>r_2$ with zero Dirichlet boundary conditions imposed at $r=r_2$, then there is a sequence of resonances exponentially close to the real-$k$ axis. 
Since $c(r)/r$ is monotonically decreasing for $r\geq r_2$, $2n + \bx\cdot\nabla n>0$ in this region, and then \emph{either} the results of \cite{Li:64}, \emph{or} \cite[Theorem 1.1]{Vo:99} combined with the bound \eqref{eq:bound3b} show that (b) holds. 
\epf

\bre[The condition $2n + \bx \cdot \nabla n>0$ and the boundary rigidity problem]
We highlight that the condition \eqref{eq:cond_c} also appears in \cite{He:05, WeZo:07} 
as a sufficient condition for the solution of the \emph{boundary rigidity problem} in the case of a sphere when $c=c(r)$. Recall that in this problem,
one seeks to find $c$ from knowing the travel times (i.e.~knowing the geodesics between any two points of the boundary); see, e.g., \cite{StUhVa:17} and the references therein.
\ere

\subsection{Semiclassical principal symbol and the associated bicharacteristics}\label{sec:RSCA}

The moral of this section is that the rays of the wave equation \eqref{eq:R1} can be understood by studying the bicharacteristics of \emph{semiclassical principal symbol} of the Helmholtz equation.

Given the differential operator (in multi-index notation)
\beq\label{eq:scpde}
P(\bx, h) = \sum_{|\balpha|\leq N} a_\balpha(\bx) 
\big(-\ri h\partial\big)^\balpha, 
\eeq
its \emph{semiclassical principal symbol} $\princP(\bx,\bxi)$ is given by 
\beq\label{eq:scps}
\princP(\bx, \bxi) =\sum_{|\balpha|\leq N} a_\balpha(\bx)\bxi^\balpha;
\eeq
see, e.g., \cite[\S12.3]{Zw:12}. 
Writing the Helmholtz equation $\Delta u + k^2 n u=0$ as 
\beq\label{eq:R3}
L u:=-h^2 \Delta u- nu=0
\eeq
with $h=k^{-1}$, we then have that the semiclassical principal symbol is given by 
\beq\label{eq:R4}
\sigma(L)=|\bxi|^2 - n(\bx).
\eeq

\begin{definition}[Semiclassical bicharacteristics of the Helmholtz equation]
The \emph{semiclassical bicharacteristics} of the Helmholtz equation \eqref{eq:R3} are the solutions $(\bx(s),\bxi(s))$, with $\bx, \bxi \in \Rea^d$, of the Hamiltonian system,
\begin{align*}
\dot{x}_i = \pdiff{}{\xi_i}\sigma(L)\big(\bx, \bxi \big), \qquad
\dot{\xi}_i = -\pdiff{}{x_i}\sigma(L)\big(\bx, \bxi \big),
\end{align*}
with $\sigma(L)(\bx,\bxi)$ given by \eqref{eq:R4}. That is, the bicharacteristics are the solution of the system
\begin{align}\label{eq:R6}
\dot{x}_i = 2\xi_i, \qquad
\dot{\xi}_i = \pdiff{n}{x_i}(\bx).
\end{align}
The \emph{null semiclassical bicharacteristics} are the semiclassical bicharacteristics for which $\sigma(L)=0$.
\end{definition}

\begin{corollary}\label{cor:ray}
The rays of the wave equation are the projection of the null semiclassical bicharacteristics of the Helmholtz equation in the $\bx$ variables.
\end{corollary}

\bpf
From the second equation in \eqref{eq:64b}, a bicharacteristic of the wave equation has $\tau$ as a constant function of $s$. 
By rescaling, we can take $\tau= -1$ without loss of generality, and then $p$ (defined by \eqref{eq:p}) equals $\sigma(L)$ (defined by \eqref{eq:R4}).
The rays are the projections of the solutions of \eqref{eq:64}, when $p=0$, in the $\bx$ variables, and, since $\tau=-1$, the rays are governed by the first equations in each of \eqref{eq:64a} and \eqref{eq:64b}. 
Since $p=\sigma(L)$, these two equations are identical to the equations of the semiclassical bicharacteristics of the Helmholtz equation \eqref{eq:R6}. 
\epf

\

For $f$ and $g$ functions of $\bx, \bxi$, and $s$, 
we define the Poisson bracket $\{f, g\}$ by
\beq\label{eq:R8}
\{f,g\}:= \sum_j \left( 
\pdiff{f}{\xi_j}\pdiff{g}{x_j}- \pdiff{f}{x_j}\pdiff{g}{\xi_j}
\right),
\eeq
(see, e.g., \cite[\S2.4]{Zw:12})
and recall that for a quantity $f(\bx(s), \bxi(s); s)$ evolving via the Hamiltonian flow associated with the Hamiltonian $H$,
\beq\label{eq:R7}
\diff{f}{s} = \{H, f\} + \pdiff{f}{s}
\eeq
by the chain rule. 
Recalling the definition of the commutator $[A,B]:= AB-BA$, we have that, for $A$ and $B$  differential operators of the form \eqref{eq:scpde},
\beq\label{eq:R13}
\sigma\big( [A, B]\big) = \frac{h}{\ri} \{ \sigma(A), \sigma(B)\};
\eeq
see, e.g., \cite[Remark (ii), Page 68]{Zw:12}.

\subsection{Proof of Lemma \ref{lem:R1}}\label{sec:Rproof}

By Corollary \ref{cor:ray}, the rays are described by the Hamiltonian flow associated with $\sigma(L)$ \eqref{eq:R4}.
By \eqref{eq:R7} and the definition of the Poisson bracket \eqref{eq:R8}, for $\bx(s)$ on a ray,
\beq\label{eq:R9aa}
\diff{}{s}\left( \half |\bx|^2\right)= \left\{|\bxi|^2 - n(\bx), \half  |\bx|^2\right\} = 2\sum_j x_j \xi_j = 2\bx\cdot \bxi
\eeq
(where we have suppressed the dependence of $\bx$ and $\bxi$ on $s$ to keep the expressions compact).
Thus,
\beq\label{eq:R9a}
\diff{^2}{s^2}\left( \half |\bx|^2\right)= \left\{
 |\bxi|^2-n(\bx),\diff{}{s}\left( \half |\bx|^2\right)
\right\}
=\Big\{ |\bxi|^2 - n(\bx),  2\bx\cdot \bxi\Big\}= 2\big(2 |\bxi|^2+\bx \cdot\nabla n(\bx)\big).
\eeq
Since 
$\sigma(L)(\bx(s), \bxi(s))=0$ for all $s$ (as the rays are the \emph{null} semiclassical bicharacteristics), we have from \eqref{eq:R4} that $|\bxi|^2 = n(\bx)$ on rays, and thus \eqref{eq:R9a} becomes \eqref{eq:R9}.

\subsection{Understanding why the condition $2n(\bx) + \bx\cdot \nabla n(\bx)>0$ arises from using Morawetz identities}

The explanation of why the condition $2n(\bx) + \bx\cdot \nabla n(\bx)>0$ arises when we use Morawetz identities to obtain bounds on the solution of $(\Delta + k^2 n)u=-f$ has three points.
\ben
\item The Morawetz/Rellich identities in \S\ref{sec:4} can be understood as commutator arguments.
\item In the setting of commutator arguments (and ignoring any contributions from infinity or boundaries), the Morawetz/Rellich identities control the first derivatives of the solution if $\rd^2 (|\bx(s)|^2/2)/\rd s^2 >0$ when $\bx(s)$ is a ray.
\item As we saw in Lemma \ref{lem:R1}, when $\bx(s)$ is a ray, $\rd^2 (|\bx(s)|^2/2)/\rd s^2  >0$ is equivalent to $2n(\bx) + \bx\cdot \nabla n(\bx)>0$.
\een
We now explain Points 1 and 2 (in \S\ref{sec:comm} and \S\ref{sec:final_sec} respectively).

\subsubsection{Morawetz/Rellich identities as commutator arguments}\label{sec:comm}

With $P$ a formally self-adjoint operator, we consider the PDE $Pu =- f$ posed on the whole space $\Rea^d$ (i.e.~we ignore any complications due to boundaries) and we proceed formally without worrying about boundedness of norms over $\Rea^d$. 
We depart slightly from the notation in the rest of the paper and use $\langle\cdot,\cdot\rangle$ to denote the $L^2$ inner-product, and $\|\cdot\|$ to denote the $L^2$ norm.

The heart of a commutator argument consists of finding 
operators $M$, $B$, and $Q$, such that
\beq\label{eq:R11}
\big\langle [M,P]u, u \big\rangle = \pm \N{B u}^2 + \big\langle QPu, u\rangle.
\eeq
Without loss of generality, assume that \eqref{eq:R11} holds with the plus sign. If $\|u\|$ can be controlled by $\|B u\|$, then \eqref{eq:R11} gives 
an estimate on $\|u\|$ in terms of $\|M f \|$, $\|M^* f \|$, and $\|Qf \|$.
Indeed, \eqref{eq:R11} and the PDE $Pu=-f$ imply that
\beq\label{eq:R11a}
\N{Bu}^2 = \big\langle Qf, u\big\rangle - \big\langle Mf ,u \big\rangle +\big\langle Mu ,f \big\rangle.
\eeq
After bounding the left-hand side below by $\|u\|$ and converting the last inner-product to $\langle u, M^* f\rangle$, an estimate on $\|u\|$ in terms of $\|M f \|$, $\|M^* f \|$, and $\|Qf \|$ can be obtained by the Cauchy-Schwarz inequality (see also the explanation of this type of argument in \cite[\S6]{RoTa:15}).

We claim that, when $P=\Delta$, \eqref{eq:R11} is satisfied with the plus sign, $M= \bx\cdot \nabla$, $B = \sqrt{2}\nabla$, and $Q=0$,
and furthermore that this is equivalent to the Rellich identity \eqref{A} with $v=u$, $A\equiv I$, and $k=0$.
Indeed, one can easily check that $[\bx \cdot \nabla, \Delta] = -2 \Delta$, and then \eqref{eq:R11} holds since
\beq\label{eq:Delta}
\langle \Delta u, u\rangle = -\| \nabla u\|^2
\eeq
by integration by parts (recalling that we're ignoring any boundary terms). Then \eqref{eq:R11a} becomes 
\beq\label{eq:R12}
2\|\nabla u\|^2 = -\big\langle \bx\cdot\nabla f, u\big\rangle +\big\langle \bx\cdot \nabla u, f\big\rangle.
\eeq
The Rellich identity \eqref{A} with $v=u$, $A\equiv I$, $k=0$, and $\Delta u =-f$, integrated over $\Rea^d$, becomes 
\beq\label{eq:R12a}
(d-2)\|\nabla u\|^2 = -2 \Re\big\langle \bx \cdot \nabla u , f\big\rangle.
\eeq
To reconcile \eqref{eq:R12} and \eqref{eq:R12a}, observe that
\beqs
\big\langle  \bx\cdot\nabla f, u\big\rangle=-\big\langle  f, \nabla\cdot(\bx \,u)\big\rangle = - d\big\langle f, u\rangle- \big\langle f, \bx\cdot \gu\big\rangle 
\eeqs
and then using both this and \eqref{eq:Delta} in  \eqref{eq:R12} we obtain \eqref{eq:R12a}.

Here we have focused on the case $P=\Delta$, but similar considerations apply when $P=\Delta +k^2$ because of the commutator relation 
$[\bx \cdot \nabla, (\Delta+k^2)] = -2 \Delta$

\subsubsection{Morawetz/Rellich identities and the condition $\rd^2 (|\bx|^2/2)/\rd s^2 >0$}\label{sec:final_sec}

As in \eqref{eq:R3}, we let $L:= - h^2 \Delta -n(\bx)$. 
The arguments in \S\ref{sec:comm} (with $Q=0$) indicate that we should try to make the commutator $[-\bx \cdot\nabla, L]$
single-signed, i.e.~a positive or negative operator
Seeking to achieve single-signedness at highest-order is then equivalent to making $\sigma([-\bx \cdot\nabla, L])$
single-signed. 

From the definition \eqref{eq:scps}, we have
\beq\label{eq:FF1}
\sigma\big( -\ri h\,\bx \cdot \nabla\big)= \bx\cdot \bxi.
\eeq
Using \eqref{eq:FF1} along with \eqref{eq:R13}, \eqref{eq:R4}, and \eqref{eq:R9a}, we find that
\beqs
\sigma\Big( \big[ -\bx\cdot\nabla, L\big]\Big)=-\frac{\ri}{h}\sigma\Big( \big[ -\ri h\,\bx\cdot\nabla, L\big]\Big) =-\Big\{\bx\cdot \bxi,\, |\bxi|^2-n(\bx)\Big\}= \half \diff{^2}{s^2}\left( \half |\bx|^2\right),
\eeqs
and thus we expect to get a bound from using the multiplier $\bx \cdot \nabla$ when $\rd^2 (|\bx(s)|^2/2)/\rd s^2>0$. 
All these calculations have ignored the contribution from infinity (which is dealt with using the Morawetz-Ludwig multiplier in Lemma \ref{lem:ML}) and from boundaries (which is where 
we require Dirichlet boundary conditions on the obstacle).

\bre[The operator  $-h^2\nabla\cdot(A(\bx) \nabla ) - n(\bx)$ with $A\not\equiv I$]\label{rem:AnotI}
When $L:= -h^2\nabla\cdot(A(\bx) \nabla) - n(\bx)$ with $A\not\equiv I$, the multiplier $\bx\cdot \nabla$ no longer has such a nice connection to rays. Indeed, 
the connection of the multiplier $\bx\cdot \nabla$ to rays when $A\equiv I$ relied on the fact that, when $\bx(s)$ is a ray,
\beqs
\half\diff{}{s}\left( \half |\bx(s)|^2\right)=\sigma\left(\frac{h}{\ri}\bx\cdot\nabla\right);
\eeqs
see \eqref{eq:R9aa} and \eqref{eq:FF1}.
However, we now have that
\beqs
\sigma(L)(\bx,\bxi) = \big(A(\bx)\bxi\big)\cdot\bxi 
- n(\bx),
\eeqs
and thus
\beqs
\half\diff{}{s}\left( \half |\bx|^2\right)= 
\half\left\{ \sigma(L) , \half |\bx|^2\right\} =\bx \cdot \big(A(\bx)\bxi\big),
\eeqs
which is not $\bx\cdot\bxi= \sigma((h/\ri)\bx\cdot \nabla)$ when $A\not\equiv I$.
Furthermore, when $A\not\equiv I$, the modulus function no longer measures distance with respect to the metric induced by $A$, and so 
the quantity 
$\rd^2 (|\bx(s) |^2/2)/\rd s^2$ 
is less connected to the geometry (and is therefore less meaningful) than when $A\equiv I$.
\ere

\begin{appendix}

\section{The truncated exterior Dirichlet problem (TEDP)}\label{sec:TEDP}

As described in the introduction, there has been recent interest in proving bounds on this problem when one or both of $A$ and $n$ are variable \cite{BrGaPe:15, Ch:15, BaChGo:17, OhVe:18, MoSp:17, SaTo:17, GrSa:18}. We therefore give in this section the analogues for the TEDP of our results for the EDP, indicating how the proofs of the EDP results need to be modified to obtain these.

\subsection{Definition of the TEDP}

\begin{definition}[Truncated Exterior Dirichlet Problem (TEDP)]\label{def:TEDP}
Let $\domain_-$ be a bounded Lipschitz open set such that the open complement $\domain_+:= \Rea^d\setminus \overline{\domain_-}$ is connected. 
Let $\widetilde{\domain}$ be a bounded connected Lipschitz open set such that $\overline{\domain_-}\subset\subset \widetilde{\domain}$. 
Let $\domain:=\widetilde{\domain}\setminus\overline{\domain_-}$, $\Gamma_D:= \partial \domain_-$, and $\Gamma_I :=\partial \widetilde{\domain}$, so that $\partial \domain= \Gamma_D \cup \Gamma_I$ and $\Gamma_D\cap \Gamma_I = \emptyset$ (see Figure \ref{fig:TEDP}).

Given 
\bit
\item $f\in L^2(\domain)$ 
\item $g_D\in H^1(\Gamma_D)$,
\item $g_I\in L^2(\Gamma_I)$
\item $\vartheta\in L^\infty(\Gamma_I, \Rea)$ with 
\beq\label{eq:thetaineq}
0<\vartheta_{\min} \leq \vartheta(\bx)\leq\vartheta_{\max}<\infty\,\, \text{ for almost every } \bx \in \Gamma_I,
\eeq
\item $n\in L^\infty(\domain,\Rea)$ satisfying \eqref{eq:nlimitsEDP} with $\domain_+$ replaced by $\domain$,
\item $A \in L^\infty(\domain , \Rea^{d\times d})$ such that $\dist(\supp(I -A),\Gamma_I)>0$, $A$ is symmetric, and there exist $0<A_{\min}\leq A_{\max}<\infty$ such that
\eqref{eq:AellEDP} holds with $\domain_+$ replaced by $\domain$,
\eit
we say $u\in H^1(\domain)$ satisfies the truncated exterior Dirichlet problem if 
\beqs
\opL_{A,n} u:= \nabla\cdot(A \gu ) + k^2 n u = -f \quad \tin \domain,
\eeqs
\beqs
\gamma u =g_D \quad\ton \Gamma_D,
\eeqs
and 
\beq\label{eq:TEDP3}
\dudnu - \ri k\vartheta  \gamma u = g_I \ton \Gamma_I.
\eeq
\end{definition}

\begin{figure}
\begin{centering}
\begin{tikzpicture}[scale=0.4,even odd rule]
    \pgfdeclarepatternformonly{owennortheast}
    {\pgfpointorigin}{\pgfpoint{1cm}{1cm}}
    {\pgfpoint{1cm}{1cm}}
    {
    \pgfpathmoveto{\pgfpointorigin}
    \pgfpathlineto{\pgfpoint{1cm}{1cm}}
    \pgfsetlinewidth{0.5\pgflinewidth}
    \pgfusepath{stroke}
    }
    
     \pgfdeclarepatternformonly{owennorthwest}
    {\pgfpointorigin}{\pgfpoint{1cm}{1cm}}
    {\pgfpoint{1cm}{1cm}}
    {
    \pgfpathmoveto{\pgfpoint{1cm}{0cm}}
    \pgfpathlineto{\pgfpoint{0cm}{1cm}}
   \pgfsetlinewidth{0.5\pgflinewidth}
    \pgfusepath{stroke}
    }

    \filldraw[pattern=owennorthwest] 
    (6,5) .. controls (5.5,3) .. 
    (7,0) -- 
    (5,-5.5) .. controls (2.5,-4) .. 
    (0,-5) .. controls (-3.75,-6.5) .. 
    (-7.5,-5) -- 
    (-7,-3) -- 
    (-8,2) -- 
    (-5,6) -- 
    (0,5) -- 
    (3,7) .. controls (4.2,5.4) ..
    cycle
    
   [scale=2/3] (4.5,0) .. controls (4.5,-0.5) and (4.25,-2) .. 
    (3,-3) .. controls (2.875,-3.1) and (2.25,-3.5) ..
    (1.5,-3) .. controls (0.75,-2.5) and (0.5,-2.5) ..
    (-0.75,-3) .. controls (-2,-3.5) and (-2.8,-3.4)..
    (-3,-3) .. controls (-3.75,-1.5) and (-4.95,-1.1) ..
    (-5.25,-1) .. controls (-6,-0.75) and (-6,0.25) ..
    (-4.5,1.5) .. controls (-4.35,1.625) and (-3.5,2.25) ..
    (-3,3.5) .. controls (0,3) and (0.6,2.9) ..
    (1.5,3.5) .. controls (2.25,4) and (3.5,3.5) ..
    (3.75,3) .. controls (4,2.5) and (4.5,0.5)..
    cycle;

    \filldraw[pattern=owennortheast,scale=2/3] 
    (4.5,0) .. controls (4.5,-0.5) and (4.25,-2) .. 
    (3,-3) .. controls (2.875,-3.1) and (2.25,-3.5) ..
    (1.5,-3) .. controls (0.75,-2.5) and (0.5,-2.5) ..
    (-0.75,-3) .. controls (-2,-3.5) and (-2.8,-3.4)..
    (-3,-3) .. controls (-3.75,-1.5) and (-4.95,-1.1) ..
    (-5.25,-1) .. controls (-6,-0.75) and (-6,0.25) ..
    (-4.5,1.5) .. controls (-4.35,1.625) and (-3.5,2.25) ..
    (-3,3.5) .. controls (0,3) and (0.6,2.9) ..
    (1.5,3.5) .. controls (2.25,4) and (3.5,3.5) ..
    (3.75,3) .. controls (4,2.5) and (4.5,0.5)..
    cycle;


\draw (0,0) node[fill=white] {$\domain_{-}$};
\draw (-3.5,-3.5) node[fill=white] {$\domain:=\widetilde{\domain}\setminus \overline{\domain_-}$};

\draw (3.4,-2.2) node[fill=white] {$\Gamma_{D}$};

\draw (3.8,7) node {$\Gamma_{I}$};
    
\end{tikzpicture}

  \caption{The domains $\Omega$ and $\Oi$, and boundaries $\Gamma_I$ and $\Gamma_D$, in the definition of the TEDP (Definition \ref{def:TEDP}).}\label{fig:TEDP}
  \end{centering}
  \end{figure}
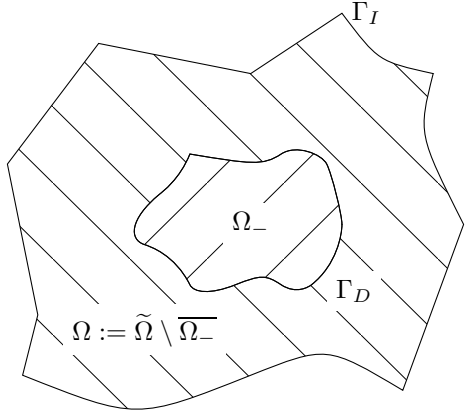

\bre[The behaviour of $A$, $n$, and $\vartheta$ on the impedance boundary $\Gamma_I$]
In the definition of the TEDP, we have made $A\equiv I$ in a neighbourhood of the impedance boundary $\Gamma_I$, but allowed $n$ to vary in this neighbourhood. With the impedance boundary condition viewed as an approximation to the Sommerfeld radiation condition \eqref{eq:src}, it would perhaps make more sense to impose the condition that $n\equiv 1$ in a neighbourhood of $\Gamma_I$ and then let the function $\vartheta$ in the impedance condition \eqref{eq:TEDP3} be equal to one. However, the interior impedance problem with $A\equiv I$, $n$ varying in the whole of the domain, and $\vartheta$ varying on the impedance boundary is considered in \cite{BaChGo:17, BrGaPe:15, Ch:15, FeLiLo:15, GrSa:18}, and so we include this situation in Definition \ref{def:TEDP} to make contact with these other works. We note that $\vartheta=n^{1/2}$ in \cite{FeLiLo:15, GrSa:18}, $\vartheta=(\varmax)^{1/2}$
in \cite{BaChGo:17, Ch:15}, and $\vartheta$ is a general function (satisfying \eqref{eq:thetaineq}) in \cite{BrGaPe:15}.
\ere

\subsection{Bounds on the TEDP for Lipschitz $A$ and $n$}

The following three sets of conditions are the TEDP-analogues of Conditions \ref{cond:1}, \ref{cond:2}, and \ref{cond:3b}.

\begin{condition}\textbf{\emph{($A$ and $n$ both Lipschitz, $g_D\equiv 0$, $\domain_-$ star-shaped, $\widetilde{\domain}$ star-shaped w.r.t.~a ball)}}\label{cond:1alt}
$\domain_-$ is star-shaped with respect to the origin, $\widetilde{\domain}$ is star-shaped with respect to a ball centred at the origin, $g_D\equiv 0$, $A\in C^{0,1}(\overline{\domain},\Rea^{d\times d})$, $n\in C^{0,1}(\overline{\domain}, \Rea)$, and there exist $\mu_1, \mu_2>0$ such that \eqref{eq:A1} and \eqref{eq:n1} hold with $\domain_+$ replaced by $\domaingen$.
\end{condition}

\begin{condition}[$A$ Lipschitz, $n\equiv 1$, $g_D\equiv 0$, $\domain_-$ star-shaped, $\widetilde{\domain}$ star-shaped w.r.t.~a ball]\label{cond:2alt}
$\domain_-$ is star-shaped with respect to the origin, $\widetilde{\domain}$ is star-shaped with respect to a ball centred at the origin,
$g_D\equiv 0$, $A\in C^{0,1}(\overline{\domain},\Rea^{d\times d})$, $n\equiv 1$, and there exists $\mu_3>0$ such that
\eqref{eq:A2}
holds with $\domain_+$ replaced by $\domaingen$.
\end{condition}

\begin{condition}[$n$ Lipschitz, $g_D\not\equiv 0$, both $\domain_-$ and $\widetilde{\domain}$ star-shaped w.r.t.~a ball]\label{cond:3balt}
$\domain_-$ is star-shaped with respect to a ball centred at the origin, $\widetilde{\domain}$ is star-shaped with respect to a ball centred at the origin,
$A\equiv I$, $n\in C^{0,1}(\overline{\domain})$, and there exists $\mu_4>0$ such that \eqref{eq:n3} holds with $\domain_+$ replaced by $\domaingen$.
\end{condition}

The following theorem is the TEDP-analogue of Theorems \ref{thm:EDP1} and \ref{thm:EDP3}.

\begin{theorem}[Bounds on the TEDP under Conditions \ref{cond:1alt}-\ref{cond:3balt}]\label{thm:TEDP}

\

\noi (i) If $\domain_-, \widetilde{\domain}, A, n, \vartheta$, $f$, $g_D$, and $g_I$ satisfy the requirements in the definition of the TEDP (Definition \ref{def:TEDP}), along with the requirements in Condition \ref{cond:1alt}, then the solution of the TEDP exists and is unique. Let $L_I:= \max_{\bx\in \Gamma_I}|\bx|$ and let $aL_I$ be the radius of the ball with respect to which $\widetilde{\domain}$ is star-shaped. 
Finally let $n_{\max,\Gamma_I}:=\esssup_{\bx\in\Gamma_I}n(\bx)$.
Then
\beqs
\mu_1 \N{\gu}^2_{L^2(\domain)} + \mu_2 k^2 \N{u}^2_{L^2(\domain)} +{a L_I}\N{\nT (\gamma u)}^2_{\LtGI} + 2L_I k^2\N{\gamma u}^2_{\LtGI}
\leq C_1 \N{f}^2_{L^2(\domain)} + \widetilde{C_1} \N{g_I}^2_{\LtGI}
\eeqs
for all $k>0$, where
\beq\label{eq:C1tilde}
C_1 := 4\left[\frac{ L_I^2}{\mu_1} + \frac{1}{\mu_2}\left(\beta+ \frac{d-1}{2k}\right)^2\right],
\qquad
\widetilde{C_1}:=
2\left[ 2 \left(1+ \frac{2}{a}\right) + \frac{\beta}{\vartheta_{\min}L_I}+ \frac{(d-1)^2}{4}\right]
L_I. 
\eeq
and
\beq\label{eq:beta}
\beta:= \frac{L_I}{\vartheta_{\min}} \left(1 + n_{\max,\Gamma_I} + \frac{1}{(kL_I)^2} + 2 (\vartheta_{\max})^2\left(1 + \frac{2}{a}\right)\right).
\eeq

\noi (ii) 
 If $\domain_-, \widetilde{\domain}, A, n,\vartheta, f$, $g_D$, and $g_I$ satisfy the requirements in the definition of the TEDP (Definition \ref{def:TEDP}), along with the requirements in Condition \ref{cond:2alt}, then the solution of the TEDP exists and is unique. Let $L_I:= \max_{\bx\in \Gamma_I}|\bx|$ and let $aL_I$ be the radius of the ball with respect to which $\widetilde{\domain}$ is star-shaped. Then
\beqs
\mu_3 \left(\N{\gu}^2_{L^2(\domain)} + k^2 \N{u}^2_{L^2(\domain)} \right)
+aL_I\N{\nT (\gamma u)}^2_{\LtGI} + L_I k^2\N{\gamma u}^2_{\LtGI}
\leq C_2 \N{f}^2_{L^2(\domain)}+ \widetilde{C_2}\N{g_I}^2_\LtGI,
\eeqs
for all $k>0$, where
\beqs
C_2 =2 \left[ \frac{2}{\mu_3}(1+2 A_{\max})^2\left(L_I^2  + \left(\beta+\frac{d}{2k}\right)^2\right)+ \frac{\mu_3}{k^2}\right],\,\,
\widetilde{C_2}= \left[ (1+2A_{\max})\widetilde{C_1}+ \frac{4\mu_3^2}{L_I k^2}\right],
\eeqs
where $\widetilde{C_1}$ is given in \eqref{eq:C1tilde} and $\beta$ is given by \eqref{eq:beta}.

\

\noi (iv) 
 If $\domain_-, \widetilde{\domain}, A, n, \vartheta, f, g_D,$ and $g_I$ satisfy the requirements in the definition of the TEDP (Definition \ref{def:TEDP}), along with the requirements in Condition \ref{cond:3balt}, then the solution of the TEDP exists and is unique. Let $L_I:= \min_{\bx\in \Gamma_I}|\bx|$ and let $a_I L_I$ be the radius of the ball with respect to which $\widetilde{\domain}$ is star-shaped.
Let $L_D:= \min_{\bx\in \Gamma_D}|\bx|$ and let $a_D L_D$ be the radius of the ball with respect to which $\domain_-$ is star-shaped. 
Then
\begin{align*}\nonumber
&\mu_4 \left(\N{\gu}^2_{L^2(\domain)} + k^2 \N{u}^2_{L^2(\domain)} \right)
+a_I L_I\N{\nT (\gamma u)}^2_{\LtGI} + L_I k^2\N{\gamma u}_{\LtGI}+\frac{a_D L_D}{2} \N{\dudnu}^2_\LtGD\\
&\leq C_3 \N{f}^2_{L^2(\domain)}+ \widetilde{C_3}\N{g_I}^2_\LtGI
+ C_4 \N{\nabla_{\Gamma_D}g_D}^2_{L^2(\Gamma_D)} +C_5 k^2 \N{g_D}^2_{L^2(\Gamma_D)},
\end{align*}
for all $k>0$, where 
\begin{align*}
C_3 =2 \left[ \frac{2}{\mu_4}\left(1+\frac{3}{2}n_{\max}\right)^2\left(L_I^2  + \left(\beta+\frac{d-2}{2k}\right)^2\right)+ \frac{\mu_4}{2k^2\varmax}\right],\,
\widetilde{C_3}= \left[ \left(1+\frac{3}{2}n_{\max}\right)\widetilde{C_1}+ \frac{2\mu_4^2}{L_I k^2}\right],
\end{align*}
\beqs
C_4:= 2\left(1 + \frac{3}{2}\varmax\right) L_D\left(1 + \frac{4}{a_D}\right),
\tand\,
C_5:= 2 \left[
\left(1 + \frac{3}{2}\varmax\right) \frac{4}{a_D L_D} \left(\beta+ \frac{d-2}{2k}\right)^2 + \frac{2 \mu_4^2}{a_D L_D}\right].
\eeqs
\noi Observe that, given $k_0>0$, each of the $C_j$ and $\widetilde{C}_j$,  $j=1,\ldots, 5$, can be bounded above, independently of $k$, for $k\geq k_0$.
\end{theorem} 
 
\noi Analogues of Remarks \ref{rem:CWM}, \ref{rem:H2}, \ref{rem:pert},  and \ref{rem:trans} and Corollary \ref{cor:H1} hold for the bounds in Theorem \ref{thm:TEDP}. 
Furthermore, just as Condition \ref{cond:1} with Lipschitz coefficients was relaxed to Condition \ref{cond:1altL} for $L^\infty$ coefficients in the case of the EDP, Condition \ref{cond:1alt} can be relaxed to allow $L^\infty$ coefficients, and an analogue of Theorem \ref{thm:EDP2} proved for the TEDP.

\subsection{Discussion of previous work on bounds for the TEDP with variable $A$ and $n$}

We first discuss results for $d=2,3$ in the case that $A$ and $n$ are continuous.
The TEDP and IIP with scalar Lipschitz $A$ and Lipschitz $n$ were considered in \cite{BrGaPe:15}. In the case $A=1$, a resolvent estimate was proved using Morawetz identities with vector field $\bx$ under the condition \eqref{eq:n3a}. In the case $A \neq 1$, the conditions on $A$ and $n$ in \cite{BrGaPe:15} are more restrictive than \eqref{eq:A1a} (this is because \cite{BrGaPe:15} do not commute $\cL_{A,n}$ with $\bx\cdot \nabla$, but rearrange equation \eqref{eq:1} into an equation for $\Delta u$ and commute this equation with $\bx\cdot \nabla$). In \cite{GaMo:17} a new variational formulation of the IIP with $A\equiv I$ and $n$ Lipschitz was obtained using Morawetz identites with the vector field $\bx$, generalising the formulation for $A\equiv I$ and $n\equiv1$ in \cite{MoSp:14}. This formulation is coercive under conditions essentially equivalent to \eqref{eq:n3a}, and this coercivity then implies that a resolvent estimate holds.

Turning to results where one of $A$ and $n$ is not continuous, a resolvent estimate for the TEDP where $A\equiv I$, $n$ is piecewise constant with jumps on $C^\infty$ convex interfaces with strictly positive curvature, and the problem is nontrapping was proved in \cite{CaVo:10} (by adapting the results of \cite{CaPoVo:99}). 
Morawetz identities with the vector field $\bx$ were used to prove resolvent estimates for the IIP where (i) $d=2$, $A$ is Lipschitz, $n$ is complex and both $A$ and $n$ have a common (nontrapping) jump  \cite{OhVe:18}, (ii) $d=2$, $A\equiv I$, and $n$ is piecewise constant with nontrapping jumps in \cite{Ch:15, BaChGo:17}, and (iii) $A\equiv I$ and $n=1+\eta$ is a random variable with $\|\eta\|_{L^\infty}\leq C/k$ (with $C$ independent of $k$) almost surely \cite{FeLiLo:15}; this last result is essentially the random-variable analogue of the bound discussed in Remark \ref{rem:pert} with $n_0=1$.

Resolvent estimates for the Helmholtz equation posed on a 1-d bounded interval with at least one of the endpoints having an impedance boundary condition were obtained for $A\equiv I$ and $n\in C^1$  in 
\cite{AzKeSt:88}, for $A\equiv I$ and piecewise-constant $n$ in \cite{Ch:15, Ch:16, SaTo:17}, and for piecewise-$C^1$ scalar $A$ 
and piecewise-$C^1$ $n$ in \cite{GrSa:18}. 
The results in \cite{AzKeSt:88, Ch:15, GrSa:18} were proved using variants of Morawetz identities (see the discussion in \cite[\S5]{GrSa:18}) and the results in \cite{SaTo:17} were proved using the explicit expression for the Green's function in 1-d.

\subsection{TEDP-analogues of the results in \S\ref{sec:3} and \S\ref{sec:4}}

\ble[Variational formulation of TEDP with $g_D\equiv 0$]
With $\domain_-, n, A, f,$ and $g_I$ as in Definition \ref{def:TEDP}, let 
\beqs
H_{0,D}^1(\domain):= \big\{ v\in H^1(\domain) : \gamma v=0 \ton \Gamma_D\big\}.
\eeqs
The variational formulation of the TEDP of Definition \ref{def:TEDP} with $g_D= 0$ is: 
\beq\label{eq:TEDPvar}
\text{ find } u \in H^1_{0,D}(\domain) \tst \quad a(u,v)=F(v) \quad \tfa v\in H^1_{0,D}(\domain),
\eeq
where
\beq\label{eq:TEDPa}
a(u,v):= \int_\domain 
\Big((A \gu)\cdot\gvb
- k^2 n u\vb\Big) - \ri k \int_{\Gamma_I}\vartheta\, \gamma u\, \overline{\gamma v} \quad\tand\quad
F(v):= \int_\domain f\, \vb + \int_{\Gamma_I} g_I \, \overline{\gamma v}.
\eeq
\ele

\bre[Continuity of $a(\cdot,\cdot)$]
The sesquilinear form $a(\cdot,\cdot)$ of the TEDP defined by \eqref{eq:TEDPa} is continuous on $H^1(\domain)$; this follows in a similar way to the EDP, except using the 
 multiplicative trace inequality
$\N{\gamma u}^2_{L^2(\partial \domain)} \leq C\N{u}_{L^2(\domain)}\N{u}_{H^1(\domain)}$
\cite[Theorem 1.5.1.10, last formula on Page 41]{Gr:85} to deal with the term on $\Gamma_I$.
\ere

\bre[The Ne\v{c}as regularity result]
The natural analgoue of the Ne\v{c}as regularity result Corollary \ref{cor:Necas} holds for the TEDP. The only difference in the proof is that, since there are different types of boundary conditions on $\Gamma_D$ and $\Gamma_I$, we introduce a $C^\infty$ boundary, $\Gamma_*$, between $\Gamma_D$ and $\Gamma_I$ and apply the Ne\v{c}as result first between $\Gamma_*$ and $\Gamma_D$, and then between $\Gamma_*$ and $\Gamma_R$ (using interior $H^2$-regularity of the operator
$\opL_{A,n}$ 
and the trace theorem \cite[Theorem 3.38]{Mc:00} to get that $\gamma u \in H^1(\Gamma_*)$ and $\dudnuAw \in L^2(\Gamma_*)$).
\ere

\begin{lemma}\textbf{\emph{(Bounding the $H^1$ semi-norm of $u$ via the $L^2$ norm and $L^2$ norm of $f$)}}
Assume there exists a solution to the TEDP of Definition \ref{def:TEDP}. 
Then, for all $k>0$,
\begin{align}\nonumber
&A_{\min}\N{\gu}_{L^2(\domain)}^2 \leq\frac{3}{2} k^2 \varmax\N{u}^2_{L^2(\domain)} + \frac{1}{2k^{2}\varmax}\N{f}^2_{L^2(\domain)}+ \N{\gamma u}_{L^2(\Gamma_D)}\N{\dudnuA}_{L^2(\Gamma_D)}\\
&\hspace{6cm}+\N{\gamma u}_{L^2(\Gamma_I)}\N{g_I}_{L^2(\Gamma_I)},\label{eq:GreenTEDP1}
\end{align}
and 
\begin{align}\nonumber
&\half k^2 \varmin\N{u}^2_{L^2(\domain)} \leq A_{\max}\N{\gu}_{L^2(\domain)}^2 +\frac{1}{2k^{2}\varmin}\N{f}^2_{L^2(\domain)}+ \N{\gamma u}_{L^2(\Gamma_D)}\N{\dudnuA}_{L^2(\Gamma_D)}\\
&\hspace{6cm}+\N{\gamma u}_{L^2(\Gamma_I)}\N{g_I}_{L^2(\Gamma_I)}\label{eq:GreenTEDP2}.
\end{align}
\ele

\bpf[Sketch proof] Applying Green's identity \eqref{eq:Green} with $\domaingen$ the domain in the Definition of the TEDP, $u$ the solution of the TEDP, and $v=u$, we find
\beqs
-\left\langle \dudnuA, \gamma u\right\rangle_{\Gamma_D} + \big\langle \ri k \vartheta\gamma u + g_I, \gamma u\big\rangle_{\Gamma_I} = \int_\domain
(A \gu)\cdot\overline{\gu} - k^2 n \nus - \ub f
\eeqs
The results \eqref{eq:GreenTEDP1} and \eqref{eq:GreenTEDP2} are obtained by taking the real part of this identity, using \eqref{eq:nlimitsEDP}, the Cauchy-Schwarz inequality, \eqref{eq:AellEDP}, and then the Cauchy inequality (on the term involving both $u$ and $f$).
\epf

\subsection{The boundary terms in the integrated Morawetz identity \eqref{eq:morid1int} under an impedance condition.}

The inequality in the following lemma can be seen as an analogue of the inequalities in Lemmas \ref{lem:2.1} and \ref{lem:2.1mod} above. Indeed, the inequalities in Lemmas \ref{lem:2.1} and \ref{lem:2.1mod} are used to deal with the contribution from $\Gamma_R:=\partial B_R$ when bounding the solution of the EDP in $\domain_R := \domain_+ \cap B_R$. The inequality \eqref{eq:impineq} below is used to deal with the contribution from $\Gamma_I$ (the boundary with the impedance condition) when bounding the solution of the TEDP in $\domaingen$.

\begin{lemma}[Inequality on $\Gamma_I$ used to deal with the impedance boundary condition]\label{lem:impineq}
Let the domains $\domain_-, \widetilde{\domain},$ and $\domaingen$, and the functions $\vartheta$ and $g_I$ be as in the definition of the TEDP (Definition \ref{def:TEDP}). Let $v \in V(\domaingen)$ satisfy the impedance boundary condition $\dvdnuw- \ri k \vartheta v =g_I$ on $\Gamma_I$. Let $L_I:= \max_{\bx\in \Gamma_I}|\bx|$ and assume that $\widetilde{\domain}$ is star-shaped with respect to the ball of radius $aL_I$. Let $\alpha\in \Rea$ be arbitrary and let $\beta$ be given by \eqref{eq:beta}.
Then
\begin{align}\nonumber
&\int_{\Gamma_I} (\bx\cdot\bnu)\left(\left|\dvdnu\right|^2 - |\nT (\gamma v)|^2 + k^2 n|\gamma v|^2\right) + 2 \Re \left( \big(\bx\cdot\overline{\nT (\gamma v)} + \ri k \beta \overline{\gamma v} + \alpha\overline{\gamma v} \big) \dvdnu\right) \\
&\quad\leq -\frac{a L_I}{2} \N{\nT (\gamma v)}^2_{\LtGI} -k^2 L_I \N{\gamma v}^2_\LtGI + 
\left( 2 \left(1+ \frac{2}{a}\right) + \frac{\beta}{\vartheta_{\min} L_I}+ \alpha^2\right)
L_I 
\N{g}^2_{\LtGI}. 
\label{eq:impineq}
\end{align}
\end{lemma}

\bpf
The impedance boundary condition implies that
\begin{align*}
I&:= 
\int_{\Gamma_I} (\bx\cdot\bnu)\left(\left|\dvdnu\right|^2 - |\nT (\gamma v)|^2 + k^2 n|\gamma v|^2\right) + 2 \Re \left( \big(\bx\cdot\overline{\nT (\gamma v)} + \ri k \beta \overline{\gamma v} + \alpha\overline{\gamma v} \big) \dvdnu\right) \\
&=\int_{\Gamma_I}(\bx\cdot\bnu) \left( |\ri k \vartheta \gamma v + g_I|^2 - |\nT (\gamma v)|^2 + k^2 n|\gamma v|^2 \right) + 2\Re\big(\bx\cdot \overline{\nT (\gamma v)} \,\ri k \vartheta \gamma v\big) \\
&\qquad\qquad-2 k^2 \beta \vartheta |\gamma v|^2 + 2 \Re\left(\big(\bx\cdot \overline{\nT (\gamma v)} + \ri k\beta \overline{\gamma v} + \alpha \overline{\gamma v}\big) g_I\right).
\end{align*}
Using the Cauchy-Schwarz and Cauchy inequalities, as well as the inequalities $aL_I \leq \bx\cdot\bnu(\bx)\leq L_I$ and \eqref{eq:thetaineq}, we find that
\begin{align*}
I&\leq 2L_I \left( (\vartheta_{\max})^2 k^2 \N{\gamma v}^2_\LtGI + \N{g_I}^2_{\LtGI} \right) - a L_I \N{\nT (\gamma v)}^2_{\LtGI} + k^2 n_{\max,\Gamma_I}L_I \N{\gamma v}^2_{\LtGI}\\
& +L_I \left(\eps_1 \N{\nT (\gamma v)}^2_{\LtGI} + \frac{k^2 (\vartheta_{\max})^2}{\eps_1}\N{\gamma v}^2_{\LtGI} \right) - 2 k^2 \beta \vartheta_{\min} \N{\gamma v}^2_\LtGI \\
&+ L_I\left( \eps_2 \N{\nT (\gamma v)}^2 + \frac{1}{\eps_2} \N{g_I}^2_\LtGI \right) + \big(\eps_3 k^2 \beta^2+ \eps_4\big) \N{\gamma v}^2_\LtGI + \left(\frac{1}{\eps_3} +\frac{\alpha^2}{\eps_4}\right)\N{g_I}^2_\LtGI,
\end{align*}
for any $\eps_j>0$, $j=1,2,3,4$. Choosing $\eps_1 =\eps_2=a/4$, we have
\begin{align*}
I\leq &- k^2 L_I\left( \frac{2\beta\vartheta_{\min}}{L_I} -n_{\max,\Gamma_I} - 2(\vartheta_{\max})^2 \left(1+ \frac{2}{a}\right) -\frac{\eps_3 \beta^2}{L_I} - \frac{\eps_4}{L_I k^2} \right) \N{\gamma v}^2_\LtGI \\
&\quad- \frac{a L_I}{2}\N{\nT (\gamma v)}^2_\LtGI + \left( 2 \left(1+ \frac{2}{a}\right) + \frac{1}{\eps_3 L_I}+ \frac{\alpha^2}{\eps_4 L_I}\right)L_I\N{g}^2_{\LtGI}.
\end{align*}
Our goal is to choose $\beta$ large enough so that the bracket in front of $\|\gamma v\|^2_{\LtGI}$ on the right-hand side is positive (in fact, we'll choose $\beta$ so that the bracket equals one), but this is only possible if $\eps_3$ is chosen appropriately. Letting $\eps_3= \vartheta_{\min}/\beta$ and $\eps_4 = 1/L_I$ we have
\begin{align*}
I\leq & - k^2 L_I\left( \frac{\beta\vartheta_{\min}}{L_I} -n_{\max,\Gamma_I} - 2(\vartheta_{\max})^2 \left(1+ \frac{2}{a}\right) - \frac{1}{(kL_I)^2} \right) \N{\gamma v}^2_\LtGI \\
&\quad- \frac{a L_I}{2}\N{\nT (\gamma v)}^2_\LtGI+ \left( 2 \left(1+ \frac{2}{a}\right) + \frac{\beta}{\vartheta_{\min}L_I}+\alpha^2\right)
L_I \N{g}^2_{\LtGI}.
\end{align*}
We then choose $\beta$ as in \eqref{eq:beta} and the result \eqref{eq:impineq} follows.
\epf

\subsection{Outline of the proof of Theorem \ref{thm:TEDP}}

This follows the proof of Theorems \ref{thm:EDP1} and \ref{thm:EDP3} very closely; the main difference that we now use Lemma \ref{lem:impineq} to deal with the terms arising from the impedance boundary $\Gamma_I$ instead of Lemmas \ref{lem:2.1} and \ref{lem:2.1mod} to deal with the terms arising from the boundary $\Gamma_R$.
Applying the Morawetz identity in $\Omega$ is justified in Part (iii) by the density result of Lemma \ref{lem:density}
(as in Part (ii) of Theorem \ref{thm:EDP3}). For Parts (i) and (ii), we first prove the bounds for $C^{1,1}$ star-shaped $\Omega_-$, and then use identical arguments to those in the proof of Theorem \ref{thm:EDP1} to extend these bounds to Lipschitz star-shaped $\Oi$. When proving the bounds for $C^{1,1}$ star-shaped $\Omega_-$, using the Morawetz identity in $\Omega$ is justifying by (a) introducing a $C^\infty$ boundary, $\Gamma_*$, between $\Gamma_D$ and $\Gamma_I$ such that $A\equiv I$ both in a neighbourhood of $\Gamma_*$ and between $\Gamma_*$ and $\Gamma_I$, (b) using $H^2$-regularity to justify applying the identity between $\Gamma_D$ and $\Gamma_*$, and (c) using Lemma \ref{lem:density} to justify applying the identity between $\Gamma_*$ and $\Gamma_I$. 

For (i), we choose $\beta$ as in \eqref{eq:beta} and let $2\alpha=d-1$.
For (ii), we choose $\beta$ as in \eqref{eq:beta} and let  $2\alpha=d$. We use \eqref{eq:GreenTEDP2} to introduce $k^2 \N{u}^2_{\LtD}$ (just as we used \eqref{eq:GreenEDP2} to introduce $k^2 \N{u}^2_{L^2(\Omega_R)}$ in the proof of Part (i) of Theorem \ref{thm:EDP3}).
For (iv), we choose $\beta$ as in \eqref{eq:beta} and let  $2\alpha=d-2$. We use \eqref{eq:GreenTEDP1} to introduce $\N{\gu}^2_{\LtD}$ (just as we used \eqref{eq:GreenEDP1} to introduce $\N{\gu}^2_{L^2(\Omega_R)}$ in the proof of Part (ii) of Theorem \ref{thm:EDP3}). The treatment of the terms on $\Gamma_D$ is essentially identical 
to that in Part (ii) of Theorem \ref{thm:EDP3}.

\section*{Acknowledgements}
For useful discussions, the authors thank Yves Capdeboscq (University of Oxford), Suresh 
Eswarathasan (Cardiff University), Andrea Moiola (Universit\`a di Pavia), Luca Rondi (Universit\`a di Trieste), and particularly Jared Wunsch (Northwestern University). We thank the Bath Institute for Mathematical Innovation for funding Jared's visit to Bath in March 2016.
We also thank the referee for their constructive comments that improved the organisation of the paper.

ORP is supported by a PhD studentship from the EPSRC Centre for Doctoral Training in Statistical Applied Mathematics at Bath (SAMBa), under the grant EP/L015684/1.
EAS is supported by EPSRC grant EP/R005591/1.

\end{appendix}

\footnotesize{


}

\end{document}